\documentclass[10pt,reqno]{amsart} 

\usepackage{appendix}
\usepackage{float} 
\usepackage{soul}
\usepackage{hhline}
\usepackage{yhmath}\usepackage{epsfig,amssymb,amsmath,version}
\usepackage{amssymb,version,graphicx,fancybox,mathrsfs}
\usepackage{subcaption}
\usepackage{enumitem}
\usepackage{bbm}
\usepackage{url,hyperref}
\usepackage{color}
\usepackage{stmaryrd}
\usepackage{multirow}
\usepackage{booktabs,siunitx}
\usepackage{multicol,epstopdf}
\usepackage{xypic}
\usepackage[]{graphicx}
\usepackage[all]{xy}
\usepackage{tikz,ifthen}
\usetikzlibrary{positioning, math, decorations.markings, arrows.meta}
\usetikzlibrary{calc,shapes,cd}
\usetikzlibrary{knots}
\usepgfmodule{decorations}
\allowdisplaybreaks
\tikzset{knotarrow/.pic={ \draw[edge, <-] (0,0) -- +(-.001,0);}}
\tikzset{edge/.style={line width=0.8}}
\tikzset{wall/.style={very thick}}
\tikzset{->-/.style n args={2}{decoration={markings, mark=at position #1 with {\arrow{#2}}}, postaction={decorate}}}

\allowdisplaybreaks[4]

\newtheorem{theorem}{Theorem}[section]
\newtheorem{lemma}[theorem]{Lemma}
\newtheorem{definition}[theorem]{Definition}
\newtheorem{corollary}[theorem]{Corollary}
\newtheorem{proposition}[theorem]{Proposition}
\newtheorem{remark}[theorem]{Remark}

\newtheorem{example}[theorem]{Example}
\newtheorem{convention}[theorem]{Convention}
\newtheorem{maintheorem}{Theorem}

\newtheorem{maincorollary}[maintheorem]{Corollary}

\graphicspath{{./}{./Figures/}}

\def\fS{\mathfrak{S}}
\def\cS{\mathscr S}
\def\SL{{\rm SL}_n}
\def\bT{\mathbb T}
\def\bP{\mathbb P}
\def\cY{\mathcal Y}
\def\Exc{\mathsf{S}}

\newcommand{\Z}{\mathbb{Z}}
\newcommand{\C}{\mathbb{C}}
\newcommand{\Q}{\mathbb{Q}}
\newcommand{\cC}{\mathcal{C}}
\newcommand{\cK}{\mathcal{K}}
\newcommand{\cQ}{\mathcal{Q}}
\def\cH{\mathcal{H}}
\newcommand{\cG}{\mathcal{G}}
\def\cD{\mathcal{D}}
\newcommand{\slh}{\widehat{\mathfrak{sl}}_n}
\newcommand{\Vc}{V_{\mathrm{core}}}
\newcommand{\Viso}{V_{\mathrm{iso}}}
\newcommand{\Lst}{V_{\mathrm{first}}}
\def\qc{{\rm qc}}
\newcommand{\core}{{\rm core}}

\begin{document}
\bibliographystyle{alpha}

\title{Categorification of the reduced stated ${\rm SL}_n$-skein algebras of polygons}

\author[Min Huang]{Min Huang}
\address{Min Huang, School of Mathematics (Zhuhai), Sun Yat-sen University, Zhuhai, China}
\email{huangm97@mail.sysu.edu.cn}

\author[Zhihao Wang]{Zhihao Wang}
\address{Zhihao Wang, School of Mathematics, Korea Institute for Advanced Study (KIAS), 85 Hoegi-ro, Dongdaemun-gu, Seoul 02455, Republic of Korea}
\email{zhihaowang@kias.re.kr}


\begin{abstract}
Let $\overline{\mathscr{S}}(\mathbb P_m)$ be the reduced stated ${\rm SL}_n$-skein algebra of a polygon with $m\geq 3$ boundary punctures. Building on the quantum cluster algebra structure of $\overline{\mathscr{S}}(\mathbb P_m)$ established by Cao and the authors~\cite{CHW}, we provide a monoidal categorification of $\overline{\mathscr{S}}(\mathbb P_m)$. Our construction utilizes the monoidal category $\mathcal{C}(\Delta^{m-2})$ of finite-dimensional modules over the quantum affine algebra $U_q'(\widehat{\mathfrak{sl}}_n)$, where $\Delta$ denotes the Garside element of ${\rm Br}_n$. This category arises in the monoidal categorification of cluster algebras, following the works of Hernandez-Leclerc \cite{HL10,HL15}, Qin \cite{Qin24a}, and Kashiwara, Kim, Oh, and Park~\cite{KKOP3}. As a consequence, $\overline{\mathscr{S}}(\mathbb P_m)$ is shown to possess a basis with positive structure constants.

A key role is played by the subalgebra $\mathcal A_{n,m}$ spanned by the stated $n$-webs whose $m-1$ edge weights are zero. We show that $\mathcal A_{n,m}$ is a quantum cluster algebra whose seed coincides with the Kashiwara--Kim--Oh--Park seed of $\Delta^{m-2}$. Hence it is the quantum Grothendieck ring of a monoidal localization, at the frozen modules, of a subcategory of $\mathcal{C}(\Delta^{m-2})$ with the same simple objects, and a quantization of the braid variety $X(\Delta^{m-1})$ \cite{CGGLSS}. The whole algebra $\overline{\mathscr{S}}(\mathbb P_m)$ is a skew Laurent extension of $\mathcal A_{n,m}$ in $(m-1)(n-1)$ frozen variables.
\end{abstract}

\maketitle

\tableofcontents

\section{Introduction}\label{sec-intro}

Throughout the paper, we fix an integer $n\ge2$ and the ground ring $R=\Z[\omega^{\pm1/2}]$, and we let $\xi=\omega^n$ and $t=\omega^{n^2}=\xi^n$ (see \eqref{eq-constants}). Unless otherwise specified, all algebras are $R$-algebras.

\subsection{Background}

\subsubsection{Reduced stated $\SL$-skein algebras and quantum cluster algebras}

Skein algebras of surfaces were introduced by Przytycki and Turaev~\cite{Prz91,Tur91} as quantizations of character varieties. Bonahon and Wong embedded the Kauffman bracket skein algebra into quantum tori attached to ideal triangulations by their quantum trace~\cite{BW11}, and Muller showed that the skein algebras of marked surfaces are quantum cluster algebras~\cite{Mul16}; for $\mathrm{SL}_3$ and $\SL$, see~\cite{IY23,HW}. Stated skein algebras of punctured bordered surfaces were introduced by L\^e~\cite{Le18} for $\mathrm{SL}_2$ (see also Costantino--L\^e~\cite{CL22}), and extended to $\SL$ by L\^e--Sikora~\cite{LS}. They are spanned by stated $n$-webs, i.e.\ $n$-webs in the thickened surface whose endpoints carry states in $\{1,\dots,n\}$, modulo the relations \eqref{w.cross}--\eqref{wzh.eight}. The \emph{reduced} stated $\SL$-skein algebra $\overline{\cS}(\fS)$ is the quotient by the ideal generated by the bad arcs (\S\ref{sub-def-reduced-skein}). L\^e and Yu~\cite{LY} constructed quantum trace maps from reduced stated skein algebras to quantum torus attached to $n$-triangulations.

Quantum cluster algebras were introduced by Berenstein and
Zelevinsky~\cite{BZ} as quantizations of the cluster algebras of Fomin and Zelevinsky~\cite{FZ02}. A quantum seed $\mathsf s=(Q,\Pi,M)$ consists of an exchange matrix $Q$, a skew-symmetric matrix $\Pi$, and a toric frame $M$ (Definition~\ref{def-quantum-seed}). Mutations produce the cluster variables, and the quantum cluster algebra $\mathscr A(\mathsf s)$ is the algebra generated by
all cluster variables together with the frozen variables (and possibly their inverses) (Definition~\ref{def-quan-cluster-algebra}).

The reduced stated skein algebra admits a quantum cluster algebra structure for
polygons \cite{CHW}. Let $\bP_m$ be the closed disk with $m\ge3$ punctures on
its boundary, and fix the star-like triangulation of $\bP_m$ at the puncture
$1$. Each triangle carries an $n$-triangulation in the sense of Fock--Goncharov~\cite{FG06}, and its small vertices form the vertex set $V$ of
a quiver $Q$ (Convention~\ref{conv:CHWlabels} and Definition~\ref{def:CHWquiver}). For each small vertex $x$ there is an explicit stated $n$-web, whose normalization is denoted by $A_x$ (Figure~\ref{fig:placeholder}). Cao and the authors showed that the $A_x$ form the initial cluster of a quantum seed $\mathsf s_{\rm skein}=(Q,\Lambda_{\rm skein},M_{\rm skein})$, which we call the \emph{skein seed}, and that
\begin{align}\label{eq-intro-CHW}
\overline{\cS}(\bP_m)=\mathscr A(\mathsf s_{\rm skein})=\mathscr U(\mathsf
s_{\rm skein})
\end{align}
(Proposition~\ref{prop:CHWseed} and Theorem~\ref{thm:CHW}; see~\cite{CHW}). The
seed $\mathsf s_{\rm skein}$ has $m(n-1)$ frozen vertices, namely the small vertices on the boundary edges of $\bP_m$. We use the skein relations of L\^e--Sikora~\cite{LS}, whose parameter $t$ is inverse to that of~\cite{CHW} (\S\ref{sub-def-reduced-skein}); with this convention, all comparisons below hold with no sign.

\subsubsection{Quantum Grothendieck rings of braid-word categories}\label{subsec-quantum-Gro}

Hernandez and Leclerc introduced monoidal categorification of cluster algebras, and realized cluster algebras as Grothendieck rings of categories of finite-dimensional modules over quantum affine algebras~\cite{HL10,HL15,HL16}. Quantum versions use the $(q,t)$-characters of Nakajima, Varagnolo--Vasserot and Hernandez~\cite{Nakajima,VV,Hernandez}, and were developed in~\cite{Qin17,Bit21,KKKO18,KKOP2}.

Qin \cite{Qin24a} and Kashiwara, Kim, Oh, and Park \cite{KKOP3} independently attached to each positive braid $\mathsf b\in\mathrm{Br}_n^+$ a monoidal subcategory $\cC(\mathsf b)$ of the Hernandez--Leclerc category of $U_q'(\slh)$-modules (\S\ref{subsec_braid_cat}). Its quantum Grothendieck ring $\cK_t(\cC(\mathsf b))$ is the $\Z[t^{\pm1/2}]$-algebra generated by the $(q,t)$-characters of its simple objects (Definition~\ref{def:Kt}). For a word $\mathbf i$ of $\mathsf b$ in which every colour occurs, they defined a quantum seed inside $\cK_t(\cC(\mathsf b))$, which we denote by $\mathsf s(\mathbf i)$ and call a \emph{KKOP seed} (Definition~\ref{def:KKOPseed}). Its vertices are the positions of $\mathbf i$, its frozen vertices are the $n-1$ last occurrences of the colours, its exchange matrix is the matrix $Q_{\mathbf i}$ of the quiver of the word (Definition~\ref{def:BFZ}), and its initial cluster variables are the $(q,t)$-characters of the affine determinantal modules $\mathsf M_s$ defined in~\eqref{eq-Ms}. Its quantum form is given by $n$ times the matrix of $R$-matrix invariants $\Lambda(\mathsf M_s,\mathsf M_{s'})$ (see Definition~\ref{def:Lambda}); the factor $n$ is introduced to accommodate the variable $\xi$, consistent with the skein side. They proved that $\cK_t(\cC(\mathsf b))$ coincides with the quantum cluster algebra associated with this seed, incorporating non-invertible frozen variables, and that the cluster monomials correspond to the $(q,t)$-characters of real simple modules (Theorem~\ref{thm:KKOP}). Furthermore, Bi utilized the Ore localizations of these rings at the frozen classes to construct a quantization of braid varieties~\cite{Bi}, see also \cite{Qin24a} for the case that $b=\Delta^{m-2}$, where $\Delta$ is the Garside element of $Br_n$.

Theorems~\ref{thm:CHW} and~\ref{thm:KKOP} constitute the two main pillars of this paper. The former identifies the reduced stated skein algebra of a polygon with a quantum cluster algebra, while the latter establishes an analogous correspondence for the quantum Grothendieck ring of a braid-word category. In this work, we compare these two quantum cluster algebra structures and lift this comparison to the level of monoidal categories.

\subsection{A categorification of the reduced stated $\SL$-skein algebra of polygons}

Categorification and positivity are recurring themes in skein theory. Thurston proved that the bracelet basis of the classical (i.e., with ground ring \(\mathbb Z\) and quantum parameter equal to $1$) \(\mathrm{SL}_2\)-skein algebra of a surface has positive structure constants~\cite{Thu14}. Asaeda, Przytycki and Sikora constructed homology theories for links in $I$-bundles over surfaces whose graded Euler characteristics recover the coefficients of a link in the standard basis of the Kauffman bracket skein module~\cite{APS04}; this is a categorification in the sense of Khovanov homology, and should be distinguished from the monoidal categorification considered in this paper. Queffelec and Wedrich later categorified skein modules of surfaces using Khovanov--Rozansky homology~\cite{QW21}. Closest to our setting, Allegretti, Kim and Shan identified the Kauffman bracket skein algebra of a genus zero surface with boundary with the Grothendieck ring of the bounded derived category of equivariant coherent sheaves on a Braverman--Finkelberg--Nakajima variety of triples, equipped with the convolution product, thereby giving a monoidal categorification of these skein algebras~\cite{AKS25}. In the context of cluster algebras, monoidal categorification yields positive bases and identifies cluster monomials with the classes of real simple objects~\cite{HL10,KKKO18,Qin17,KKOP3}.

By \eqref{eq-intro-CHW}, the reduced stated $\SL$-skein algebra of a polygon
carries a quantum cluster algebra structure with invertible frozen variables.
The aim of this paper is to give a monoidal categorification of this algebra.
After deleting the frozen vertices on the boundary edges
$(1,2),(2,3),\dots,(m-1,m)$, the quiver of the skein seed coincides with that
of the KKOP seed of the word $u^{m-2}$ \eqref{eq-intro-u}
(Lemma~\ref{lem:M}). It is therefore natural to expect that
$\cC(\Delta^{m-2})$ gives rise to a monoidal category categorifying
$\overline{\cS}(\bP_m)$. However, three obstructions must be overcome:
\begin{enumerate}
\item the skein seed has $(m-1)(n-1)$ more frozen vertices than the KKOP seed,
and these vertices are connected to mutable vertices by arrows;
\item the quantum forms of the two seeds differ: on the skein side it is
defined by webs, whereas on the representation side it is defined by
$R$-matrices;
\item the frozen variables of the reduced skein algebra are invertible, while
the frozen modules of $\cC(\Delta^{m-2})$ are not.
\end{enumerate}
We address (1) via a quasi-isomorphism of quantum cluster algebras, which isolates the extra frozen vertices (Proposition~\ref{prop:sqc}). Remarkably, under this quasi-isomorphism, the two quantum forms coincide on the remaining vertices (Theorem~\ref{thmA}), thereby resolving (2). We handle (3) by employing a monoidal localization of the category at the frozen modules (Proposition~\ref{prop:loc}). Finally, the extra frozen variables are realized through a grading of the localized category; since these variables are not central in the skein algebra, this necessitates a twist of the multiplication (Theorem~\ref{thmB}).

Let
\begin{align}\label{eq-intro-u}
u=(1\,2\cdots n{-}1)(1\,2\cdots n{-}2)\cdots(1\,2)(1)
\end{align}
be the reduced word of the longest element of the symmetric group $S_n$ of Definition~\ref{def:u}. Then $u^{m-2}$ is a word for $\Delta^{m-2}$, where $\Delta$ is the Garside element of $\mathrm{Br}_n$ (\S\ref{subsec-wordu}). Let $\Viso\subset V$ be the set of the $(m-1)(n-1)$ frozen vertices of $\mathsf s_{\rm skein}$ on the boundary edges $(1,2),(2,3),\dots,(m-1,m)$, and let $\Vc=V\setminus \Viso$ (Definition~\ref{def:core}). In \S\ref{subsec-qc} we use the edge weights of stated webs to construct an explicit quasi-isomorphic seed $\mathsf s^\qc_{\rm skein}$ of $\mathsf s_{\rm skein}$ in which all vertices of $\Viso$ are isolated (Definition~\ref{def:sqc} and Proposition~\ref{prop:sqc}). Its restriction to $\Vc$ is a sub-seed, which we denote by $\mathsf s^\qc_\core$.

\begin{maintheorem}[Theorem~\ref{thm:seedmatch}]\label{thmA}
There is an explicit bijection $\psi$ from $\Vc$ to the set of positions of the word $u^{m-2}$ (Lemma~\ref{lem:M}) which identifies the sub-seed $\mathsf s^\qc_\core$ with the KKOP seed $\mathsf s(u^{m-2})$: the two seeds have the same exchange matrix and the same quantum form. In particular, $\mathscr A(\mathsf s^\qc_\core)\cong\mathscr A(\mathsf s(u^{m-2}))$ as quantum cluster algebras, matching mutations, cluster variables, and cluster monomials.
\end{maintheorem}

Let $F_1,\dots,F_{n-1}$ be the frozen modules of $\cC=\cC(\Delta^{m-2})$. We construct a monoidal localization $\widetilde\cC$ (see~\eqref{eq-C-Delta} and Proposition~\ref{prop:loc}), in which the $F_i$ become invertible, of the full subcategory $\cC^{\mathrm s}\subseteq\cC$ of subquotients of tensor products of cuspidal modules (Definition~\ref{def:Cs}), which has the same simple objects as $\cC$, and its quantum Grothendieck ring $\cK_t(\widetilde\cC)$ (Definition~\ref{def:Ktloc}). Let $\widetilde\cC_{\Viso}$ be the category of $\Z^{\Viso}$-graded objects of $\widetilde\cC$ (\S\ref{subsec-CI}). Its quantum Grothendieck ring is $\cK_t(\widetilde\cC_{\Viso})=\cK_t(\widetilde\cC)[Z_g^{\pm1}:g\in \Viso]$, where the $Z_g$ are central. We regard $\Z[t^{\pm1/2}]$ as a subring of $R$ via $t^{1/2}=\xi^{n/2}$. Finally, let
\begin{align}\label{eq-intro-Vfirst}
\Lst=\{(j,1)^1\mid 1\le j\le n-1\}\subseteq\Vc
\end{align}
be the set of the $n-1$ small vertices of the first triangle $\Delta_1$ adjacent to the edge $(1,2)$ (see~\eqref{eq-Lst}); the bijection $\psi$ of Theorem~\ref{thmA} maps $\Lst$ onto the first occurrences of the colours in $u^{m-2}$, and $\Z^{\Lst}$ is the cokernel of the exchange matrix of $\mathsf s^\qc_\core$ (Lemma~\ref{lem:Lst}).

\begin{maintheorem}[Theorem~\ref{thm:main}]\label{thmB}
There is an explicit skew-symmetric integer matrix $\Xi$ on $\Z^{\Lst}\oplus\Z^{\Viso}$ (Definition~\ref{def:Xi}) and an isomorphism of $R$-algebras
\begin{align}\label{eq-intro-Phi}
\Phi\colon\bigl(\cK_t(\widetilde\cC_{\Viso})\otimes_{\Z[t^{\pm1/2}]}R\bigr)^{\Xi}\xrightarrow{\ \sim\ }\overline{\cS}(\bP_m),
\end{align}
where $(-)^\Xi$ is the cocycle twist \eqref{eq-twist-product}, for a grading in which the $(q,t)$-characters of all simple objects of $\widetilde\cC_{\Viso}$ are homogeneous. Moreover:
\begin{enumerate}[label={\rm (\alph*)}]
\item the cluster monomials of the KKOP seed $\mathsf s(u^{m-2})$ and of its mutations are, up to powers of $t^{1/2}$, $(q,t)$-characters of real simple objects, and $\Phi$ maps them to cluster monomials of $\overline{\cS}(\bP_m)$, up to Laurent monomials in the frozen variables and powers of $\xi^{1/2}$;
\item at $\omega^{1/2}=1$, $\Phi$ induces a ring isomorphism $K(\widetilde\cC_{\Viso})\cong\overline{\cS}(\bP_m)\otimes_{\mathbb Z[\omega^{\pm 1/2}]} \mathbb Z$, where $K(\widetilde\cC_{\Viso})$ is the Grothendieck ring of the monoidal category $\widetilde\cC_{\Viso}$, and $\Z[\omega^{\pm1/2}]$ acts on $\Z$ by $\omega^{1/2}\mapsto1$. So the specialization of the reduced skein algebra at $\omega^{1/2}=1$ is the Grothendieck ring of an abelian monoidal category, in which the cluster monomials are classes of real simple objects and the frozen variables are classes of invertible objects.
\end{enumerate}
\end{maintheorem}

The matrix $\Xi$ vanishes on $\Z^{\Lst}\times\Z^{\Lst}$, and it is the restriction of the quantum form of $\mathsf s^\qc_{\rm skein}$ to the remaining pairs (Definition~\ref{def:Xi}, \eqref{eq-Xi}); the vanishing reflects the equality of the quantum forms on $\Vc$ in Theorem~\ref{thm:seedmatch}. So the twist only concerns the variables $Z_g$, $g\in \Viso$, which are central in $\cK_t(\widetilde\cC_{\Viso})$, while the corresponding frozen variables $A_g$ of the skein algebra are not central. Since all simple $(q,t)$-characters are homogeneous (Proposition~\ref{prop:homog-simple} and Theorem~\ref{thm:locDelta}(b)), the twist only rescales structure constants by powers of $\xi^{1/2}$, and we obtain a positive basis.

\begin{maincorollary}[Corollary~\ref{cor:positivity}]\label{thmC}
The images under $\Phi$ of the $(q,t)$-characters of the simple objects of $\widetilde\cC_{\Viso}$ form an $R$-basis of $\overline{\cS}(\bP_m)$, whose structure constants lie in $\Z_{\ge0}[\xi^{\pm1/2}]$.
\end{maincorollary}

We briefly outline the proof of Theorems~\ref{thmA} and~\ref{thmB}.

\vspace{2mm}

The first step is combinatorial. Delete from the quiver $Q$ the vertices of $\Viso$. We show that what remains, the \emph{core}, is equal to the quiver of the braid word $u^{m-2}$, with no mutation and no sign (Lemma~\ref{lem:M}); the proof compares two explicit lists of arrows (Lemmas~\ref{lem:corequiver} and~\ref{lem:Quk}).

The second step isolates the vertices of $\Viso$. By Definition~\ref{def:edgewt}, every stated web carries an edge weight in the weight lattice $\mathsf P$ of $\mathfrak{sl}_n$ on each boundary edge. As shown in Figure~\ref{fig:placeholder}, the edge weights of the initial webs are determined explicitly (Lemma~\ref{lem:legwt}); these weights are compatible with the exchange matrix (Lemma~\ref{lem:wtkills}), and the weights of the frozen variables $A_g$ ($g\in \Viso$) form a basis of $\mathsf P^{m-1}$ (Lemma~\ref{lem:theta}). This structure yields an admissible change of frozen coordinates (Definition~\ref{def:admissible}) that isolates the vertices in $\Viso$ (Lemma~\ref{lem:balP}). In the sense of~\cite{KWQ,CHL}, this constitutes a quasi-isomorphism---a quantum analogue of the quasi-homomorphisms introduced by Fraser~\cite{Fraser} (Corollary~\ref{cor:Pchange}).

The third step compares the quantum forms. Since both forms are compatible with the same exchange matrix and the same scalar $2n$ (Proposition~\ref{prop:sqc}(c), Theorem~\ref{thm:KKOP}(a)), they coincide as soon as they agree on the $n-1$ vertices of $\Lst$, which span the cokernel of the exchange matrix (Lemma~\ref{lem:Lst}). Agreement on these vertices is verified by direct computation (Lemmas~\ref{lem:XiK} and~\ref{lem:XiS}). This completes the proof of Theorem~\ref{thmA}.

The final step involves the localization of the category. To realize the inverses of frozen variables, we construct braiders for the frozen modules using renormalized universal $R$-matrices. This construction is carried out on the full subcategory $\cC^{\mathrm s}$, consisting of subquotients of tensor products of cuspidal modules (Definition~\ref{def:Cs}, Propositions~\ref{prop:braider} and~\ref{prop:commfamily}). Applying the localization theorem of~\cite{KKOP-loc} yields the abelian monoidal category $\widetilde\cC$ (Proposition~\ref{prop:loc}), whose quantum Grothendieck ring is identified as the Ore localization of $\cK_t(\cC)$ at the frozen classes (Theorem~\ref{thm:Ktloc}). This construction applies to every positive braid in all simply-laced types, provided the defining word contains every color (Remark~\ref{rem:ADE}). Consequently, it recovers the localized quantum Grothendieck rings of Bi~\cite{Bi} as quantum Grothendieck rings of abelian monoidal categories (see the discussion after Theorem~\ref{thmD}). Finally, Theorem~\ref{thmB} follows from Theorem~\ref{thmA}, Theorem~\ref{thm:Ktloc}, and the skew Laurent decomposition established in Proposition~\ref{prop:skewLaurent}.

\subsection{Braid varieties and the reduced stated $\SL$-skein algebra}

The braid variety $X(\beta)$ of a positive braid word $\beta$ parametrizes sequences of flags in relative positions prescribed by $\beta$ (\S\ref{subsec:braidvar}). Casals, Gorsky, Gorsky, Le, Shen and Simental proved that its coordinate ring is a cluster algebra, with seeds given by Demazure weaves~\cite{CGGLSS}; see also~\cite{GLSBS,CGGSSS}. Braid varieties are also augmentation varieties of Legendrian links~\cite{Kal06,CN22,CGGS}.

It is natural to ask whether the skein algebra contains a quantization of the braid variety $X(\Delta^{m-1})$. We show that it does, and that this subalgebra has a simple description in terms of webs. The skein algebra is graded by edge weights (Lemma~\ref{lem:edgegrading}). Let
\begin{align}\label{eq-intro-Anm}
\mathcal A_{n,m}=\operatorname{span}_R\bigl\{\alpha\ \text{stated $n$-web}\mid\ \operatorname{wt}_e(\alpha)=0\ \text{for } e=(1,2),(2,3),\dots,(m-1,m)\bigr\}\subseteq\overline{\cS}(\bP_m)
\end{align}
be the span of the webs that are balanced on these edges (Definition~\ref{def:Anm}). Classically, it is the algebra of invariants of the torus $H^{m-1}$ of decorations, which rescales the states on these edges (Remark~\ref{rem:TI} and Proposition~\ref{prop:Aquant-more}(c)).

\begin{maintheorem}[Theorem~\ref{thm:Aquant}]\label{thmD}
\leavevmode
\begin{enumerate}[label={\rm (\alph*)}]
\item $\overline{\cS}(\bP_m)=\bigoplus_{\mathbf c\in\Z^{\Viso}}\mathcal A_{n,m}\,A^{\mathbf c}$ is a skew Laurent extension of $\mathcal A_{n,m}$ in the frozen variables $A_g$, $g\in \Viso$.
\item $\mathcal A_{n,m}=\mathscr A(\mathsf s^\qc_\core)=\mathscr U(\mathsf s^\qc_\core)$. In particular, $\mathcal A_{n,m}$ is a quantum cluster algebra whose seed coincides with the KKOP seed $\mathsf s(u^{m-2})$.
\item There is an isomorphism of $R$-algebras $\cK_t(\widetilde\cC)\otimes_{\Z[t^{\pm1/2}]}R\cong\mathcal A_{n,m}$, with no twist, which matches cluster monomials. In particular, $\mathcal A_{n,m}$ is a free $R$-module, and the simple objects of $\widetilde\cC$ give a basis of $\mathcal A_{n,m}$ with structure constants in $\Z_{\ge0}[t^{\pm1/2}]$.
\item $\mathcal A_{n,m}\otimes_R\C\cong\C[X(\Delta^{m-1})]$ at $\omega^{1/2}=1$, compatibly with the cluster structure of~\cite{CGGLSS}.
\end{enumerate}
\end{maintheorem}

Hence $\mathcal A_{n,m}$ is a cluster quantization of $\C[X(\Delta^{m-1})]$ (Corollary~\ref{cor:Aquantization}).
Qin \cite{Qin24a} and Bi \cite{Bi} quantized $X(\Delta^{m-1})$ by localizing the quantum Grothendieck ring of a KKOP category at the frozen classes. This localization is algebraic: the inverses of the frozen classes are adjoined formally and do not correspond to objects. By Theorem~\ref{thmD}(c) and~\cite[Lemma~8.14]{Bi}, together with the fact that the quantum Grothendieck ring of $\cC(\Delta^{m-2})$ and its initial seed do not depend on the choice of $Q$-datum~\cite[Rem.~2.9, Thms.~9.7, 9.12]{KKOP3}, Their quantization is isomorphic to $\mathcal A_{n,m}$ after extending scalars to $R$. We realize this localization categorically: it is the quantum Grothendieck ring of the monoidal localization $\widetilde\cC$, in which the frozen modules become invertible objects (Theorem~\ref{thm:Ktloc}). 
In addition, our approach realizes the quantization inside the skein algebra by balanced webs.

Theorem~\ref{thmD}(a) and (d) give the following.

\begin{maincorollary}[Corollary~\ref{cor:braidvar}]\label{thmE}
The classical reduced skein algebra is the coordinate ring of a braid variety times a torus:
\begin{align}\label{eq-intro-braidvar}
\overline{\cS}(\bP_m)\otimes_{\mathbb Z[\omega^{\pm 1/2}]} \mathbb C\ \cong\ \C[X(\Delta^{m-1})]\otimes\C[Z_g^{\pm1}:g\in \Viso],
\end{align}
where $\Z[\omega^{\pm1/2}]$ acts on $\C$ by $\omega^{1/2}\mapsto1$, compatibly with the cluster structures of~\cite{CHW} and~\cite{CGGLSS}.
\end{maincorollary}

The dimensions match: $\dim X(\Delta^{m-1})+(m-1)(n-1)$ is the number of small vertices of $\bP_m$, which is the Gelfand--Kirillov dimension of the skein algebra~\cite{LY}. Corollary~\ref{thmE} is the skein-theoretic counterpart of the description of decorated local systems on polygons by Ishibashi--Oya--Shen~\cite[\S6.2.2]{IOS}, and its proof is algebraic.

We briefly outline the proof of Theorem~\ref{thmD}. The skein algebra is graded by edge weights (Lemma~\ref{lem:edgegrading}), and $\mathcal A_{n,m}$ is the weight-zero component with respect to the edges $(1,2),\dots,(m-1,m)$ in the coarser grading~\eqref{eq-coarse-grading}. In the seed $\mathsf s^\qc_{\rm skein}$ from the second step above, every cluster variable at a vertex of $\Vc$ is balanced: this holds by construction in the initial seed (Proposition~\ref{prop:sqc}(d)), and is preserved under mutation because the vertices of $\Viso$ are isolated and the skein algebra is a domain. The frozen variables $A_g$, $g\in\Viso$, have linearly independent edge weights (Lemma~\ref{lem:theta}), so the skein algebra decomposes as $\bigoplus_{\mathbf c}\mathcal A_{n,m}A^{\mathbf c}$. Consequently, $\mathcal A_{n,m}$ is the quantum cluster algebra associated with the sub-seed $\mathsf s^\qc_\core$ (Proposition~\ref{prop:skewLaurent}). This establishes parts (a) and (b). Part (c) follows from Theorem~\ref{thmA} and Theorem~\ref{thm:Ktloc}, with no twist appearing because $\mathsf s^\qc_\core$ and $\mathsf s(u^{m-2})$ share the same quantum form (Theorem~\ref{thm:locDelta}). For part (d), we show that a quantum cluster algebra that equals its upper cluster algebra specializes, at $\omega^{1/2}=1$, to the corresponding classical cluster algebra (Lemma~\ref{lem:specialize}). By Lemma~\ref{lem:M}, the exchange matrix of $\mathsf s^\qc_\core$ coincides, up to a global sign, with that of the seed of $X(\Delta^{m-1})$ arising from the right inductive weave of the word $\Delta\cdot u^{m-2}$ (Lemma~\ref{lem:braidseed}), and the conclusion follows from~\cite{CGGLSS}.

\subsection{Organization of the paper}
In Section~\ref{sec-pre} we recall quantum cluster algebras, and we prove the facts about admissible changes of frozen coordinates that we need. In Section~\ref{sub-def-reduced-skein} we recall the reduced stated $\SL$-skein algebras and the quantum cluster structure of~\cite{CHW} for polygons, and we construct the seed $\mathsf s^\qc_{\rm skein}$ and the subalgebra $\mathcal A_{n,m}$. In Section~\ref{sec-KKOP} we recall the KKOP seeds and their monoidal categorification, and we construct the monoidal localization of the subcategory $\cC^{\mathrm s}(\mathsf b)\subseteq\cC(\mathsf b)$ at the frozen modules. In Section~\ref{sec-cat} we prove Theorems~\ref{thmA} and~\ref{thmB}. In Section~\ref{sec-braid} we prove Theorem~\ref{thmD} and Corollary~\ref{thmE}.

\subsection*{Acknowledgements}
 This work was partially supported by the National Natural Science Foundation
of China (No.12471023) (M.H.).
Part of this research was performed while the second author (Z.W.) was visiting the Institute for Pure and Applied Mathematics (IPAM), which is supported by the National Science Foundation (Grant Nos. DMS-1925919 and DMS-2422832).
Z.W. was supported by a KIAS Individual Grant (MG104701) at the Korea Institute for Advanced Study. We would like to thank Jian-Rong Li and Fan Qin for insightful and inspiring discussions.

\subsection*{AI Usage Disclosure.}
 AI tools were used to assist in the review of \cite{KKOP3} in Section~\ref{sec-KKOP}, in the development of the proof of Proposition~\ref{prop:loc}, and in proofreading and improving the clarity and readability of the manuscript.
 We also used AI to draw Figures~\ref{fig:CHWquiver} and \ref{fig:placeholder}.
 All mathematical arguments and results were verified by the authors.

\section{Preliminaries}\label{sec-pre}

In this section, we review quantum cluster algebras following Berenstein--Zelevinsky~\cite{BZ}, and quasi-homomorphisms of quantum cluster algebras following~\cite{Fraser,KWQ,CHL}.

Let $R=\Z[\omega^{\pm 1/2}]$, and set
\begin{align}\label{eq-constants}
    \xi=\omega^n,\qquad t=\omega^{n^2}=\xi^n,
\end{align}
so that $\xi^{1/2}=\omega^{n/2}$ and $t^{1/2}=\xi^{n/2}$ lie in $R$. We regard $\Z[t^{\pm1/2}]$ and $\Z[\xi^{\pm1/2}]$ as subrings of $R$. The variable $\xi$ is the one in which both the skein seed (Proposition~\ref{prop:CHWseed}) and the KKOP seed (Definition~\ref{def:KKOPseed}) are compatible, and $t$ is the parameter of the skein relations (\S\ref{sub-def-reduced-skein}) and of $(q,t)$-characters (\S\ref{subsec:qGr}).

\subsection{Quantum cluster algebras}\label{sec-mutation-quantum}

Fix a non-empty finite set $\mathcal{V}$, and a subset $\mathcal{V}_{\rm mut}$ of $\mathcal{V}$. We allow $\mathcal V_{\rm mut}=\emptyset$; then there are no mutations, and the quantum cluster algebra with invertible frozen variables is the based quantum torus of the seed (this happens for $n=2$ and $m=3$). The elements of $\mathcal{V}$, $\mathcal{V}_{\rm mut}$ and $\mathcal{V}\setminus \mathcal{V}_{\rm mut}$ are called {\bf vertices}, {\bf mutable vertices} and {\bf frozen vertices}, respectively. By an (ice) quiver $\Gamma$ we mean a directed graph whose set of vertices is $\mathcal{V}$
and equipped with weight on the edges so that the weight on each edge is $1$ unless the edge connects two frozen vertices, in which case the weight is $\frac{1}{2}$. An edge of weight $1$ will be called an {\bf arrow}, and an edge of weight $\frac{1}{2}$ a {\bf half-arrow}. We denote by $Q = (Q(u,v))_{u,v\in \mathcal{V}}$ the signed adjacency matrix of $\Gamma$, i.e.\ $Q(u,v)$ is the total weight of the edges from $u$ to $v$ minus the total weight of the edges from $v$ to $u$; it is called an {\bf exchange matrix}. (For the quivers of Section~\ref{sub-def-reduced-skein}, which are glued from triangles, half-arrows may also join two mutable vertices on a common diagonal; they cancel after summing over the two triangles, see Lemma~\ref{lem:sideedges}, and the convention above applies to the summed quiver.)
We will identify the weighted quiver with its signed adjacency matrix.

Suppose that $k\in\mathcal V_{
\rm mut}$. The {\bf mutation}
$\mu_{k}$ for the quiver $\Gamma$ at the
mutable vertex $k\in\mathcal V_{
\rm mut}$ is defined by
 the following procedures:
\begin{enumerate}
    \item for each pair of arrows $i\rightarrow k$ and $k\rightarrow j$ draw an arrow $i\rightarrow j$,

    \item reverse all the arrows incident to the vertex $k$,

    \item delete a maximal set of pairs of arrows $i\rightarrow j$ and $i\leftarrow j$ going in the opposite directions.
\end{enumerate}
We use $Q' = (Q'(u,v))_{u,v\in\mathcal V}$ to denote the signed adjacency matrix of $\Gamma'=\mu_k(\Gamma)$. Then we have
\begin{align}\label{eq-mutation-Q}
Q'(u,v) = \begin{cases}
    - Q(u,v) & k\in\{u,v\},\\
    Q(u,v) +\frac{1}{2}\Big(Q(u,k)|Q(k,v)|+
    |Q(u,k)|Q(k,v)\Big) & k\notin\{u,v\}.
\end{cases}
\end{align}

For $u\in\mathcal V$, define ${\bf e}_u\in\mathbb Z^{\mathcal V}$ by
\begin{align}\label{eq-def-vector-ei-ele}
    {\bf e}_u(v)=\delta_{u,v} \text{ for $v\in\mathcal V$}.
\end{align}
We regard elements of $\Z^{\mathcal V}$ as row vectors. For a skew-symmetric matrix $\Pi=(\Pi(u,v))_{u,v\in\mathcal V}$ and ${\bf k},{\bf t}\in\Z^{\mathcal V}$ we write
\begin{align}\label{eq-bilinear}
\Pi({\bf k},{\bf t})={\bf k}\,\Pi\,{\bf t}^T=\sum_{u,v\in\mathcal V}k_u\Pi(u,v)t_v ,
\end{align}
so that $\Pi({\bf e}_u,{\bf e}_v)=\Pi(u,v)$. For $v\in\mathcal V_{\rm mut}$, the {\bf column vector} of $Q$ at $v$ is
\begin{align}\label{eq-column-vector}
{\bf q}_v=\sum_{u\in\mathcal V}Q(u,v)\,{\bf e}_u\in\Z^{\mathcal V}.
\end{align}

Let $\mathcal F$ be a skew-field containing $R=\mathbb Z[\omega^{\pm\frac{1}{2}}]$ in its centre.
\begin{definition}[\cite{BZ}]\label{def-quantum-seed}
    A {\bf quantum seed} in $ \mathcal F$ is a triple $\mathsf s=(Q,\Pi,M)$, where
\begin{enumerate}
    \item[$(1)$] $Q=(Q(u,v))_{u,v\in\mathcal V}$ is an exchange matrix;

    \item[$(2)$] $\Pi=(\Pi(u,v))_{u,v\in\mathcal V}$ is a skew-symmetric matrix with integral entries satisfying the compatibility relation
    \begin{align}\label{eq-compatible}
    \sum_{k\in\mathcal V}Q(k,u) \Pi(k,v)=\delta_{u,v} d_u
    \end{align}
    for all $u\in\mathcal V_{\text{mut}}$
    and $v\in\mathcal V$, where $d_u$ is a positive integer for $u\in\mathcal V_{\text{mut}}$; equivalently, $\Pi({\bf q}_u,{\bf e}_v)=\delta_{u,v}d_u$;

    \item[$(3)$] $M\colon \mathbb Z^{\mathcal V}\rightarrow\mathcal F\setminus\{0\}$ is a function such that
    \begin{align}\label{eq-M-relation}
        M({\bf k}) M({\bf t}) = \xi^{\frac{1}{2} {\bf k} \Pi {\bf t}^T} M({\bf k} + {\bf t})
    \end{align}
    for row vectors ${\bf k},{\bf t}\in \mathbb Z^{\mathcal V}$, and $M(\mathbb Z^{\mathcal V})$ is a basis of a quantum torus over $R$ whose skew-field of fractions is $\mathcal F$.
\end{enumerate}
If $d_u=d$ for all $u\in\mathcal V_{\rm mut}$, we say that $(Q,\Pi)$ is {\bf compatible with scalar $d$}.
\end{definition}

By \eqref{eq-M-relation}, $M({\bf k})M({\bf t})=\xi^{\Pi({\bf k},{\bf t})}M({\bf t})M({\bf k})$. In all quantum seeds of this paper, $(Q,\Pi)$ is compatible with scalar $2n$.

The function $M$ is uniquely determined by the values $A_u:=M({\bf e}_u)$, which we call {\bf (quantum) cluster variables}.
A (quantum) cluster variable $A_u$ is called a {\bf frozen (quantum) cluster variable} if $u\in\mathcal V\setminus\mathcal V_{\text{mut}}$, and an {\bf exchangeable (quantum) cluster variable} otherwise.

We label the vertex set $\mathcal V$ as
$\{v_1,\cdots,v_r\}$.
For any ${\bf k}=(k_v)_{v\in\mathcal V}\in \mathbb Z^{\mathcal V}$, we have
\begin{align}\label{eq-A-power}
    M({\bf k}) = A^{\bf k}:= \left[\prod_{v\in \mathcal V}A_v^{k_v}\right],
\end{align}
where $\left[\prod_{v\in \mathcal V}A_v^{k_v}\right]$ is the Weyl-ordered product
\begin{align}\label{Weyl-A}
    \left[\prod_{v\in \mathcal V}A_v^{k_v}\right]= \xi^{-\frac{1}{2}\sum_{i<j} k_{v_i} k_{v_j}\Pi(v_i,v_j)} A_{v_1}^{k_{v_1}} A_{v_2}^{k_{v_2}} \cdots A_{v_r}^{k_{v_r}}.
\end{align}
Note that the Weyl-ordered product in \eqref{Weyl-A} is independent of how we label the vertex set $\mathcal V$, and that the skew-symmetric matrix $\Pi$ is uniquely determined by the function $M$ because of \eqref{eq-M-relation}. More precisely, \eqref{Weyl-A} shows the following.

\begin{lemma}\label{lem:frame-unique}
Let $M,M'\colon\Z^{\mathcal V}\to\mathcal F\setminus\{0\}$ satisfy \eqref{eq-M-relation} for the same matrix $\Pi$, and assume that $M({\bf e}_u)=M'({\bf e}_u)$ for all $u\in\mathcal V$. Then $M=M'$.
\end{lemma}

Given a quantum seed $\mathsf s=(Q,\Pi,M)$ in $\mathcal F$, the {\bf mutation} at vertex $k\in\mathcal V_{\text{mut}}$ produces a new quantum seed $(Q',\Pi',M')=\mu_{k}(\mathsf s)$ according to the following rules:
\begin{enumerate}[label={\rm (\roman*)}]
    \item $Q'=\mu_k(Q)$;

    \item the mutation on the quantum cluster variables is
    \begin{equation}\label{eq-exchange}
        A_i':=M'({\bf e}_i)=
        \begin{cases}
            M({\bf e}_i) & i\neq k,\\
            M(-{\bf e}_k + \sum_{j\in\mathcal V}[Q(j,k)]_{+}{\bf e}_j)
            + M(-{\bf e}_k + \sum_{j\in\mathcal V}[-Q(j,k)]_{+}{\bf e}_j)
            & i=k,
        \end{cases}
    \end{equation}
    where $[a]_+ := \max\{a,0\} = {\textstyle \frac{1}{2}(a+|a|)}$ for $a\in\mathbb R$;

    \item $\Pi'=(\Pi'(i,j))_{i,j\in\mathcal V}$, where, for each pair $i,j\in\mathcal V$, the entry $\Pi'(i,j)$ is the unique integer such that
    $$M'({\bf e}_i)M'({\bf e}_j)
    =\xi^{\Pi'(i,j)} M'({\bf e}_j)M'({\bf e}_i).$$
\end{enumerate}
By~\cite[\S\S3--4]{BZ}, $\mu_k(\mathsf s)$ is a quantum seed, with the same integers $d_u$, and $\mu_k$ is an involution.

\begin{remark}\label{rem:frozen-entries}
The entries $Q(u,v)$ with $u,v$ both frozen do not enter the compatibility relation~\eqref{eq-compatible}, nor the mutation~\eqref{eq-exchange} of cluster variables, nor, by~\eqref{eq-mutation-Q}, the mutation of the entries $Q(u,v)$ with $v$ mutable. Hence two quantum seeds which differ only in these entries have the same cluster variables after any sequence of mutations. We say that two exchange matrices {\bf agree} if they agree on $\mathcal V\times\mathcal V_{\rm mut}$.
\end{remark}

\begin{definition}\label{def-quantum-class}
Two quantum seeds in $\mathcal F$ are said to be
{\bf mutation-equivalent} if they are transformed to each other by a finite sequence of quantum
mutations. An equivalence class of quantum seeds is called a {\bf quantum
mutation class}.
\end{definition}

Let $\mathsf{S}$ be a quantum mutation class of quantum seeds.
For any $\mathsf s=(Q, \Pi, M)\in \mathsf{S}$, define
\begin{align}
\textbf{A}^{+}(\mathsf s)&:=\{M({\bf k})\mid {\bf k}\in\mathbb Z_{\geq 0}^{\mathcal V}\}\subset
\textbf{A}(\mathsf s):=\{M({\bf k})\mid {\bf k}\in\mathbb Z^{\mathcal V}\}\subset \mathcal F,
\label{eq-A-w}\\
\mathbb T(\mathsf s)&:= \text{span}_R(\textbf{A}(\mathsf s))\subset \mathcal F.
    \label{eq-T-w}
\end{align}
We call $\mathbb T(\mathsf s)$ the {\bf based quantum torus} associated to the quantum seed $\mathsf s$.
We use ${\rm Frac}(\mathbb T(\mathsf s))$ to denote
 the skew-field of fractions of $\mathbb T(\mathsf s)$; for the existence of ${\rm Frac}(\mathbb T(\mathsf s))$, see \cite{Cohn}.

We will define two versions of quantum (upper) cluster algebras, depending on whether the frozen variables are inverted or not.

\begin{definition}\label{def-quan-cluster-algebra}
Let $\mathsf{S}$ be a quantum mutation class of quantum seeds.

\begin{enumerate}[label={\rm (\alph*)}]\itemsep0.3em
\item
The associated {\bf quantum cluster algebra} is the $R$-subalgebra
$\mathscr A(\mathsf S)\subset \mathcal F$
generated by the set $\bigcup_{\mathsf s\in\Exc} {\bf A}^{+}(\mathsf s)$ (see~\eqref{eq-A-w}) together with the inverses of all frozen quantum cluster variables.
We also denote by
$\mathscr A^+(\mathsf S)\subset \mathcal F$
the $R$-subalgebra generated by $\bigcup_{\mathsf s\in\Exc} {\bf A}^{+}(\mathsf s)$.

\item Each exchangeable quantum cluster variable of any seed $\mathsf s\in\mathsf S$ is called an {\bf exchangeable quantum cluster variable} of $\mathscr A(\mathsf S)$.

\item The associated {\bf quantum upper cluster algebra} is
\[
\mathscr U(\mathsf S)
:= \bigcap_{\mathsf s\in\Exc} \mathbb T(\mathsf s)
\;\subset\; \mathcal F,
\]
where $\mathbb T(\mathsf s)$ is as in~\eqref{eq-T-w}.

\item Let $\mathsf s=(Q, \Pi, M)\in \mathsf{S}$ be a quantum seed. Let ${\bf k}\in\mathbb Z^{\mathcal V}$ such that
${\bf k}|_{\mathcal V_{\rm mut}}\in\mathbb Z_{\geq 0}^{\mathcal V_{\rm mut}}$. Then $M({\bf k})$ is called a \textbf{cluster monomial} in $\mathscr A(\mathsf S)$. If, in addition, ${\bf k}\in\mathbb Z_{\geq 0}^{\mathcal V}$, then $M({\bf k})$ is called a \textbf{cluster monomial} in $\mathscr A^+(\mathsf S)$.
\end{enumerate}
\end{definition}

For a quantum seed $\mathsf s\in \mathsf S$, we write $\mathscr A(\mathsf s)$, $\mathscr A^+(\mathsf s)$ and $\mathscr U(\mathsf s)$ for $\mathscr A(\mathsf S)$, $\mathscr A^+(\mathsf S)$ and $\mathscr U(\mathsf S)$, respectively, and call them the quantum cluster algebra and the quantum upper cluster algebra associated with $\mathsf s$.

The following theorem is the so-called Laurent phenomenon.

\begin{theorem}[\cite{BZ,goodearl2017quantum}]\label{thm-inclusion-quantum}
    We have $\mathscr A^+(\mathsf S)\subset\mathscr A(\mathsf S)\subset \mathscr U(\mathsf S)$.
\end{theorem}

For $\mathsf s=(Q,\Pi,M)\in\mathsf S$ and $g\in\mathcal V\setminus\mathcal V_{\rm mut}$, the frozen variable $A_g=M({\bf e}_g)$ is a cluster variable of every seed in $\mathsf S$, so it quasi-commutes with every cluster variable of every seed. Hence $A_g\mathscr A^+(\mathsf s)=\mathscr A^+(\mathsf s)A_g$, and we get the following.

\begin{lemma}\label{lem:normal}
The frozen variables are normal elements in $\mathscr A^+(\mathsf s)$. The multiplicative set generated by them forms an Ore set, and $\mathscr A(\mathsf s)$ is the corresponding Ore localization of $\mathscr A^+(\mathsf s)$.
\end{lemma}

We next compare the quantum cluster algebra with its classical limit. Let $h=\omega^{1/2}-1\in R$, so that $R/hR \cong \Z$. For a quantum seed $\mathsf s=(Q,\Pi,M)$, the quotient $\bT(\mathsf s)/h\bT(\mathsf s)$ is the Laurent polynomial ring $\Z[A_v^{\pm1}:v\in\mathcal V]$. Then the \emph{classical limit} of $\mathsf s$, denoted by $\mathsf s\mid_{\omega^{1/2}=1}$, is the commutative seed with exchange matrix $Q$ and initial cluster $(A_v)_{v\in\mathcal V}$. Let $\mathscr A_1(\mathsf s)\subseteq\Z[A_v^{\pm1}]$ denote the commutative cluster algebra associated with $\mathsf s\mid_{\omega^{1/2}=1}$ over $\Z$. Since the exchange relations~\eqref{eq-exchange} reduce to the commutative exchange relations modulo $h$, the image of every quantum cluster variable under the canonical projection coincides with the corresponding classical cluster variable. The following lemma can be viewed as a variant of the specialization results in~\cite{GLS-spec}.


\begin{lemma}\label{lem:specialize}
Let $\mathsf s$ be a quantum seed satisfying $\mathscr A(\mathsf s)=\mathscr U(\mathsf s)$. The natural map
\[\iota:\mathscr A(\mathsf s)\otimes_R\Z\to\bT(\mathsf s)\otimes_R\Z=\Z[A_v^{\pm1}:v\in\mathcal V]\] is injective with image $\mathscr A_1(\mathsf s)$, where $R$ acts on $\Z$ via $\omega^{1/2}\mapsto1$.

Consequently, for any commutative ring $\mathbf k'$ that is flat over $\Z$ (for instance, $\mathbf k'=\C$), regarded as an $R$-algebra via $\omega^{1/2}\mapsto1$, we have an isomorphism $\mathscr A(\mathsf s)\otimes_R\mathbf k'\cong\mathscr A_1(\mathsf s)\otimes_\Z\mathbf k'$. This identifies the former with the $\mathbf k'$-subalgebra of $\mathbf k'[A_v^{\pm1}:v\in\mathcal V]$ generated by the commutative cluster variables and the inverses of the frozen variables.
\end{lemma}

\begin{proof}
Let $\mathsf s_\tau$ be a seed mutation equivalent to $\mathsf s$. Let $\mathsf s_{\tau'}=\mu_k(\mathsf s_\tau)$ for some $k\in \mathcal V_{\rm mut}$. Denote the toric frames by $M_\tau$ and $M_{\tau'}$, and the cluster variables by $A^{(\tau)}_v$ and $A^{(\tau')}_v$, respectively. Let ${\bf q}_k$ be the $k$-th column of the exchange matrix of $\mathsf s_\tau$, and set ${\bf b}^\pm=[\pm{\bf q}_k]_+$. These vectors have vanishing $k$-coordinate, so $M_\tau({\bf b}^\pm)=M_{\tau'}({\bf b}^\pm)$. Since ${\bf b}^+-{\bf b}^-={\bf q}_k$, the compatibility relation~\eqref{eq-compatible} implies $\Pi({\bf b}^+,{\bf b}^-)=\Pi({\bf q}_k,{\bf b}^-)=0$. Consequently, $M_{\tau'}({\bf b}^+)$ and $M_{\tau'}({\bf b}^-)$ commute.

For $a\in\frac{1}{2}\mathbb{Z}$, let
\[
M^a=M_{\tau'}({\bf b}^+)+\xi^{a}M_{\tau'}({\bf b}^-)\in\bT(\mathsf s_{\tau'}),
\]
and let $S\subseteq\bT(\mathsf s_{\tau'})$ be the multiplicative set generated by the elements $M^a$ with $a\in\frac{1}{2}\mathbb{Z}$. Clearly, $S$ is commutative. By~\eqref{eq-M-relation}, every monomial $M_{\tau'}({\bf c})$ satisfies $M_{\tau'}({\bf c})\,M^a=\xi^{e}\,M^{a'}\,M_{\tau'}({\bf c})$ for some $e,a'\in\frac{1}{2}\mathbb{Z}$. Furthermore, by~\eqref{eq-exchange}, we have $A^{(\tau)}_k=\xi^{e_0}\,M^{a_0}(A^{(\tau')}_k)^{-1}$ for some $e_0,a_0\in\frac{1}{2}\mathbb{Z}$. Thus $A^{(\tau)}_k\in\bT(\mathsf s_{\tau'})$, and in the fraction field $\mathcal F$, $(A^{(\tau)}_k)^{-1}=\xi^{-e_0}A^{(\tau')}_k(M^{a_0})^{-1}$. By moving the factors $(M^{a_0})^{-1}$ to the right, each power $(A^{(\tau)}_k)^{-N}$ ($N\ge1$) can be written in the form $y\,\sigma^{-1}$ with $y\in\bT(\mathsf s_{\tau'})$ and $\sigma\in S$. Since the other cluster variables of $\mathsf s_\tau$ and $\mathsf s_{\tau'}$ coincide, every $z\in\bT(\mathsf s_\tau)$ admits a representation $z=\sum_j y_j\sigma_j^{-1}$ with $y_j\in\bT(\mathsf s_{\tau'})$ and $\sigma_j\in S$. As $S$ is commutative, taking $\sigma=\prod_j\sigma_j\in S$ ensures that $z\sigma\in\bT(\mathsf s_{\tau'})$.

We now use these facts to prove that $\mathscr U(\mathsf s)\cap h\bT(\mathsf s)=h\,\mathscr U(\mathsf s)$.

Let $x\in\mathscr U(\mathsf s)\cap h\bT(\mathsf s_\tau)$. Write $x=hz$ with $z\in\bT(\mathsf s_\tau)$, and choose $\sigma$ as constructed above. Then $x\sigma=hz\sigma\in h\bT(\mathsf s_{\tau'})$. Let $\iota: \bT(\mathsf s_{\tau'}) \to \bT(\mathsf s_{\tau'})/h\bT(\mathsf s_{\tau'})$ be the canonical projection. Then $\iota(x\sigma)=0$. Note that $\iota(\sigma)$ is a power of $\iota(M_{\tau'}({\bf b}^+)+M_{\tau'}({\bf b}^-))$, which is a sum of two monomials with coefficient $1$ and hence non-zero. Since $\bT(\mathsf s_{\tau'})/h\bT(\mathsf s_{\tau'})$ is a commutative Laurent polynomial ring over $\mathbb{Z}$ (and thus an integral domain), it follows that $\iota(x)=0$, i.e., $x\in h\bT(\mathsf s_{\tau'})$. By induction along mutation sequences, we deduce that $x\in h\bT(\mathsf s_\tau)$ for every seed in the mutation class, which implies $x\in h\mathscr U(\mathsf s)$. This proves the inclusion $\mathscr U(\mathsf s)\cap h\bT(\mathsf s)\subseteq h\,\mathscr U(\mathsf s)$; the reverse inclusion is immediate. Thus, $\mathscr U(\mathsf s)\cap h\bT(\mathsf s)=h\,\mathscr U(\mathsf s)$.

Consequently, the map $\mathscr U(\mathsf s)/h\mathscr U(\mathsf s)\to\bT(\mathsf s)/h\bT(\mathsf s)=\mathbb{Z}[A_v^{\pm1}]$ is injective. Since $\mathscr A(\mathsf s)=\mathscr U(\mathsf s)$ and $\mathscr A(\mathsf s)\otimes_R\mathbb{Z} \cong \mathscr A(\mathsf s)/h\mathscr A(\mathsf s)$, the injectivity of $\iota$ is established. The image is generated by the images of the generators of $\mathscr A(\mathsf s)$, namely the quantum cluster variables and the inverses of the frozen variables, which correspond to the commutative cluster variables and the inverses of the commutative frozen variables, respectively. Therefore, the image is $\mathscr A_1(\mathsf s)$. The final statement follows by tensoring the inclusion $\mathscr A_1(\mathsf s)\subseteq\mathbb{Z}[A_v^{\pm1}]$ with the flat $\mathbb{Z}$-module $\mathbf k'$.
\end{proof}

The above lemma says that, when a quantum cluster algebra equals its upper cluster algebra, setting $\omega^{1/2}=1$ loses no information and gives exactly the commutative cluster algebra of the classical limit. We apply it to the seeds $\mathsf s_{\rm skein}$ and $\mathsf s^\qc_{\rm skein}$, for which $\mathscr A=\mathscr U$ by Theorem~\ref{thm:CHW} and Proposition~\ref{prop:sqc}(b), and to the seed $\mathsf s^\qc_\core$, for which $\mathscr A=\mathscr U$ by Theorem~\ref{thm:Aquant}(b).

The following is clear.

\begin{lemma}\label{lem:seediso}
Let $\mathsf s=(Q,\Pi,M)$ and $\mathsf s'=(Q',\Pi,M')$ be quantum seeds in skew-fields $\mathcal F$ and $\mathcal F'$, on the same vertex set, with the same mutable vertices and the same skew-symmetric matrix $\Pi$, and such that $Q$ and $Q'$ agree in the sense of Remark~\ref{rem:frozen-entries}. Then there is a unique isomorphism of skew-fields $\mathcal F\to\mathcal F'$ with $M({\bf k})\mapsto M'({\bf k})$ for ${\bf k}\in\Z^{\mathcal V}$. It maps the cluster variables of $\mu_{k_l}\cdots\mu_{k_1}(\mathsf s)$ to those of $\mu_{k_l}\cdots\mu_{k_1}(\mathsf s')$ for every sequence $k_1,\dots,k_l$ of mutable vertices, and it restricts to isomorphisms $\mathscr A^+(\mathsf s)\cong\mathscr A^+(\mathsf s')$, $\mathscr A(\mathsf s)\cong\mathscr A(\mathsf s')$, $\mathscr U(\mathsf s)\cong\mathscr U(\mathsf s')$, which is compatible with mutation.
\end{lemma}


Let $\mathsf s=(Q,\Pi,M)$ be a quantum seed. Suppose $\mathcal V'=\mathcal V'_1\sqcup \mathcal V'_2$ is a subset of $\mathcal V$ with $\mathcal V_1'\subseteq\mathcal V_{\rm mut}$ and $\mathcal V'_2\subseteq\mathcal V\setminus\mathcal V_{\rm mut}$, such that there are no arrows between vertices in $\mathcal V'_1$ and those in $\mathcal V\setminus \mathcal V'$. Let $Q'$ be the subquiver of $Q$ induced by the vertices in $\mathcal V'$, with the vertices in $\mathcal V'_1$ mutable and those in $\mathcal V'_2$ frozen. Let $\Pi'$ and $M'$ be the restrictions of $\Pi$ and $M$ to $\mathcal V'$ and to $\Z^{\mathcal V'}\subseteq\Z^{\mathcal V}$, respectively. The following is immediate from \eqref{eq-compatible} and \eqref{eq-exchange}.

\begin{lemma}\label{lem:subseed}
   Using the above notation, we have:
   \begin{enumerate}[label={\rm (\alph*)}]\itemsep0.3em
       \item The triple $\mathsf s'=(Q',\Pi',M')$ constitutes a quantum seed in the skew-field of fractions of $\bT(\mathsf s')\subseteq\mathcal F$, with the same integers $d_u$, $u\in\mathcal V'_1$.
       \item The quantum cluster algebra $\mathscr A(\mathsf s')$ is a subalgebra of $\mathscr A(\mathsf s)$. Specifically, for any sequence of mutable vertices $v_1,\dots,v_l$ in $Q'$, it holds that
       $\widetilde A'_{u} = \widetilde A_{u}$ for all $u\in\mathcal V'$.
       Here, $\widetilde A'_{u}$ and $\widetilde A_{u}$ are the cluster variables associated with the vertex $u$ of $\mu_{v_l}\cdots\mu_{v_1}(\mathsf s')$ and $\mu_{v_l}\cdots\mu_{v_1}(\mathsf s)$, respectively.
    \end{enumerate}
\end{lemma}

\begin{definition}\label{def:subseed}
With the above notation, we call the seed $\mathsf s'=(Q',\Pi',M')$ the \textbf{sub-seed} of $\mathsf s$ supported on the vertex set $\mathcal V'$, with the vertices in $\mathcal V'_2$ frozen.
\end{definition}

\subsection{Quasi-isomorphism of quantum cluster algebras}\label{sec-quasi}
Let $\mathsf s=(Q,\Pi,M)$ and  $\mathsf s'=(Q',\Pi',M')$ be two quantum seeds. Assume that $Q$ and $Q'$ have the same mutable vertices $\mathcal V_{\rm mut}$. Denote the vertices of $Q$ and $Q'$ by $\mathcal V$ and $\mathcal V'$, respectively.

\begin{definition}[\cite{KWQ,CHL}]\label{def:QQH}
We call an $R$-algebra homomorphism $\varphi: \mathscr A(\mathsf s) \to \mathscr A(\mathsf s')$ a {\bf quantum quasi-homomorphism} (with respect to $\mathsf s$ and $\mathsf s'$) if the following conditions are satisfied:
\begin{enumerate}

\item[$(1)$] The mutable subquivers of $Q$ and $Q'$ are the same, i.e, $Q\mid_{\mathcal V_{\rm mut}}=Q'\mid_{\mathcal V_{\rm mut}}$;

\item[$(2)$] For $v\in \mathcal V_{\rm mut}$, $\varphi(M(\mathbf e_v))=M'({\bf f}_v)$ for some ${\bf f}_v=(f_{v,w})_{w\in \mathcal V'}$ with $f_{v,v}=1$ and $f_{v,w}=0$ for all $w(\neq v)\in \mathcal V_{\rm mut}$.

\item[$(3)$] For $v\in \mathcal V\setminus \mathcal V_{\rm mut}$, $\varphi(M(\mathbf e_v))=M'({\bf f}_v)$ for some ${\bf f}_v=(f_{v,w})_{w\in \mathcal V'}$ with $f_{v,w}=0$ for all $w\in \mathcal V_{\rm mut}$.

\item[$(4)$] For all $v\in \mathcal V_{\rm mut}$, $\varphi(M({\bf q}_v))=M'({\bf q}'_v)$, where ${\bf q}_v$ and ${\bf q}'_v$ are the column vectors~\eqref{eq-column-vector} of $Q$ and $Q'$ at $v$.
\end{enumerate}
A quantum quasi-homomorphism is a {\bf quantum quasi-isomorphism} if it is an isomorphism of algebras.
\end{definition}

\begin{proposition}[{\cite[Prop.~2.7]{CHL}}]\label{prop:quasi-iso}
Assume that $Q|_{\mathcal{V}_{\mathrm{mut}}} = Q'|_{\mathcal{V}_{\mathrm{mut}}}$.
A map
$$\varphi: \{M(\mathbf{e}_v) \mid v \in \mathcal{V}\} \to \mathbb{T}(\mathsf{s}')$$
induces a quasi-homomorphism $\varphi: \mathscr A(\mathsf{s}) \to \mathscr A(\mathsf{s}')$ if and only if the following conditions hold:
\begin{enumerate}
    \item[$(1)$] For each $v \in \mathcal{V}_{\mathrm{mut}}$, we have $\varphi(M(\mathbf{e}_v)) = M'(\mathbf{f}_v)$ for some vector $\mathbf{f}_v = (f_{v,w})_{w \in \mathcal{V}'}$ satisfying $f_{v,v} = 1$ and $f_{v,w} = 0$ for all $w \in \mathcal{V}_{\mathrm{mut}} \setminus \{v\}$.

    \item[$(2)$] For each $v \in \mathcal{V} \setminus \mathcal{V}_{\mathrm{mut}}$, we have $\varphi(M(\mathbf{e}_v)) = M'(\mathbf{f}_v)$ for some vector $\mathbf{f}_v = (f_{v,w})_{w \in \mathcal{V}'}$ such that $f_{v,w} = 0$ for all $w \in \mathcal{V}_{\mathrm{mut}}$.

    \item[$(3)$] For each $v \in \mathcal{V}_{\mathrm{mut}}$, the following two equivalent conditions hold:
    \[
        \sum_{w \in \mathcal{V}} Q(w,v)\mathbf{f}_w = \sum_{w' \in \mathcal{V}'} Q'(w',v)\mathbf{e}_{w'},
    \qquad
        \varphi\Bigl(\Bigl[\prod_{w \in \mathcal{V}}A_w^{Q(w,v)}\Bigr]\Bigr) = \Bigl[\prod_{w' \in \mathcal{V'}}A_{w'}^{Q'(w',v)}\Bigr].
    \]

    \item[$(4)$] For all $v_1, v_2 \in \mathcal{V}$, $\Pi(\mathbf{e}_{v_1}, \mathbf{e}_{v_2}) = \Pi'(\mathbf{f}_{v_1}, \mathbf{f}_{v_2})$.
\end{enumerate}
\end{proposition}

We need the following special case, in which $\mathcal V=\mathcal V'$, the map $\varphi$ is the identity of a subalgebra of $\mathcal F$, and the vectors ${\bf f}_v$ come from a linear automorphism of $\Z^{\mathcal V}$.

\begin{definition}\label{def:admissible}
Let $\mathsf s=(Q,\Pi,M)$ be a quantum seed. A $\Z$-linear automorphism $\mathsf G$ of $\Z^{\mathcal V}$ is called {\bf admissible} if
\begin{align}\label{eq-admissible}
\mathsf G({\bf e}_v)-{\bf e}_v\in\Z^{\mathcal V\setminus\mathcal V_{\rm mut}}\ \ (v\in\mathcal V_{\rm mut}),\qquad \mathsf G\bigl(\Z^{\mathcal V\setminus\mathcal V_{\rm mut}}\bigr)=\Z^{\mathcal V\setminus\mathcal V_{\rm mut}} .
\end{align}
For an admissible $\mathsf G$ we define $\mathsf s^{\mathsf G}=(Q^{\mathsf G},\Pi^{\mathsf G},M^{\mathsf G})$ by
\begin{align}\label{eq-sP}
M^{\mathsf G}({\bf k})=M(\mathsf G({\bf k})),\qquad \Pi^{\mathsf G}({\bf k},{\bf t})=\Pi(\mathsf G({\bf k}),\mathsf G({\bf t})),\qquad {\bf q}^{\mathsf G}_v=\mathsf G^{-1}({\bf q}_v)\ \ (v\in\mathcal V_{\rm mut}),
\end{align}
where ${\bf q}^{\mathsf G}_v$ is the column vector of $Q^{\mathsf G}$ at $v$. This determines $Q^{\mathsf G}(u,v)$ for $v\in\mathcal V_{\rm mut}$; for $u\in\mathcal V_{\rm mut}$ and $v$ frozen we put $Q^{\mathsf G}(u,v)=-Q^{\mathsf G}(v,u)$, and the entries of $Q^{\mathsf G}$ between frozen vertices are irrelevant (Remark~\ref{rem:frozen-entries}), so we set them to be $0$.
\end{definition}

Admissible automorphisms form a group, and $\mathsf G^{-1}({\bf e}_v)-{\bf e}_v\in\Z^{\mathcal V\setminus\mathcal V_{\rm mut}}$ for $v\in\mathcal V_{\rm mut}$. In particular, ${\bf q}^{\mathsf G}_v$ and ${\bf q}_v$ have the same $\mathcal V_{\rm mut}$-coordinates, i.e.\ $Q^{\mathsf G}$ and $Q$ have the same mutable part.

The next lemma says that an admissible change of coordinates transforms a quantum seed into a quantum seed, and commutes with mutation up to another admissible change.

\begin{lemma}[\cite{CHW}]\label{lem:Pchange}
Let $\mathsf s=(Q,\Pi,M)$ be a quantum seed with integers $d_u$, and let $\mathsf G$ be admissible. Then $\mathsf s^{\mathsf G}$ is a quantum seed with the same integers $d_u$, in the same skew-field. Moreover, for every $k\in\mathcal V_{\rm mut}$ there is an admissible $\mathsf G'$, with $\mathsf G'|_{\Z^{\mathcal V\setminus\mathcal V_{\rm mut}}}=\mathsf G|_{\Z^{\mathcal V\setminus\mathcal V_{\rm mut}}}$ and $\mathsf G'({\bf e}_v)=\mathsf G({\bf e}_v)$ for $v\in\mathcal V_{\rm mut}\setminus\{k\}$, such that
\begin{align}\label{eq-Pchange-mutation}
\mu_k(\mathsf s^{\mathsf G})=\mu_k(\mathsf s)^{\mathsf G'} ,
\end{align}
where the two exchange matrices agree in the sense of Remark~\ref{rem:frozen-entries}. If moreover $F\subseteq\mathcal V\setminus\mathcal V_{\rm mut}$ is a set of frozen vertices such that $\mathsf G({\bf e}_v)-{\bf e}_v\in\Z^F$ for all $v\in\mathcal V$ and $\mathsf G({\bf e}_g)={\bf e}_g$ for all $g\in F$, then $\mathsf G'$ has the same two properties.
\end{lemma}

Applying Lemma~\ref{lem:Pchange} along mutation sequences, we get the following.


\begin{corollary}\label{cor:Pchange}
Let $\mathsf s=(Q,\Pi,M)$ be a quantum seed and let $\mathsf G$ be admissible. Then
\begin{align}\label{eq-Pchange-equal}
\mathscr A(\mathsf s^{\mathsf G})=\mathscr A(\mathsf s)\qquad\text{and}\qquad\mathscr U(\mathsf s^{\mathsf G})=\mathscr U(\mathsf s)
\end{align}
as subalgebras of $\mathcal F$, and the identity map $\mathscr A(\mathsf s^{\mathsf G})\to\mathscr A(\mathsf s)$ is a quantum quasi-isomorphism (Definition~\ref{def:QQH}), with ${\bf f}_v=\mathsf G({\bf e}_v)$. Moreover, up to a power of $\xi^{1/2}$, each cluster variable of a mutated seed of $\mathsf s^{\mathsf G}$ is the corresponding cluster variable of the same mutated seed of $\mathsf s$ times a Laurent monomial in the frozen variables, which only involves the $A_g$, $g\in F$, if $\mathsf G$ satisfies the hypothesis of the last statement of Lemma~\ref{lem:Pchange} for $F$.
\end{corollary}

\begin{proof}
By Lemma~\ref{lem:Pchange} and induction on the length of a mutation sequence $k_1,\dots,k_l$, there is an admissible $\mathsf G_l$, which agrees with $\mathsf G$ on $\Z^{\mathcal V\setminus\mathcal V_{\rm mut}}$ and satisfies the hypothesis on $F$ if $\mathsf G$ does, such that $\mu_{k_l}\cdots\mu_{k_1}(\mathsf s^{\mathsf G})=\bigl(\mu_{k_l}\cdots\mu_{k_1}(\mathsf s)\bigr)^{\mathsf G_l}$ in the sense of Remark~\ref{rem:frozen-entries}; let $M_l$ be the toric frame of $\mu_{k_l}\cdots\mu_{k_1}(\mathsf s)$. For $v$ mutable, the cluster variable $M_l(\mathsf G_l({\bf e}_v))=M_l({\bf e}_v+\boldsymbol\pi)$, with $\boldsymbol\pi$ frozen, equals a power of $\xi^{1/2}$ times $M_l(\boldsymbol\pi)M_l({\bf e}_v)$. The frozen variables of $\mathsf s^{\mathsf G}$ are the $M(\mathsf G({\bf e}_g))$, which generate the same group of frozen Laurent monomials as the $A_g$ up to powers of $\xi^{1/2}$, since $\mathsf G$ restricts to an automorphism of $\Z^{\mathcal V\setminus\mathcal V_{\rm mut}}$. So both algebras in~\eqref{eq-Pchange-equal} are generated by the same elements up to units of the frozen Laurent ring. For the upper algebras, the based quantum tori $\bT$ of corresponding seeds are equal, because $\mathsf G_l$ is invertible. The conditions of Proposition~\ref{prop:quasi-iso} hold with ${\bf f}_v=\mathsf G({\bf e}_v)$: (1) and (2) are~\eqref{eq-admissible}, (3) is ${\bf q}_v=\mathsf G({\bf q}^{\mathsf G}_v)$, and (4) is the definition of $\Pi^{\mathsf G}$.
\end{proof}

\section{(Reduced) stated $\SL$-skein algebras and quantum cluster algebra structure}\label{sub-def-reduced-skein}

\def\Si{\fS}

Stated $\SL$-skein algebras were introduced by L\^e and Sikora~\cite{LS}, extending the stated $\mathrm{SL}_2$-skein algebras of L\^e~\cite{Le18} and Costantino--L\^e~\cite{CL22}, and their reduced versions were introduced by L\^e and Yu~\cite{LY}. In this section, we review their definitions in \S\ref{subsec-def-skein}, and the quantum cluster algebra structure on the reduced stated $\SL$-skein algebra $\overline{\cS}(\bP_m)$ of a polygon established in~\cite{CHW} in \S\ref{subsec-CHW}. In \S\ref{subsec-qc} we use the edge weights of stated $n$-webs to construct a quantum seed $\mathsf s^\qc_{\rm skein}$ which is quasi-isomorphic to the seed of~\cite{CHW} and in which $(m-1)(n-1)$ frozen vertices are isolated, and we define the subalgebra $\mathcal A_{n,m}$. In Section~\ref{sec-cat} the seed $\mathsf s^\qc_{\rm skein}$ is compared with the KKOP seed of Section~\ref{sec-KKOP}, and in Section~\ref{sec-braid} the subalgebra $\mathcal A_{n,m}$ is identified with a quantization of a braid variety.

\subsection{(Reduced) stated $\SL$-skein algebras}\label{subsec-def-skein}

A \emph{punctured bordered surface} (or pb surface) $\fS$ is obtained from a compact oriented surface by removing finitely many points, called \emph{punctures}, such that every boundary component is an open interval. For $m\ge1$, the polygon $\bP_m$ is the closed disk with $m$ punctures removed from its boundary. We label the punctures of $\bP_m$ by $1,2,\dots,m$ in clockwise order~\cite[\S5]{CHW}.

Consider the 3-manifold $\Si \times (-1,1)$, the thickened surface. For a point $(x,t) \in \Si \times (-1,1)$, the value $t$ is called the {\bf height} of this point.
\begin{definition}[\cite{LS}]\label{def-n-web}
    An {\bf $n$-web} $\alpha$ in $\Si\times(-1,1)$ is a disjoint union of oriented closed curves and a directed finite graph properly embedded into $\Si\times(-1,1)$, satisfying the following requirements:
\begin{enumerate}
    \item[$(1)$] $\alpha$ only contains $1$-valent or $n$-valent vertices. Each $n$-valent vertex is a source or a  sink. The set of $1$-valent vertices is denoted as $\partial \alpha$, which are called \textbf{endpoints} of $\alpha$. For any boundary component $c$ of $\Si$, we require that the points of  $\partial\alpha\cap (c\times(-1,1))$ have mutually distinct heights.
    \item[$(2)$] Every edge of the graph is an embedded oriented  closed interval  in $\Si\times(-1,1)$.
    \item[$(3)$] $\alpha$ is equipped with a transversal \textbf{framing}.
    \item[$(4)$] The set of half-edges at each $n$-valent vertex is equipped with a  cyclic order.
    \item[$(5)$] $\partial \alpha$ is contained in $\partial\Si\times (-1,1)$ and the framing at these endpoints is given by the positive direction of $(-1,1)$.
\end{enumerate}
We will consider $n$-webs up to (ambient) \textbf{isotopy} which are continuous deformations of $n$-webs in their class.
The empty $n$-web, denoted by $\emptyset$, is also considered as an $n$-web, with the convention that $\emptyset$ is only isotopic to itself.

A {\bf state} for $\alpha$ is a map $s\colon\partial\alpha\rightarrow \{1,2,\cdots,n\}$. A {\bf stated $n$-web} in $\Si\times(-1,1)$ is an $n$-web equipped with a state.
\end{definition}

Webs of this kind go back to Kuperberg~\cite{Kup96} for $n=3$, and to Sikora~\cite{Sik05} and Cautis--Kamnitzer--Morrison~\cite{CKM14} for general $n$.

By regarding $\fS$ as $\fS\times\{0\}$, there is a projection $\text{pr}\colon\fS\times(-1,1)\rightarrow \fS.$
We say the (stated) $n$-web $\alpha$ is in {\bf vertical position} if
\begin{enumerate}
    \item the framing at everywhere is given by the positive direction of $(-1,1)$,
    \item $\alpha$ is in general position with respect to the projection  $\text{pr}\colon \Si\times(-1,1)\rightarrow \Si\times\{0\}$,
    \item at every $n$-valent vertex, the cyclic order of half-edges as the image of $\text{pr}$ is given by the positive orientation of $\Si$ (drawn counter-clockwise in pictures).
\end{enumerate}
For every (stated) $n$-web $\alpha$, we can isotope $\alpha$ to be in vertical position. For each boundary component $c$ of $\fS$, the heights of $\partial\alpha\cap (c\times(-1,1))$ determine a linear order on  $c\cap \text{pr}(\alpha)$.
Then a {\bf (stated) $n$-web diagram} of $\alpha$ is $\text{pr}(\alpha)$ equipped with the usual over/underpassing information at each double point (called a crossing) and a linear order on $c\cap \text{pr}(\alpha)$ for each boundary component $c$ of $\fS$.

In this paper, we always assume that the orientations of depicted
surfaces point to the readers, i.e., the surfaces are equipped with the orientations in counterclockwise.

Let $S_n$ denote the permutation group on the set $\{1,2,\cdots,n\}$.
For an integer $i\in\{1,2,\cdots,n\}$, we use $\bar{i}$ to denote $n+1-i$. Define the following constants:
\begin{align}\label{eq-skein-constants}
\mathbbm{c}_i=(-t)^{n-i}\,t^{\frac{n-1}{2n}}\quad(1\le i\le n),\qquad
\mathbbm{t}=(-1)^{n-1}\,t^{\frac{n^2-1}{n}},\qquad
\mathbbm{a}=t^{\frac{n+1-2n^2}{4}},\qquad [n]=\frac{t^n-t^{-n}}{t-t^{-1}} .
\end{align}
Since $t=\omega^{n^2}$, all of them, as well as $t^{\pm\frac1n}=\omega^{\pm n}$, lie in $R$.

The \textbf{stated $\SL$-skein algebra} $\cS(\fS)$ of $\fS$ is
the quotient module of the $R$-module freely generated by the set
 of all isotopy classes of stated
$n$-webs in $\fS\times (-1,1)$ subject to  relations \eqref{w.cross}-\eqref{wzh.eight}.

\begin{equation}\label{w.cross}
t^{\frac{1}{n}}
\raisebox{-.20in}{
\begin{tikzpicture}
\filldraw[draw=white,fill=gray!20] (-0,-0.2) rectangle (1, 1.2);
\draw [line width =1pt,decoration={markings, mark=at position 0.5 with {\arrow{>}}},postaction={decorate}](0.6,0.6)--(1,1);
\draw [line width =1pt,decoration={markings, mark=at position 0.5 with {\arrow{>}}},postaction={decorate}](0.6,0.4)--(1,0);
\draw[line width =1pt] (0,0)--(0.4,0.4);
\draw[line width =1pt] (0,1)--(0.4,0.6);
\draw[line width =1pt] (0.4,0.6)--(0.6,0.4);
\end{tikzpicture}
}
- t^{-\frac{1}{n}}
\raisebox{-.20in}{
\begin{tikzpicture}
\filldraw[draw=white,fill=gray!20] (-0,-0.2) rectangle (1, 1.2);
\draw [line width =1pt,decoration={markings, mark=at position 0.5 with {\arrow{>}}},postaction={decorate}](0.6,0.6)--(1,1);
\draw [line width =1pt,decoration={markings, mark=at position 0.5 with {\arrow{>}}},postaction={decorate}](0.6,0.4)--(1,0);
\draw[line width =1pt] (0,0)--(0.4,0.4);
\draw[line width =1pt] (0,1)--(0.4,0.6);
\draw[line width =1pt] (0.6,0.6)--(0.4,0.4);
\end{tikzpicture}
}
= (t-t^{-1})
\raisebox{-.20in}{
\begin{tikzpicture}
\filldraw[draw=white,fill=gray!20] (-0,-0.2) rectangle (1, 1.2);
\draw [line width =1pt,decoration={markings, mark=at position 0.5 with {\arrow{>}}},postaction={decorate}](0,0.8)--(1,0.8);
\draw [line width =1pt,decoration={markings, mark=at position 0.5 with {\arrow{>}}},postaction={decorate}](0,0.2)--(1,0.2);
\end{tikzpicture}
},
\end{equation}
\begin{equation}\label{w.twist}
\raisebox{-.15in}{
\begin{tikzpicture}
\filldraw[draw=white,fill=gray!20] (-1,-0.35) rectangle (0.6, 0.65);
\draw [line width =1pt,decoration={markings, mark=at position 0.5 with {\arrow{>}}},postaction={decorate}](-1,0)--(-0.25,0);
\draw [color = black, line width =1pt](0,0)--(0.6,0);
\draw [color = black, line width =1pt] (0.166 ,0.08) arc (-37:270:0.2);
\end{tikzpicture}}
= \mathbbm{t}
\raisebox{-.15in}{
\begin{tikzpicture}
\filldraw[draw=white,fill=gray!20] (-1,-0.5) rectangle (0.6, 0.5);
\draw [line width =1pt,decoration={markings, mark=at position 0.5 with {\arrow{>}}},postaction={decorate}](-1,0)--(-0.25,0);
\draw [color = black, line width =1pt](-0.25,0)--(0.6,0);
\end{tikzpicture}},
\end{equation}
\begin{equation}\label{w.unknot}
\raisebox{-.20in}{
\begin{tikzpicture}
\filldraw[draw=white,fill=gray!20] (0,0) rectangle (1,1);
\draw [line width =1pt,decoration={markings, mark=at position 0.5 with {\arrow{>}}},postaction={decorate}](0.45,0.8)--(0.55,0.8);
\draw[line width =1pt] (0.5 ,0.5) circle (0.3);
\end{tikzpicture}}
= (-1)^{n-1} [n]\
\raisebox{-.20in}{
\begin{tikzpicture}
\filldraw[draw=white,fill=gray!20] (0,0) rectangle (1,1);
\end{tikzpicture}},
\end{equation}
\begin{equation}\label{wzh.four}
\raisebox{-.30in}{
\begin{tikzpicture}
\filldraw[draw=white,fill=gray!20] (-1,-0.7) rectangle (1.2,1.3);
\draw [line width =1pt,decoration={markings, mark=at position 0.5 with {\arrow{>}}},postaction={decorate}](-1,1)--(0,0);
\draw [line width =1pt,decoration={markings, mark=at position 0.5 with {\arrow{>}}},postaction={decorate}](-1,0)--(0,0);
\draw [line width =1pt,decoration={markings, mark=at position 0.5 with {\arrow{>}}},postaction={decorate}](-1,-0.4)--(0,0);
\draw [line width =1pt,decoration={markings, mark=at position 0.5 with {\arrow{<}}},postaction={decorate}](1.2,1)  --(0.2,0);
\draw [line width =1pt,decoration={markings, mark=at position 0.5 with {\arrow{<}}},postaction={decorate}](1.2,0)  --(0.2,0);
\draw [line width =1pt,decoration={markings, mark=at position 0.5 with {\arrow{<}}},postaction={decorate}](1.2,-0.4)--(0.2,0);
\node  at(-0.8,0.5) {$\vdots$};
\node  at(1,0.5) {$\vdots$};
\end{tikzpicture}}=(-t^{-1})^{-\frac{n(n-1)}{2}}\cdot \sum_{\sigma\in S_n}
(-t^{\frac{1-n}n})^{\ell(\sigma)} \raisebox{-.30in}{
\begin{tikzpicture}
\filldraw[draw=white,fill=gray!20] (-1,-0.7) rectangle (1.2,1.3);
\draw [line width =1pt,decoration={markings, mark=at position 0.5 with {\arrow{>}}},postaction={decorate}](-1,1)--(0,0);
\draw [line width =1pt,decoration={markings, mark=at position 0.5 with {\arrow{>}}},postaction={decorate}](-1,0)--(0,0);
\draw [line width =1pt,decoration={markings, mark=at position 0.5 with {\arrow{>}}},postaction={decorate}](-1,-0.4)--(0,0);
\draw [line width =1pt,decoration={markings, mark=at position 0.5 with {\arrow{<}}},postaction={decorate}](1.2,1)  --(0.2,0);
\draw [line width =1pt,decoration={markings, mark=at position 0.5 with {\arrow{<}}},postaction={decorate}](1.2,0)  --(0.2,0);
\draw [line width =1pt,decoration={markings, mark=at position 0.5 with {\arrow{<}}},postaction={decorate}](1.2,-0.4)--(0.2,0);
\node  at(-0.8,0.5) {$\vdots$};
\node  at(1,0.5) {$\vdots$};
\filldraw[draw=black,fill=gray!20,line width =1pt]  (0.1,0.3) ellipse (0.4 and 0.7);
\node  at(0.1,0.3){$\sigma_{+}$};
\end{tikzpicture}},
\end{equation}
where the ellipse enclosing $\sigma_+$  is the minimum crossing positive braid representing a permutation $\sigma\in S_n$ and $\ell(\sigma)=\#\{(i,j)\mid 1\leq i<j\leq n,\ \sigma(i)>\sigma(j)\}$ is the length of $\sigma\in S_n$,
\begin{equation}\label{wzh.five}
   \raisebox{-.30in}{
\begin{tikzpicture}
\filldraw[draw=white,fill=gray!20] (-1,-0.7) rectangle (0.2,1.3);
\draw [line width =1pt](-1,1)--(0,0);
\draw [line width =1pt](-1,0)--(0,0);
\draw [line width =1pt](-1,-0.4)--(0,0);
\draw [line width =1.5pt](0.2,1.3)--(0.2,-0.7);
\node  at(-0.8,0.5) {$\vdots$};
\filldraw[fill=white,line width =0.8pt] (-0.5 ,0.5) circle (0.07);
\filldraw[fill=white,line width =0.8pt] (-0.5 ,0) circle (0.07);
\filldraw[fill=white,line width =0.8pt] (-0.5 ,-0.2) circle (0.07);
\end{tikzpicture}}
   =
   \mathbbm{a} \sum_{\sigma \in S_n} (-t^{-1})^{-\ell(\sigma)}\,  \raisebox{-.30in}{
\begin{tikzpicture}
\filldraw[draw=white,fill=gray!20] (-1,-0.7) rectangle (0.2,1.3);
\draw [line width =1pt](-1,1)--(0.2,1);
\draw [line width =1pt](-1,0)--(0.2,0);
\draw [line width =1pt](-1,-0.4)--(0.2,-0.4);
\draw [line width =1.5pt,decoration={markings, mark=at position 1 with {\arrow{>}}},postaction={decorate}](0.2,1.3)--(0.2,-0.7);
\node  at(-0.8,0.5) {$\vdots$};
\filldraw[fill=white,line width =0.8pt] (-0.5 ,1) circle (0.07);
\filldraw[fill=white,line width =0.8pt] (-0.5 ,0) circle (0.07);
\filldraw[fill=white,line width =0.8pt] (-0.5 ,-0.4) circle (0.07);
\node [right] at(0.2,1) {$\sigma(n)$};
\node [right] at(0.2,0) {$\sigma(2)$};
\node [right] at(0.2,-0.4){$\sigma(1)$};
\end{tikzpicture}},
\end{equation}
\begin{equation} \label{wzh.six}
\raisebox{-.20in}{
\begin{tikzpicture}
\filldraw[draw=white,fill=gray!20] (-0.7,-0.7) rectangle (0,0.7);
\draw [line width =1.5pt,decoration={markings, mark=at position 1 with {\arrow{>}}},postaction={decorate}](0,0.7)--(0,-0.7);
\draw [color = black, line width =1pt] (0 ,0.3) arc (90:270:0.5 and 0.3);
\node [right]  at(0,0.3) {$i$};
\node [right] at(0,-0.3){$j$};
\filldraw[fill=white,line width =0.8pt] (-0.5 ,0) circle (0.07);
\end{tikzpicture}}   = \delta_{\bar j,i }\,  \mathbbm{c}_{i} \raisebox{-.20in}{
\begin{tikzpicture}
\filldraw[draw=white,fill=gray!20] (-0.7,-0.7) rectangle (0,0.7);
\draw [line width =1.5pt](0,0.7)--(0,-0.7);
\end{tikzpicture}},
\end{equation}
\begin{equation} \label{wzh.seven}
\raisebox{-.20in}{
\begin{tikzpicture}
\filldraw[draw=white,fill=gray!20] (-0.7,-0.7) rectangle (0,0.7);
\draw [line width =1.5pt](0,0.7)--(0,-0.7);
\draw [color = black, line width =1pt] (-0.7 ,-0.3) arc (-90:90:0.5 and 0.3);
\filldraw[fill=white,line width =0.8pt] (-0.55 ,0.26) circle (0.07);
\end{tikzpicture}}
= \sum_{i=1}^n  (\mathbbm{c}_{\bar i})^{-1}\, \raisebox{-.20in}{
\begin{tikzpicture}
\filldraw[draw=white,fill=gray!20] (-0.7,-0.7) rectangle (0,0.7);
\draw [line width =1.5pt,decoration={markings, mark=at position 1 with {\arrow{>}}},postaction={decorate}](0,0.7)--(0,-0.7);
\draw [line width =1pt](-0.7,0.3)--(0,0.3);
\draw [line width =1pt](-0.7,-0.3)--(0,-0.3);
\filldraw[fill=white,line width =0.8pt] (-0.3 ,0.3) circle (0.07);
\filldraw[fill=black,line width =0.8pt] (-0.3 ,-0.3) circle (0.07);
\node [right]  at(0,0.3) {$i$};
\node [right]  at(0,-0.3) {$\bar{i}$};
\end{tikzpicture}},
\end{equation}
\begin{equation}\label{wzh.eight}
\raisebox{-.20in}{
\begin{tikzpicture}
\filldraw[draw=white,fill=gray!20] (-0,-0.2) rectangle (1, 1.2);
\draw [line width =1.5pt,decoration={markings, mark=at position 1 with {\arrow{>}}},postaction={decorate}](1,1.2)--(1,-0.2);
\draw [line width =1pt](0.6,0.6)--(1,1);
\draw [line width =1pt](0.6,0.4)--(1,0);
\draw[line width =1pt] (0,0)--(0.4,0.4);
\draw[line width =1pt] (0,1)--(0.4,0.6);
\draw[line width =1pt] (0.4,0.6)--(0.6,0.4);
\filldraw[fill=white,line width =0.8pt] (0.2 ,0.2) circle (0.07);
\filldraw[fill=white,line width =0.8pt] (0.2 ,0.8) circle (0.07);
\node [right]  at(1,1) {$i$};
\node [right]  at(1,0) {$j$};
\end{tikzpicture}
} =t^{-\frac{1}{n}}\left(\delta_{{j<i} }(t-t^{-1})\raisebox{-.20in}{
\begin{tikzpicture}
\filldraw[draw=white,fill=gray!20] (-0,-0.2) rectangle (1, 1.2);
\draw [line width =1.5pt,decoration={markings, mark=at position 1 with {\arrow{>}}},postaction={decorate}](1,1.2)--(1,-0.2);
\draw [line width =1pt](0,0.8)--(1,0.8);
\draw [line width =1pt](0,0.2)--(1,0.2);
\filldraw[fill=white,line width =0.8pt] (0.2 ,0.8) circle (0.07);
\filldraw[fill=white,line width =0.8pt] (0.2 ,0.2) circle (0.07);
\node [right]  at(1,0.8) {$i$};
\node [right]  at(1,0.2) {$j$};
\end{tikzpicture}
}+t^{\delta_{i,j}}\raisebox{-.20in}{
\begin{tikzpicture}
\filldraw[draw=white,fill=gray!20] (-0,-0.2) rectangle (1, 1.2);
\draw [line width =1.5pt,decoration={markings, mark=at position 1 with {\arrow{>}}},postaction={decorate}](1,1.2)--(1,-0.2);
\draw [line width =1pt](0,0.8)--(1,0.8);
\draw [line width =1pt](0,0.2)--(1,0.2);
\filldraw[fill=white,line width =0.8pt] (0.2 ,0.8) circle (0.07);
\filldraw[fill=white,line width =0.8pt] (0.2 ,0.2) circle (0.07);
\node [right]  at(1,0.8) {$j$};
\node [right]  at(1,0.2) {$i$};
\end{tikzpicture}
}\right),
\end{equation}
where
$\delta_{j<i}=1$ if $j<i$ and $0$ otherwise, $\delta_{i,j}=1$ if $i=j$ and $0$ otherwise, and small white dots represent an arbitrary orientation of the edges (left-to-right or right-to-left), consistent for the entire equation. The black dot represents the opposite orientation. When a boundary edge of a shaded area is directed, the direction indicates the height order of the endpoints of the diagrams on that directed line, where going along the direction increases the height, and the involved endpoints are consecutive in the height order. The height order outside the drawn part can be arbitrary.

\begin{remark}[Conventions]\label{rem:conventions}
The relations \eqref{w.cross}--\eqref{wzh.eight} and the constants~\eqref{eq-skein-constants} are those of L\^e--Sikora~\cite{LS} and L\^e--Yu~\cite{LY}, whose parameter $q$ is our $t$. The paper~\cite{CHW}, on whose results we build, uses the inverse parameter: its relations are \eqref{w.cross}--\eqref{wzh.eight} with $t$ replaced by $t^{-1}$, its parameter $\omega^{\frac12}$ is the inverse of ours, and its quiver and its form $P$ are the negatives of ours (see~\cite[Rems.~2.3, 2.4]{CHW}). The simultaneous change $\omega^{\frac12}\mapsto\omega^{-\frac12}$, $Q\mapsto-Q$, $P\mapsto-P$ preserves the quantum torus relations, the compatibility relation~\eqref{eq-compatible} and the mutation rules, so the results of~\cite{CHW} hold verbatim in our convention. With our convention, the comparison with the KKOP seed in Section~\ref{sec-cat} holds without any sign or passage to the opposite algebra. We use the letter $t$, rather than $q$, for two reasons: the letter $q$ is used in Section~\ref{sec-KKOP} for the quantum parameter of $U_q'(\slh)$, and Theorem~\ref{thm:main} shows that our $t$ is exactly the parameter $t$ of the $(q,t)$-characters of Section~\ref{sec-KKOP}. All results remain valid for a general commutative ring $R'$ with an invertible element $\omega^{\frac12}$, by applying the functor $-\otimes_{\mathbb{Z}[\omega^{\pm \frac{1}{2}}]}R'$.
\end{remark}

The algebra structure for $\cS(\fS)$ is given by stacking the stated $n$-webs, i.e. for any two stated $n$-webs $\alpha,\alpha'
 \subset\fS\times(-1,1)$, the product $\alpha\alpha'$ is defined by stacking $\alpha$ above $\alpha'$. That is, if $\alpha \subset \fS\times(0,1)$ and $\alpha' \subset \fS\times (-1,0)$, we have $\alpha \alpha' = \alpha \cup \alpha'$.
The \emph{reflection} $\ast$ is the $\Z$-linear anti-involution of $\cS(\fS)$ that sends $\omega^{1/2}\mapsto\omega^{-1/2}$, switches all crossings and reverses the height orders~\cite{LS}.

\begin{figure}[h]
    \centering
\begin{tikzpicture}[scale=1]
\fill[gray!20] (0,0)--(-1.6,1.4)--(1.6,1.4)--cycle;
\draw[line width=1.5pt] (-1.6,1.4)--(0,0)--(1.6,1.4);
\filldraw[fill=white] (0,0) circle (0.07);
\node[below] at (0,0) {$p$};
\draw[line width=1pt,decoration={markings, mark=at position 0.5 with {\arrow{>}}},postaction={decorate}] (-0.8,0.7) .. controls (-0.2,1.1) and (0.2,1.1) .. (0.8,0.7);
\node[left] at (-0.8,0.55) {$j$}; \node[right] at (0.8,0.55) {$i$};
\begin{scope}[xshift=4.5cm]
\fill[gray!20] (0,0)--(-1.6,1.4)--(1.6,1.4)--cycle;
\draw[line width=1.5pt] (-1.6,1.4)--(0,0)--(1.6,1.4);
\filldraw[fill=white] (0,0) circle (0.07);
\node[below] at (0,0) {$p$};
\draw[line width=1pt,decoration={markings, mark=at position 0.5 with {\arrow{<}}},postaction={decorate}] (-0.8,0.7) .. controls (-0.2,1.1) and (0.2,1.1) .. (0.8,0.7);
\node[left] at (-0.8,0.55) {$j$}; \node[right] at (0.8,0.55) {$i$};
\end{scope}
\end{tikzpicture}
    \caption{The left is $C(p)_{ij}$ and the right is $\overline{C}(p)_{ij}$.}\label{Fig;badarc}
\end{figure}
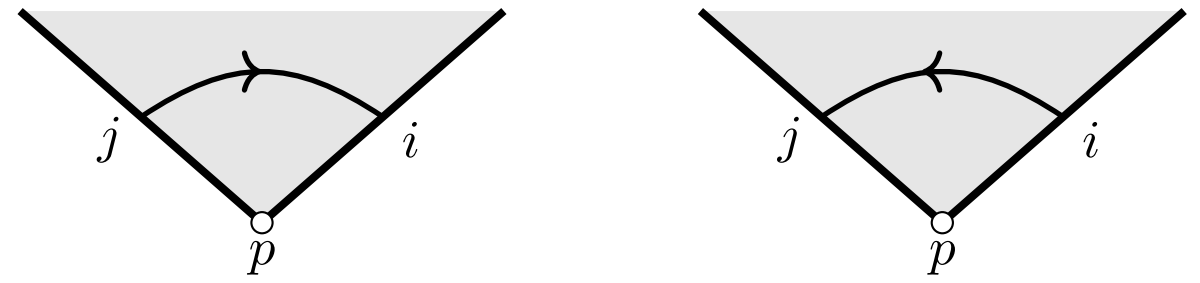

For a boundary puncture $p$ of a pb surface $\fS$,
the corner arcs $C(p)_{ij}$ and $\overline{C}(p)_{ij}$ are the stated arcs depicted in Figure~\ref{Fig;badarc}: they turn around $p$ once, from the state $i$ on one boundary edge adjacent to $p$ to the state $j$ on the other, counterclockwise and clockwise around $p$, respectively.
For a boundary puncture $p$ which is not on a $\mathbb P_1$ component of $\fS$, set
$$C_p=\{C(p)_{ij}\mid i<j\},\quad\overline{C}_p=\{\overline{C}(p)_{ij}\mid i<j\}.$$
Each element of $C_p\cup \overline{C}_p$ is called a \textbf{bad arc} at $p$.

\begin{definition}[{\cite{LY}}]\label{def:reduced}
For a pb surface $\fS$,
\begin{align}\label{eq-def-reduced}
\overline{\cS}(\fS) = \cS(\fS)/I^{\text{bad}}
\end{align}
is called the \textbf{reduced stated $\SL$-skein algebra}, where $I^{\text{bad}}$ is the two-sided ideal of $\cS(\fS)$ generated by all bad arcs. The reflection $\ast$ descends to $\overline{\cS}(\fS)$. We write $\overline\cS_1(\fS)=\overline{\cS}(\fS)\otimes_R\Z$, where $R$ acts on $\Z$ by $\omega^{1/2}\mapsto1$, for the \textbf{classical reduced skein algebra}.
\end{definition}

\subsection{The quantum cluster algebra structure of $\overline{\cS}(\mathbb P_m)$}\label{subsec-CHW}

In this subsection we fix $m\ge3$ and recall the quantum seed of $\overline{\cS}(\bP_m)$ constructed in~\cite{CHW}, building on~\cite{LY}.

\begin{convention}[Small vertices, {\cite[\S2.4, \S5]{CHW}}]\label{conv:CHWlabels}
Fix the star-like triangulation of $\bP_m$ with diagonals $(1,i)$, $3\le i\le m-1$, and triangles
\begin{align}\label{eq-triangles}
\Delta_i=(1,i+1,i+2),\qquad 1\le i\le m-2.
\end{align}
Each $\Delta_i$ is subdivided into $n^2$ small triangles, as in the $n$-triangulations of Fock--Goncharov~\cite{FG06}. Its \emph{small vertices} are the lattice points other than the three corners. We label them $(j,s)^i$, with $0\le s\le j\le n$ and $(j,s)\notin\{(0,0),(n,0),(n,n)\}$. Here $j$ is the lattice distance from the puncture $1$; $s=0$ on the edge $(1,i+1)$; $s=j$ on the edge $(1,i+2)$; and $j=n$ on the edge $(i+1,i+2)$, with $s$ increasing towards $i+2$. In barycentric coordinates with respect to the corners $(1,i+1,i+2)$, the point $(j,s)^i$ is
\begin{align}\label{eq-barycentric}
(a,b,c)=(n-j,\ j-s,\ s).
\end{align}
Small vertices of adjacent triangles on a common diagonal are identified: $(j,j)^i=(j,0)^{i+1}$ for $1\le i\le m-3$. We write every small vertex in \emph{normal form}: with $s\ge1$ whenever possible, i.e.\ we use $(j,j)^{i-1}$ instead of $(j,0)^i$ for $i\ge2$. The set $V$ of all small vertices of $\bP_m$ has $(m-2)\binom n2+(m-1)(n-1)$ elements. The \emph{frozen} vertices are those on boundary edges of $\bP_m$, namely $(j,0)^1$, $(n,s)^i$ and $(j,j)^{m-2}$ for $1\le j,s\le n-1$ and $1\le i\le m-2$; there are $m(n-1)$ of them. All other vertices are \emph{mutable}, and we denote the set of mutable vertices by $V_{\rm mut}$.
\end{convention}

\begin{definition}[{\cite[\S2.4, (28)]{CHW}}]\label{def:CHWquiver}
In each triangle $\Delta_i$, draw an arrow along every small edge. The small edges parallel to a side of $\Delta_i$ are oriented like the counterclockwise orientation of that side, i.e.\ along $i+1\to1$, $i+2\to i+1$ and $1\to i+2$. Explicitly, the arrows are
\begin{align}\label{eq-CHW-arrows}
(j+1,s)^i\to(j,s)^i,\qquad (j,s+1)^i\to(j,s)^i,\qquad (j-1,s-1)^i\to(j,s)^i,
\end{align}
whenever both endpoints are small vertices. An arrow lying on a side of $\Delta_i$ gets weight $\tfrac12$, and all other arrows get weight $1$. We denote by $Q=(Q(x,y))_{x,y\in V}$ the signed adjacency matrix of the resulting weighted quiver, summed over all triangles. This is the quiver of~\cite[\S2.4, (28)]{CHW} with all arrows reversed.
\end{definition}

\begin{example}\label{ex:CHWquiver}
Let $n=3$ and $m=5$, with triangles $\Delta_1=(1,2,3)$, $\Delta_2=(1,3,4)$ and $\Delta_3=(1,4,5)$. Figure~\ref{fig:CHWquiver} shows the resulting quiver; frozen vertices are drawn as squares, and the dashed arrows are half-arrows between frozen vertices.

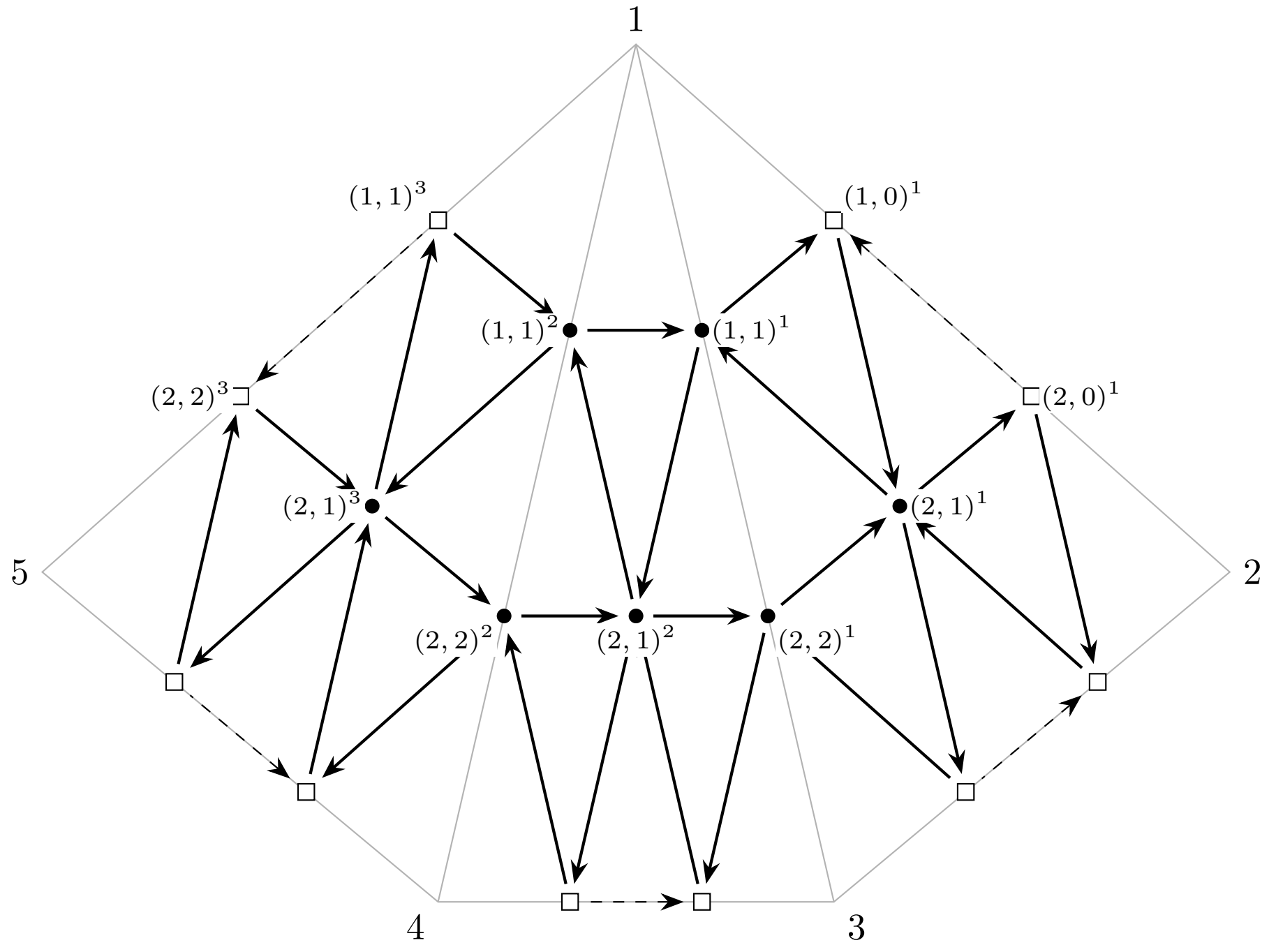
\begin{figure}[ht]
\centering
\begin{tikzpicture}[scale=1.25,
  mut/.style={circle,fill=black,inner sep=1.4pt},
  fro/.style={rectangle,draw=black,fill=white,inner sep=2.2pt},
  ar/.style={{Stealth[length=2.2mm]}-,thick,shorten >=2.5pt,shorten <=2.5pt},
  har/.style={{Stealth[length=2.2mm]}-,dashed,shorten >=2.5pt,shorten <=2.5pt}]
\coordinate (P1) at (0,4); \coordinate (P2) at (4.5,0); \coordinate (P3) at (1.5,-2.5);
\coordinate (P4) at (-1.5,-2.5); \coordinate (P5) at (-4.5,0);
\draw[gray!60] (P1)--(P2)--(P3)--(P4)--(P5)--cycle;
\draw[gray!60] (P1)--(P3); \draw[gray!60] (P1)--(P4);
\node[above] at (P1) {$1$}; \node[right] at (P2) {$2$}; \node[below right] at (P3) {$3$};
\node[below left] at (P4) {$4$}; \node[left] at (P5) {$5$};
\node[fro] (F1) at (1.5,2.667) {}; \node[fro] (F2) at (3,1.333) {};
\node[mut] (d11) at (0.5,1.833) {}; \node[mut] (c1) at (2,0.5) {}; \node[mut] (d12) at (1,-0.333) {};
\node[fro] (b11) at (3.5,-0.833) {}; \node[fro] (b12) at (2.5,-1.667) {};
\node[mut] (d21) at (-0.5,1.833) {}; \node[mut] (c2) at (0,-0.333) {}; \node[mut] (d22) at (-1,-0.333) {};
\node[fro] (b21) at (0.5,-2.5) {}; \node[fro] (b22) at (-0.5,-2.5) {};
\node[fro] (e1) at (-1.5,2.667) {}; \node[mut] (c3) at (-2,0.5) {}; \node[fro] (e2) at (-3,1.333) {};
\node[fro] (b31) at (-2.5,-1.667) {}; \node[fro] (b32) at (-3.5,-0.833) {};
\foreach \A/\B/\C/\D/\E/\F/\G in {F1/F2/d11/c1/d12/b11/b12, d11/d12/d21/c2/d22/b21/b22, d21/d22/e1/c3/e2/b31/b32}{
  \draw[ar] (\A)--(\C); \draw[ar] (\C)--(\D); \draw[ar] (\B)--(\D); \draw[ar] (\D)--(\E);
  \draw[ar] (\D)--(\A); \draw[ar] (\D)--(\F); \draw[ar] (\E)--(\G); \draw[ar] (\F)--(\B); \draw[ar] (\G)--(\D);
  \draw[har] (\F)--(\G);}
\draw[har] (F1)--(F2); \draw[har] (e2)--(e1);
\node[above right,font=\scriptsize,fill=white,inner sep=0.5pt,outer sep=2pt] at (F1) {$(1,0)^1$}; \node[right,font=\scriptsize,fill=white,inner sep=0.5pt,outer sep=2pt] at (F2) {$(2,0)^1$};
\node[right,font=\scriptsize,fill=white,inner sep=0.5pt,outer sep=2pt] at (d11) {$(1,1)^1$}; \node[right,font=\scriptsize,fill=white,inner sep=0.5pt,outer sep=2pt] at (c1) {$(2,1)^1$};
\node[below right,font=\scriptsize,fill=white,inner sep=0.5pt,outer sep=2pt] at (d12) {$(2,2)^1$};
\node[left,font=\scriptsize,fill=white,inner sep=0.5pt,outer sep=2pt] at (d21) {$(1,1)^2$}; \node[below,font=\scriptsize,fill=white,inner sep=0.5pt,outer sep=2pt] at (c2) {$(2,1)^2$};
\node[below left,font=\scriptsize,fill=white,inner sep=0.5pt,outer sep=2pt] at (d22) {$(2,2)^2$};
\node[above left,font=\scriptsize,fill=white,inner sep=0.5pt,outer sep=2pt] at (e1) {$(1,1)^3$}; \node[left,font=\scriptsize,fill=white,inner sep=0.5pt,outer sep=2pt] at (c3) {$(2,1)^3$};
\node[left,font=\scriptsize,fill=white,inner sep=0.5pt,outer sep=2pt] at (e2) {$(2,2)^3$};
\end{tikzpicture}
\caption{The quiver $Q$ for $n=3$ and $m=5$ (Example~\ref{ex:CHWquiver}). The unlabelled squares near the edges $(2,3)$, $(3,4)$, $(4,5)$ are the frozen vertices $(3,1)^i$ and $(3,2)^i$, with $(3,1)^i$ closer to the puncture $i+1$.}\label{fig:CHWquiver}
\end{figure}
\end{example}

The next lemma, whose proof is straightforward, records which arrows lie on the sides of a triangle; this is how half-arrows combine along diagonals. It is used in the proofs of Lemmas~\ref{lem:wtkills} and~\ref{lem:corequiver}.

\begin{lemma}\label{lem:sideedges}
In the triangle $\Delta_i$, an arrow $(j+1,s)\to(j,s)$ lies on a side iff $s=0$, namely the side $(1,i+1)$; an arrow $(j,s+1)\to(j,s)$ lies on a side iff $j=n$, namely the side $(i+1,i+2)$; and an arrow $(j-1,s-1)\to(j,s)$ lies on a side iff $s=j$, namely the side $(1,i+2)$. Consequently $Q(x,y)\in\Z$ whenever $x$ or $y$ is mutable.
\end{lemma}

For a triangle $\Delta$, let $P_{\Delta}$ be the $\Z_3$-invariant skew-symmetric matrix on its small vertices defined as follows; it is the negative of the matrix of~\cite[(30)]{CHW}, and it agrees with~\cite{LY}. If two small vertices $x$ and $x'$ of $\Delta$ have barycentric coordinates $(a,b,c)$ and $(a',b',c')$, and either $a\le a'$, $b\le b'$ or $a\ge a'$, $b\ge b'$, then
\begin{align}\label{eq-PDelta}
P_{\Delta}(x,x')=n(ab'-ba').
\end{align}
For each pair, some cyclic rotation of the coordinates brings it into this position.

For a small vertex $x\in V$ and $1\le i\le m-2$, the \emph{skeleton} $\mathrm{sk}_{\Delta_{i}}(x)$ is a formal $\Z$-linear combination of small vertices of $\Delta_{i}$. It is obtained by elongating the $Y$-shaped graph of $x$ through the neighbouring triangles, turning left~{\cite[(29)]{CHW}, \cite{LY}}. For the star-like triangulation it is explicit. Let $x=(j,s)^l$, written in normal form (Convention~\ref{conv:CHWlabels}), with barycentric coordinates $(a,b,c)$ as in \eqref{eq-barycentric} in $\Delta_l$. Then $\mathrm{sk}_{\Delta_l}(x)$ contains $x$ itself; if $a>0$, then for each $i<l$, $\mathrm{sk}_{\Delta_i}(x)$ contains the vertex with coordinates $(a,0,n-a)$, i.e.\ $(j,j)^i$; and if $c>0$ and $l<m-2$, then $\mathrm{sk}_{\Delta_{l+1}}(x)$ contains the vertex with coordinates $(n-c,c,0)$, i.e.\ $(s,0)^{l+1}$. All these terms have coefficient $1$, and $\mathrm{sk}_{\Delta_i}(x)=0$ for all other $i$.

\begin{definition}[{\cite[(31), (50)]{CHW}, \cite{LY}}]\label{def:Pbar}
For $x,y\in V$, put
\begin{align}\label{eq-P}
P(x,y)=\sum_{i=1}^{m-2}P_{\Delta_{i}}\bigl(\mathrm{sk}_{\Delta_{i}}(x),\mathrm{sk}_{\Delta_{i}}(y)\bigr),
\end{align}
where $P_{\Delta_i}$ is extended bilinearly to formal linear combinations. All entries of $P$ are divisible by $n$, by \eqref{eq-PDelta}, and we put
\begin{align}\label{eq-LambdaCHW}
\Lambda_{\rm skein}=\tfrac1n P .
\end{align}
\end{definition}

\begin{figure}
    \centering
    \includegraphics[width=0.4\linewidth]{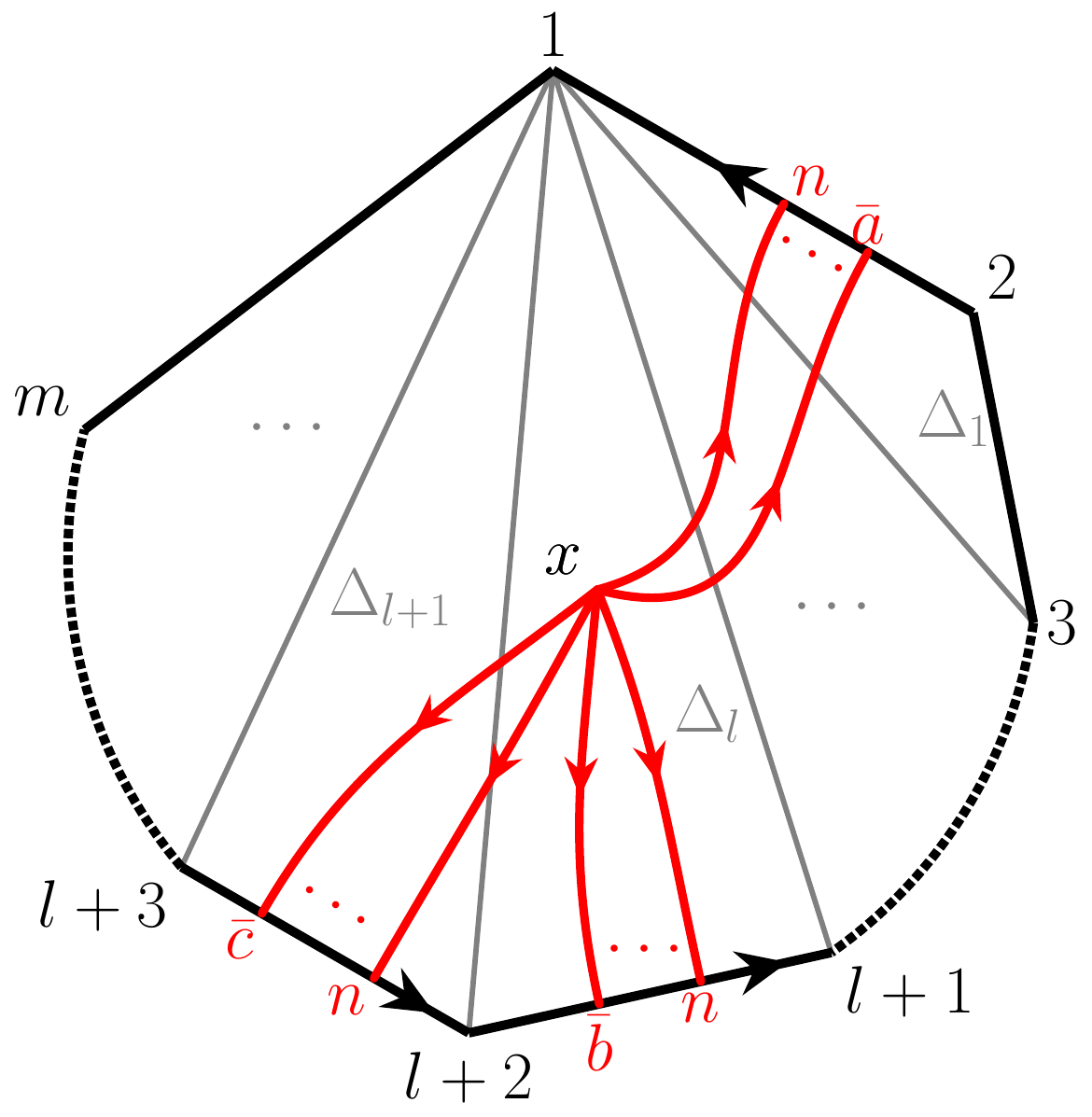}
    \caption{The stated $n$-web for the cluster variable $A_x$.}
    \label{fig:placeholder}
\end{figure}

\begin{proposition}[{\cite[\S2.5]{CHW}, \cite{LY}}]\label{prop:CHWseed}
For each $x\in V$, there exists an element $A_{x}\in\overline{\cS}(\bP_m)$, given by the image of the reflection-invariant stated $n$-web shown in Figure~\ref{fig:placeholder}, such that the following hold:
\begin{enumerate}[label={\rm (\alph*)}]
\item $A_{x}A_{y}=\omega^{P(x,y)}A_{y}A_{x}=\xi^{\Lambda_{\rm skein}(x,y)}A_yA_x$ for all $x,y\in V$;
\item the elements $A_{x}$ generate a quantum torus inside the skew-field of fractions $\mathcal F_{\rm skein}$ of $\overline{\cS}(\bP_m)$, and the Weyl-ordered products $M_{\rm skein}({\bf k})=\bigl[\prod_xA_x^{k_x}\bigr]$ (see \eqref{Weyl-A}) satisfy \eqref{eq-M-relation} with $\Pi=\Lambda_{\rm skein}$;
\item $\sum_{z\in V}Q(z,x)\Lambda_{\rm skein}(z,y)=2n\,\delta_{x,y}$ for all $x\in V_{\rm mut}$ and $y\in V$.
\end{enumerate}
Consequently, $\mathsf s_{\rm skein}=(Q,\Lambda_{\rm skein},M_{\rm skein})$ is a quantum seed in $\mathcal F_{\rm skein}$, compatible with the scalar $2n$.
\end{proposition}

The following theorem from~\cite{CHW} serves as the starting point of this paper. All subsequent results concerning $\overline{\cS}(\bP_m)$ are derived from it via operations on quantum seeds.

\begin{theorem}[{\cite[Thm.~5.1]{CHW}}]\label{thm:CHW}
For every $m\ge3$,
\begin{align}\label{eq-CHW}
\overline{\cS}(\bP_m)=\mathscr A(\mathsf s_{\rm skein})=\mathscr U(\mathsf s_{\rm skein})
\end{align}
as subalgebras of $\mathcal F_{\rm skein}$. In particular, all frozen variables are invertible in $\overline{\cS}(\bP_m)$, and $\overline{\cS}(\bP_m)$ is a domain.
\end{theorem}

\subsection{A quasi-isomorphism for $\overline{\cS}(\mathbb P_m)$}\label{subsec-qc}

In this subsection we define the subalgebra $\mathcal A_{n,m}$ of $\overline{\cS}(\bP_m)$, and an admissible change of frozen coordinates of $\mathsf s_{\rm skein}$ which isolates the frozen vertices on the boundary edges $(1,2),\dots,(m-1,m)$. Both are defined by the edge weights of stated $n$-webs, introduced by L\^e and Yu~\cite[\S4.3]{LY}. The resulting quantum seed $\mathsf s^\qc_{\rm skein}$ (Definition~\ref{def:sqc}) is used in \S\ref{subsec-ident} to identify the skein seed with a KKOP seed, and its cluster algebra decomposes $\overline{\cS}(\bP_m)$ as a skew Laurent extension of $\mathcal A_{n,m}$ (Proposition~\ref{prop:skewLaurent}); this decomposition is used in \S\ref{subsec-CI} to categorify $\overline{\cS}(\bP_m)$, and in Section~\ref{sec-braid} to identify $\mathcal A_{n,m}$ with a quantized braid variety.

Throughout this subsection we abbreviate $\overline{\cS}(\bP_m)$ to $\overline{\cS}$, and we label the boundary edges of $\bP_m$ as
\begin{align}\label{eq-edges}
E_0=(1,2),\qquad E_i=(i+1,i+2)\ \ (1\le i\le m-2),\qquad E_{m-1}=(1,m).
\end{align}

\begin{definition}\label{def:core}
Let $\Viso\subset V$ be the set of small vertices on the edges $E_0,E_1,\dots,E_{m-2}$, i.e.\
\begin{align}\label{eq-I}
\Viso=\{(j,0)^1\mid 1\le j\le n-1\}\ \sqcup\ \{(n,s)^i\mid 1\le s\le n-1,\ 1\le i\le m-2\},
\end{align}
and let $\Vc=V\setminus \Viso=\{(j,s)^i\mid 1\le s\le j\le n-1,\ 1\le i\le m-2\}$, written in normal form (Convention~\ref{conv:CHWlabels}).
\end{definition}

Thus $|\Viso|=(m-1)(n-1)$, all vertices of $\Viso$ are frozen, and the frozen vertices of $\Vc$ are the $n-1$ vertices $(j,j)^{m-2}$ on the edge $E_{m-1}$. The subscript ``iso'' refers to the fact that these vertices become isolated after a change of frozen coordinates (Proposition~\ref{prop:sqc}(a)), and ``core'' to the fact that the vertices of $\Vc$ carry the KKOP seed (Theorem~\ref{thm:seedmatch}).

\emph{Edge weights.} Let $\varpi_i$, $1\leq i\leq n-1$, be the fundamental weights of $\mathfrak{sl}_n$, and let $\mathsf P=\bigoplus_{1\leq i\leq n-1}\Z\varpi_i$ be the weight lattice. Put $\varpi_0=\varpi_n=0$ and $\varepsilon_i=\varpi_i-\varpi_{i-1}$ for $1\le i\le n$, so that $\varepsilon_1+\dots+\varepsilon_n=0$.

\begin{definition}[{\cite[\S4.3]{LY}}]\label{def:edgewt}
Let $\alpha$ be a stated $n$-web in $\bP_m\times(-1,1)$. The weight of an endpoint of $\alpha$ with state $i$ is $\varepsilon_{\bar i}$ if the edge of $\alpha$ at this endpoint points towards the boundary, i.e.\ out of the surface, and $-\varepsilon_{i}$ if it points away from it, i.e.\ into the surface. For a boundary edge $E_i$, the \emph{edge weight} $\operatorname{wt}_{E_i}(\alpha)\in\mathsf P$ is the sum of the weights of the endpoints of $\alpha$ on $E_i$. We put $\operatorname{wt}(\alpha)=(\operatorname{wt}_{E_0}(\alpha),\dots,\operatorname{wt}_{E_{m-1}}(\alpha))\in\mathsf P^{m}$.
\end{definition}

\begin{remark}\label{rem:LYwt}
Definition~\ref{def:edgewt} is the edge weight of L\^e and Yu~\cite[\S4.3]{LY}, whose $\mathsf w_i$ is our $\varepsilon_i$.
\end{remark}

The following lemma is the reduced version of~\cite[Prop.~4.1]{LY}. It follows directly from their argument, since the ideal $I^{\rm bad}$ is generated by bad arcs, each of which is a single stated $n$-web and is therefore homogeneous.

\begin{lemma}[{\cite[Prop.~4.1]{LY}}]\label{lem:edgegrading}
The algebra $\overline{\cS}$ is $\mathsf P^m$-graded, i.e.\
\[
\overline{\cS}=\bigoplus_{\mathbf w\in\mathsf P^m}\overline{\cS}^{\mathbf w},\qquad \overline{\cS}^{\mathbf w}\,\overline{\cS}^{\mathbf w'}\subseteq\overline{\cS}^{\mathbf w+\mathbf w'},
\]
where $\overline{\cS}^{\mathbf w}$ is the $R$-span of the images of the stated $n$-webs $\alpha$ with $\operatorname{wt}(\alpha)=\mathbf w$.
\end{lemma}

Let $\operatorname{res}\colon\mathsf P^m\to\mathsf P^{m-1}$ denote the map that forgets the component supported on the edge $E_{m-1}=(1,m)$. For
$\nu\in\mathsf P^{m-1}$,
\begin{align}\label{eq-coarse-grading}
\overline{\cS}[\nu]=\bigoplus_{\mathbf w\in\mathsf P^m:\ \operatorname{res}(\mathbf w)=\nu}\overline{\cS}^{\mathbf w},\qquad\text{so that}\qquad \overline{\cS}=\bigoplus_{\nu\in\mathsf P^{m-1}}\overline{\cS}[\nu],\quad \overline{\cS}[\nu]\overline{\cS}[\nu']\subseteq\overline{\cS}[\nu+\nu'] .
\end{align}

\begin{definition}\label{def:Anm}
We say that a stated $n$-web $\alpha$ is \emph{balanced} if $\operatorname{wt}_{E_i}(\alpha)=0$ for $0\le i\le m-2$. Define
\begin{align}\label{eq-Anm}
\mathcal A_{n,m}=\overline{\cS}[0]=\operatorname{span}_R\bigl\{\alpha\ \mid\ \alpha\text{ is a balanced stated $n$-web}\bigr\}\subseteq\overline{\cS} .
\end{align}
By \eqref{eq-coarse-grading}, $\mathcal A_{n,m}$ is an $R$-subalgebra of $\overline{\cS}$.
\end{definition}

For example, a web without endpoints on $E_0,\dots,E_{m-2}$ lies in $\mathcal A_{n,m}$, and so does a web whose endpoints on each of these edges come in pairs of opposite orientations with states $i$ and $\bar i$.

\emph{Edge weights of the initial variables.} Let $x=(j,s)^l\in V$ be in normal form (Convention~\ref{conv:CHWlabels}), with barycentric coordinates $(a,b,c)$ as in~\eqref{eq-barycentric}. The web representing $A_x$ in Figure~\ref{fig:placeholder} features a single $n$-valent source at $x$ and three legs of weights $a$, $b$, and $c$. A leg of weight $w$ consists of $w$ parallel strands oriented towards the boundary, with endpoints carrying states $\bar w,\bar w+1,\dots,n$, where the strand with state $n$ is the innermost. Specifically, the leg of weight $a$ turns left around the puncture $1$, passes through $\Delta_{l-1},\dots,\Delta_1$, and terminates on $E_0$. The leg of weight $b$ terminates on $E_l$. The leg of weight $c$ crosses the diagonal $(1,l+2)$ into $\Delta_{l+1}$ and turns left around the puncture $l+2$ before terminating on $E_{l+1}$ (provided $l<m-2$); if $l=m-2$, it terminates on $E_{m-1}$. Legs of weight $0$ are absent. By Definition~\ref{def:edgewt}, a leg of weight $w$ contributes
\begin{align}\label{eq-f}
f(w):=\sum_{k=n-w+1}^{n}\varepsilon_{\bar k}=\sum_{i=1}^{w}\varepsilon_i=\varpi_w
\end{align}
to the edge weight of its terminal edge. Note that $f(0)=f(n)=0$, and $f(1),\dots,f(n-1)$ form a $\Z$-basis of $\mathsf P$. For $0\le i\le m-1$ and $w\in\{0,\dots,n\}$, let $f_{E_i}(w)\in\mathsf P^m$ denote the vector with $f(w)$ in the $E_i$-component and $0$ elsewhere. We summarize this as follows.

\begin{lemma}\label{lem:legwt}
For each $x=(j,s)^l\in V$, the element $A_x$ lies in $\overline{\cS}^{\operatorname{wt}(A_x)}$, where
\begin{align}\label{eq-wtAx}
\operatorname{wt}(A_x)=f_{E_0}(a)+f_{E_l}(b)+f_{E_{l+1}}(c),\qquad (a,b,c)=(n-j,j-s,s),
\end{align}
with the convention that $E_{l+1}=E_{m-1}$ when $l=m-2$. This weight is independent of the choice of representative for a diagonal vertex, i.e., $(j,j)^l=(j,0)^{l+1}$. Furthermore, if $g$ is a frozen vertex, then $A_g^{-1}\in\overline{\cS}^{-\operatorname{wt}(A_g)}$.
\end{lemma}

\begin{proof}
The first statement is the computation above. For a diagonal vertex both representatives give the leg of weight $n-j$ on $E_0$ and the leg of weight $j$ on $E_{l+1}$. For the last statement, write $A_g^{-1}=\sum_{\mathbf w}h_{\mathbf w}$ with $h_{\mathbf w}\in\overline{\cS}^{\mathbf w}$, which is possible by Theorem~\ref{thm:CHW} and Lemma~\ref{lem:edgegrading}. Then $1=\sum_{\mathbf w}A_gh_{\mathbf w}$ with $A_gh_{\mathbf w}\in\overline{\cS}^{\operatorname{wt}(A_g)+\mathbf w}$, and $A_gh_{\mathbf w}\ne0$ if $h_{\mathbf w}\neq0$, since $\overline{\cS}$ is a domain. Comparing components, $h_{\mathbf w}=0$ unless $\mathbf w=-\operatorname{wt}(A_g)$.
\end{proof}

Define the homomorphism
\begin{align}\label{eq-wt-hom}
\varrho\colon\Z^V\to\mathsf P^m,\qquad {\bf k}\mapsto\sum_{x\in V}k_x\operatorname{wt}(A_x).
\end{align}
By Lemma~\ref{lem:legwt}, $M_{\rm skein}({\bf k})\in\overline{\cS}^{\varrho({\bf k})}$ whenever $k_x\ge0$ for all mutable $x$.

\begin{lemma}\label{lem:wtkills}
For every mutable vertex $y$, $\varrho({\bf q}_y)=\sum_{x\in V}Q(x,y)\operatorname{wt}(A_x)=0$, where ${\bf q}_y$ denotes the column vector of $Q$ at $y$ defined in~\eqref{eq-column-vector}.
\end{lemma}

\begin{proof}
Assume that $y$ is a vertex in $\Delta_i$. By \eqref{eq-wtAx}, a small vertex of $\Delta_i$ with coordinates $(a,b,c)$ has weight $F_1(a)+F_2(b)+F_3(c)$, where $F_1=f_{E_0}$, $F_2=f_{E_i}$, and $F_3=f_{E_{i+1}}$ (with the convention that $E_{i+1}=E_{m-1}$ if $i=m-2$). Note that these functions vanish at $0$ and $n$.

Suppose first that $y=(a,b,c)$ lies in the interior of $\Delta_i$. By \eqref{eq-CHW-arrows}, the six neighbours of $y$ are $(a+1,b-1,c)$, $(a,b+1,c-1)$, and $(a-1,b,c+1)$, with arrows directed away from $y$; and $(a-1,b+1,c)$, $(a,b-1,c+1)$, and $(a+1,b,c-1)$, with arrows directed towards $y$. All these neighbours have weight $1$. The column sum is therefore given by
\begin{multline*}
\bigl[F_1(a{+}1)+F_2(b{-}1)+F_3(c)\bigr]+\bigl[F_1(a)+F_2(b{+}1)+F_3(c{-}1)\bigr]+\bigl[F_1(a{-}1)+F_2(b)+F_3(c{+}1)\bigr]\\
-\bigl[F_1(a{-}1)+F_2(b{+}1)+F_3(c)\bigr]-\bigl[F_1(a)+F_2(b{-}1)+F_3(c{+}1)\bigr]-\bigl[F_1(a{+}1)+F_2(b)+F_3(c{-}1)\bigr]=0 .
\end{multline*}

Next, consider the case where $y$ lies on the diagonal $(1,i+2)$, with coordinates $(a,0,c)$ in $\Delta_i$ and $(a,c,0)$ in $\Delta_{i+1}$. The two half-arrows along the diagonal cancel each other out (see Lemma~\ref{lem:sideedges}). The remaining neighbours are $(a,1,c-1)$ (pointing away from $y$) and $(a-1,1,c)$ (pointing towards $y$) in $\Delta_i$, together with $(a-1,c,1)$ (away) and $(a,c-1,1)$ (towards) in $\Delta_{i+1}$. Their respective weights are $F_1(a)+F_2(1)+F_3(c-1)$, $F_1(a-1)+F_2(1)+F_3(c)$, $F_1(a-1)+F_3(c)+F'(1)$, and $F_1(a)+F_3(c-1)+F'(1)$, where $F'=f_{E_{i+2}}$ (or $f_{E_{m-1}}$). Here, we use the fact that the legs of a diagonal vertex terminate on the same edges in both triangles. The signed sum of these weights is again $0$.
\end{proof}

\begin{lemma}\label{lem:theta}
Let $\theta\colon\Z^{\Viso}\to\mathsf P^{m-1}$ be the homomorphism defined by $\theta({\bf e}_g)=\operatorname{res}(\operatorname{wt}(A_g))$ for $g\in \Viso$. Then $\theta$ is an isomorphism.
\end{lemma}

\begin{proof}
By \eqref{eq-wtAx}, the restricted weights $\operatorname{res}(\operatorname{wt}(A_g))$ are given as follows: for $g=(n,s)^{m-2}$, the weight is $f_{E_{m-2}}(n-s)$ (since the leg of weight $s$ terminates on $E_{m-1}$); for $g=(n,s)^i$ with $i<m-2$, it is $f_{E_i}(n-s)+f_{E_{i+1}}(s)$; and for $g=(j,0)^1$, it is $f_{E_0}(n-j)+f_{E_1}(j)$. Ordering the edges as $E_{m-2},E_{m-3},\dots,E_0$ and the vertices of $\Viso$ accordingly, the matrix of $\theta$ becomes block lower triangular. Its diagonal blocks map the $n-1$ vertices attached to $E_i$ bijectively onto the basis $\{f(w)\}_{1\le w\le n-1}$ of $\mathsf P$.
\end{proof}

We now define the \emph{balanced $\Viso$-degree}
\begin{align}\label{eq-pinatural}
\deg_{\rm iso}=\theta^{-1}\circ\operatorname{res}\circ\varrho\colon\Z^V\longrightarrow\Z^{\Viso} ,
\end{align}
where $\varrho$ is the homomorphism in \eqref{eq-wt-hom}. By Lemma~\ref{lem:wtkills}, $\deg_{\rm iso}({\bf q}_y)=0$ for every mutable $y$, while by the definition of $\theta$, $\deg_{\rm iso}({\bf c})={\bf c}$ for any ${\bf c}\in\Z^{\Viso}\subseteq\Z^V$. Furthermore, Lemma~\ref{lem:legwt} implies that $M_{\rm skein}({\bf k})\in\overline{\cS}[\theta(\deg_{\rm iso}({\bf k}))]$ whenever $k_x\ge0$ for all mutable $x$; in particular, $M_{\rm skein}({\bf k})\in\mathcal A_{n,m}$ if $\deg_{\rm iso}({\bf k})=0$.

\begin{lemma}\label{lem:balP}
Let $\mathsf G_{\rm iso}$ be the $\Z$-linear endomorphism of $\Z^V$ defined by
\begin{align}\label{eq-Pnatural}
\mathsf G_{\rm iso}({\bf e}_x)={\bf e}_x-\deg_{\rm iso}({\bf e}_x)\quad(x\in\Vc),\qquad \mathsf G_{\rm iso}({\bf e}_g)={\bf e}_g\quad(g\in \Viso),
\end{align}
where $\deg_{\rm iso}({\bf e}_x)$ is viewed as an element of $\Z^{\Viso}\subseteq\Z^V$. Then $\mathsf G_{\rm iso}$ is admissible for $\mathsf s_{\rm skein}$ (Definition~\ref{def:admissible}), and satisfies
\begin{align}\label{eq-Pnatural-degree}
\deg_{\rm iso}\bigl(\mathsf G_{\rm iso}({\bf k})\bigr)={\bf k}|_{\Viso}\qquad({\bf k}\in\Z^V).
\end{align}
Moreover, for every mutable $y$, the column vector $(\mathsf G_{\rm iso})^{-1}({\bf q}_y)$ has the same $\Vc$-coordinates as ${\bf q}_y$ and vanishing $\Viso$-coordinates.
\end{lemma}

\begin{proof}
For $x\in\Vc$, $\mathsf G_{\rm iso}({\bf e}_x)-{\bf e}_x\in\Z^{\Viso}$, and all vertices of $\Viso$ are frozen. In the basis ordered as $\Vc\sqcup \Viso$, $\mathsf G_{\rm iso}$ is unitriangular, so it is invertible and $\mathsf G_{\rm iso}$ restricts to an automorphism of $\Z^{V\setminus V_{\rm mut}}$. Hence $\mathsf G_{\rm iso}$ is admissible. For $x\in\Vc$, $\deg_{\rm iso}(\mathsf G_{\rm iso}({\bf e}_x))=\deg_{\rm iso}({\bf e}_x)-\deg_{\rm iso}({\bf e}_x)=0$, and for $g\in \Viso$, $\deg_{\rm iso}(\mathsf G_{\rm iso}({\bf e}_g))={\bf e}_g$; this gives \eqref{eq-Pnatural-degree}. Let $\hat{\bf q}=(\mathsf G_{\rm iso})^{-1}({\bf q}_y)$. Then $\hat{\bf q}|_{\Viso}=\deg_{\rm iso}(\mathsf G_{\rm iso}(\hat{\bf q}))=\deg_{\rm iso}({\bf q}_y)=0$. Since $\mathsf G_{\rm iso}$ and $(\mathsf G_{\rm iso})^{-1}$ are the identity modulo $\Z^{\Viso}$, $\hat{\bf q}$ and ${\bf q}_y$ have the same $\Vc$-coordinates.
\end{proof}

The seed $\mathsf s_{\rm skein}$ cannot be compared directly with a KKOP seed: it has $m(n-1)$ frozen vertices instead of $n-1$, and the frozen vertices of $\Viso$ are joined by arrows to mutable vertices. We therefore change the frozen coordinates by $\mathsf G_{\rm iso}$. The choice of $\mathsf G_{\rm iso}$ is dictated by two requirements. First, the vertices of $\Viso$ should become isolated (Lemma~\ref{lem:balP}), so that the vertices of $\Vc$ carry a sub-seed. Isolation alone determines $\mathsf G_{\rm iso}({\bf e}_x)-{\bf e}_x$, $x\in\Vc$, only up to a homomorphism $\Z^{\Vc}\to\Z^{\Viso}$ vanishing on the span of the $\Vc$-components of the column vectors of $Q$, i.e.\ up to a homomorphism defined on $\Z^{\Lst}$ (Lemma~\ref{lem:Lst}). Second, the new cluster variables $A^\qc_x$, $x\in\Vc$, should be balanced, i.e.\ lie in $\mathcal A_{n,m}$; this fixes the remaining freedom. The first requirement makes the exchange matrices of the sub-seed on $\Vc$ and of the KKOP seed agree (Lemma~\ref{lem:M}). The second requirement makes their quantum forms agree as well (Lemma~\ref{lem:XiS} and Theorem~\ref{thm:seedmatch}), and it gives the description $\mathcal A_{n,m}=\mathscr A(\mathsf s^\qc_\core)$ (Proposition~\ref{prop:skewLaurent}). The seed $\mathsf s^\qc_{\rm skein}$ is used in \S\ref{subsec-ident} to identify the skein seed with the KKOP seed $\mathsf s(u^{m-2})$, in \S\ref{subsec-CI}, where its quantum form defines the twist $\Xi$ of Theorem~\ref{thm:main}, and in Section~\ref{sec-braid}.

\begin{definition}\label{def:sqc}
We define the quantum seed
\begin{align}\label{eq-sqc}
\mathsf s^\qc_{\rm skein}=(Q^\qc,\Lambda^\qc,M^\qc):=(\mathsf s_{\rm skein})^{\mathsf G_{\rm iso}}
\end{align}
of Definition~\ref{def:admissible}, i.e.\ $M^\qc({\bf k})=M_{\rm skein}(\mathsf G_{\rm iso}({\bf k}))$ and $\Lambda^\qc({\bf k},{\bf t})=\Lambda_{\rm skein}(\mathsf G_{\rm iso}({\bf k}),\mathsf G_{\rm iso}({\bf t}))$. We write $A^\qc_x=M^\qc({\bf e}_x)$ for $x\in V$; thus $A^\qc_g=A_g$ for $g\in \Viso$. Let $\mathsf s^\qc_\core$ be the sub-seed of $\mathsf s^\qc_{\rm skein}$ supported on $\Vc$, with the vertices $(j,j)^{m-2}$, $1\le j\le n-1$, frozen (Definition~\ref{def:subseed}).
\end{definition}

For $\mathbf c\in\Z^{\Viso}$ we write
\begin{align}\label{eq-Ac}
A^{\mathbf c}=M^\qc(\mathbf c)=M_{\rm skein}(\mathbf c)
\end{align}
for the Weyl-ordered monomial in the frozen variables $A_g$, $g\in \Viso$.

\begin{proposition}\label{prop:sqc}
The seed $\mathsf s^\qc_{\rm skein}$ is a quantum seed compatible with scalar $2n$, and the following hold.
\begin{enumerate}[label={\rm (\alph*)}]
\item $Q^\qc(g,y)=0$ for $g\in \Viso$ and $y\in V_{\rm mut}$, i.e.\ all vertices of $\Viso$ are isolated, and $Q^\qc(x,y)=Q(x,y)$ for $x\in\Vc$ and $y\in V_{\rm mut}$.
\item $\overline{\cS}(\bP_m)=\mathscr A(\mathsf s^\qc_{\rm skein})=\mathscr U(\mathsf s^\qc_{\rm skein})$, and the identity map is a quantum quasi-isomorphism $\mathscr A(\mathsf s^\qc_{\rm skein})\to\mathscr A(\mathsf s_{\rm skein})$ (Definition~\ref{def:QQH}).
\item $\mathsf s^\qc_\core$ is a quantum seed compatible with scalar $2n$, and for every sequence of mutable vertices, the cluster variables of the mutated seed of $\mathsf s^\qc_\core$ are those of the mutated seed of $\mathsf s^\qc_{\rm skein}$ at the vertices of $\Vc$.
\item Each $A^\qc_x$, $x\in\Vc$, is the unique element of $\mathcal A_{n,m}$ of the form $M_{\rm skein}({\bf e}_x+{\bf c})$, ${\bf c}\in\Z^{\Viso}$. In other words, $A^\qc_x$ is $A_x$ multiplied by the unique Laurent monomial in the $A_g$, $g\in \Viso$, which makes it balanced, up to the power of $\xi^{1/2}$ given by the Weyl normalization.
\end{enumerate}
\end{proposition}

\begin{proof}
By Lemma~\ref{lem:Pchange} and Proposition~\ref{prop:CHWseed}, $\mathsf s^\qc_{\rm skein}$ is a quantum seed compatible with scalar $2n$. Part (a) is Lemma~\ref{lem:balP}, and (b) follows from Theorem~\ref{thm:CHW} and Corollary~\ref{cor:Pchange}. By (a), there are no arrows between the mutable vertices and the vertices of $\Viso$, so (c) is Lemma~\ref{lem:subseed}. For (d), $\deg_{\rm iso}({\bf e}_x+{\bf c})=\deg_{\rm iso}({\bf e}_x)+{\bf c}$, which vanishes only for ${\bf c}=-\deg_{\rm iso}({\bf e}_x)$, and $M_{\rm skein}({\bf e}_x+{\bf c})\in\overline{\cS}[\theta(\deg_{\rm iso}({\bf e}_x+{\bf c}))]$.
\end{proof}

The next proposition decomposes the skein algebra along the frozen variables $A_g$, $g\in \Viso$. It is used in the proof of Theorem~\ref{thm:main} and in Theorem~\ref{thm:Aquant}.

\begin{proposition}\label{prop:skewLaurent}
We have $\mathcal{A}_{n,m}=\mathscr{A}(\mathsf{s}^\mathrm{qc}_\mathrm{core})$, and
\begin{align}\label{eq-skewLaurent}
\overline{\cS}(\bP_m)=\bigoplus_{\mathbf{c}\in\Z^{\Viso}}\mathcal{A}_{n,m}\,A^{\mathbf{c}},
\end{align}
where each summand $\mathcal{A}_{n,m}A^{\mathbf{c}}=\overline{\cS}[\theta(\mathbf{c})]$ is a free left $\mathcal{A}_{n,m}$-module of rank one. Moreover, $A_g\mathcal{A}_{n,m}=\mathcal{A}_{n,m}A_g$ for all $g\in \Viso$.
\end{proposition}

\begin{proof}
We first show that $\mathscr{A}(\mathsf{s}^\mathrm{qc}_\mathrm{core})\subseteq\overline{\cS}[0]=\mathcal{A}_{n,m}$. By Proposition~\ref{prop:sqc}(c), the algebra $\mathscr{A}(\mathsf{s}^\mathrm{qc}_\mathrm{core})$ is generated by the cluster variables at the vertices of $\Vc$ from the seeds obtained by mutating $\mathsf{s}^\qc_{\rm skein}$ at mutable vertices, together with the inverses of the frozen variables among them. Since $\overline{\cS}[0]$ is a subalgebra, it suffices to show that all these generators lie in $\overline{\cS}[0]$. We proceed by induction on the length of the mutation sequence.

For the initial seed, $A^\qc_x\in\overline{\cS}[0]$ for every $x\in\Vc$ by Proposition~\ref{prop:sqc}(d). Furthermore, if $x$ is frozen, then $A^\qc_x$ is, up to a power of $\xi^{1/2}$, a Laurent monomial in the frozen variables of $\mathsf{s}_{\rm skein}$. Consequently, it is invertible in $\overline{\cS}$, and $(A^\qc_x)^{-1}\in\overline{\cS}[0]$ by Lemma~\ref{lem:legwt} and~\eqref{eq-coarse-grading}. Since frozen variables remain unchanged under mutation, it remains to consider the newly generated exchangeable cluster variables.

Let $\mathsf{s}'$ be a seed obtained from $\mathsf{s}^\qc_{\rm skein}$ by a sequence of mutations, with toric frame $M'$, and assume that its cluster variables at the vertices of $\Vc$ lie in $\overline{\cS}[0]$. Let $A'_k$ be the cluster variable at $k$ of $\mu_k(\mathsf{s}')$, where $k$ is mutable. By~\eqref{eq-mutation-Q}, an isolated vertex remains isolated after mutation; thus, by Proposition~\ref{prop:sqc}(a), the column vector ${\bf q}_k$ of the exchange matrix of $\mathsf{s}'$ is supported on $\Vc$. Let ${\bf b}^\pm=[\pm{\bf q}_k]_+\in\Z^{\Vc}_{\ge0}$. Multiplying~\eqref{eq-exchange} on the left by $A_k=M'({\bf e}_k)$ and applying~\eqref{eq-M-relation}, we obtain
\[
A_kA'_k=\xi^{a_+}M'({\bf b}^+)+\xi^{a_-}M'({\bf b}^-)
\]
for some $a_\pm\in\frac{1}{2}\Z$. Each $M'({\bf b}^\pm)$ is, up to a power of $\xi^{1/2}$, a product of cluster variables of $\mathsf{s}'$ at the vertices of $\Vc$, implying that $A_kA'_k\in\overline{\cS}[0]$. Since $A'_k\in\mathscr{A}(\mathsf{s}^\qc_{\rm skein})=\overline{\cS}$ by Proposition~\ref{prop:sqc}(b), we can write $A'_k=\sum_\nu h_\nu$ with $h_\nu\in\overline{\cS}[\nu]$. By~\eqref{eq-coarse-grading}, $A_kh_\nu\in\overline{\cS}[\nu]$, and $A_kh_\nu\neq 0$ whenever $h_\nu\neq 0$, because $\overline{\cS}$ is a domain (Theorem~\ref{thm:CHW}). Comparing the components of $A_kA'_k$ in~\eqref{eq-coarse-grading}, we deduce that $h_\nu=0$ for all $\nu\neq 0$, i.e., $A'_k\in\overline{\cS}[0]$.

Next, we observe that $\overline{\cS}=\mathscr{A}(\mathsf{s}^\qc_{\rm skein})$ is generated by the cluster variables of the mutated seeds at the vertices of $\Vc$, the inverses of the frozen variables, and the elements $A_g^{\pm1}$ for $g\in \Viso$. Since each $A_g$ is a cluster variable in every seed of the mutation class, it quasi-commutes with all cluster variables (Lemma~\ref{lem:normal}). Hence,
\[
\overline{\cS}=\sum_{\mathbf{c}\in\Z^{\Viso}}\mathscr{A}(\mathsf{s}^\mathrm{qc}_\mathrm{core})\,A^{\mathbf{c}}.
\]
By Lemma~\ref{lem:legwt}, $A^{\mathbf{c}}\in\overline{\cS}[\theta(\mathbf{c})]$, which implies $\mathscr{A}(\mathsf{s}^\mathrm{qc}_\mathrm{core})A^{\mathbf{c}}\subseteq\overline{\cS}[\theta(\mathbf{c})]$. Given that $\theta$ is bijective (Lemma~\ref{lem:theta}) and $\overline{\cS}=\bigoplus_\nu\overline{\cS}[\nu]$, we conclude that $\overline{\cS}[\theta(\mathbf{c})]=\mathscr{A}(\mathsf{s}^\mathrm{qc}_\mathrm{core})A^{\mathbf{c}}$ for every $\mathbf{c}$. Setting $\mathbf{c}=0$ yields $\mathcal{A}_{n,m}=\overline{\cS}[0]=\mathscr{A}(\mathsf{s}^\mathrm{qc}_\mathrm{core})$, which establishes~\eqref{eq-skewLaurent}. Since right multiplication by the invertible element $A^{\mathbf{c}}$ is injective, each summand is free of rank one. Finally, by~\eqref{eq-coarse-grading}, we have $A_g\overline{\cS}[0]A_g^{-1}\subseteq\overline{\cS}[0]$ and $A_g^{-1}\overline{\cS}[0]A_g\subseteq\overline{\cS}[0]$.
\end{proof}

Since the elements $A_g$ for $g \in \Viso$ quasi-commute with one another, $\overline{\cS}(\bP_m)$ forms an iterated skew Laurent extension of $\mathcal{A}_{n,m}$; that is, it constitutes a ``quantum torus over $\mathcal{A}_{n,m}$'' in the frozen variables $A_g$ for $g \in \Viso$. Throughout the sequel, this is precisely what we mean by a skew Laurent extension.

\section{Category $\mathcal C({\mathsf b})$ and its localization}\label{sec-KKOP}

Qin \cite{Qin24a} and Kashiwara, Kim, Oh, and Park \cite{KKOP3} independently attached to each positive braid $\mathsf b\in\mathrm{Br}_n^+$ a monoidal subcategory $\cC(\mathsf b)$ of the Hernandez--Leclerc category of $U_q'(\slh)$-modules, and proved that its quantum Grothendieck ring possesses a quantum cluster algebra structure, wherein the cluster monomials correspond to the $(q,t)$-characters of real simple modules~\cite{Qin24a}, \cite{Qin24b}, \cite{KKOP3}. In this section, we review their quantum seed $\mathsf{s}(\mathbf{i})$ associated with a word $\mathbf{i}$ of $\mathsf{b}$ in Section \ref{subsec-KKOPseed}, along with their monoidal categorification theorem in Section \ref{subsec-moncat}. Rather than reproducing the (lengthy) constructions, we precisely state the objects and results required for our arguments, providing relevant references. In Section \ref{subsec-braider}, we construct braiders for the frozen modules within a full subcategory $\mathcal{C}^{\mathrm{s}}(\mathsf{b})\subseteq\mathcal{C}(\mathsf{b})$ sharing the same simple objects, as well as the corresponding monoidal localization $\widetilde{\mathcal{C}}(\mathsf{b})$. The quantum Grothendieck ring of this localization is the quantum cluster algebra $\mathscr{A}(\mathsf{s}(\mathbf{i}))$ with invertible frozen variables (Theorem~\ref{thm:Ktloc}). The results in this section hold for an arbitrary positive braid $\mathsf{b}\in\mathrm{Br}_n^+$, given by a word in which every colour occurs. In Section~\ref{sec-cat}, these results are applied to $\mathsf{b}=\Delta^{m-2}$, whose seed is identified with the skein seed $\mathsf{s}^\mathrm{qc}_\mathrm{core}$ defined in Definition~\ref{def:sqc}.

We work with the untwisted quantum affine algebra $U_q'(\slh)$ of type $A^{(1)}_{n-1}$. The finite type is $\mathsf g=A_{n-1}$ with index set $\mathsf I=\{1,\dots,n-1\}$, all symmetrizers equal $1$, and the Dynkin involution is $\imath^*=n-\imath$.

\subsection{KKOP seed ${\mathsf s}(\mathbf i)$ associated with a positive braid word}\label{subsec-KKOPseed}

\subsubsection{Quantum affine algebras and the Hernandez--Leclerc category}\label{subsec:HL}

Let $\mathsf W=S_n$ be the Weyl group of $\mathfrak{sl}_n$, generated by the simple reflections $s_\imath$, $\imath\in\mathsf I$. Let $\alpha_\imath$ be the simple roots and $\varpi_\imath$ the fundamental weights of $\mathfrak{sl}_n$, and let $\mathsf P=\bigoplus_{\imath\in\mathsf I}\Z\varpi_\imath$ be the weight lattice and $\mathsf Q=\bigoplus_{\imath\in\mathsf I}\Z\alpha_\imath$ the root lattice. We write $(\ ,\ )$ for the $\mathsf W$-invariant symmetric bilinear form on $\mathsf P\otimes\Q$ with $(\alpha_\imath,\alpha_\jmath)=a_{\imath\jmath}$, where $(a_{\imath\jmath})$ is the Cartan matrix of $A_{n-1}$. Then $(\varpi_\imath,\alpha_\jmath)=\delta_{\imath\jmath}$, $(\varpi_\imath,\varpi_\jmath)=\min(\imath,\jmath)-\frac{\imath\jmath}{n}$, and, with $\varepsilon_i=\varpi_i-\varpi_{i-1}$ as in \S\ref{subsec-qc}, $(\varepsilon_a,\varepsilon_b)=\delta_{ab}-\frac1n$.

Let $q$ be an indeterminate, and let $\mathbf k$ be the algebraic closure of $\C(q)$. For a ring $A$ we write $A^\times$ for its group of units; thus $\mathbf k[z^{\pm1}]^\times=\{c\,z^m\mid c\in\mathbf k^\times,\ m\in\Z\}$. Let $U_q'(\slh)$ be the quantum affine algebra without derivation~\cite{Dri85,Jim85}, \cite[\S12]{CP-book}, \cite[\S1.2]{KKOP3}, i.e.\ the $\mathbf k$-algebra generated by $e_i,f_i,t_i^{\pm1}$, $i\in\{0,1,\dots,n-1\}$, subject to the relations of the Drinfeld--Jimbo quantum group of the affine Cartan matrix of type $A^{(1)}_{n-1}$. Let $\cC_{\mathfrak g}$ be the category of finite-dimensional integrable $U_q'(\slh)$-modules~\cite[\S1.2]{KKOP3}, \cite{Kas02}. It is a rigid monoidal abelian category, not braided in general, with exact tensor product, and its simple objects are parametrized by Drinfeld polynomials~\cite{CP91,CP95}. We write $\mathscr D$ for the right duality functor.

For a module $M$, its \emph{affinization} is $M_z=\mathbf k[z^{\pm1}]\otimes M$, with $e_i(u_z)=z^{\delta_{i,0}}(e_iu)_z$, $f_i(u_z)=z^{-\delta_{i,0}}(f_iu)_z$ and $t_i(u_z)=(t_iu)_z$, where $u_z=1\otimes u$. We have $(M\otimes N)_z\cong M_z\otimes_{\mathbf k[z^{\pm1}]}N_z$. For $x\in\mathbf k^\times$, the \emph{spectral shift} of $M$ by $x$ is $M_x=M_z/(z-x)M_z$. For $\imath\in\mathsf I$, let $V(\varpi_\imath)$ be the $\imath$-th fundamental module~\cite[\S5.2]{Kas02}.

A \emph{$Q$-datum}~\cite{FO21}, \cite[Def.~2.2]{KKOP3} for $U_q'(\slh)$ is a height function $\xi_\cQ\colon\mathsf I\to\Z$ with $|\xi_\cQ(\imath)-\xi_\cQ(\jmath)|=1$ whenever $|\imath-\jmath|=1$; we denote it by $\cQ$. Its repetition quiver has vertex set $\widehat\triangle_0=\{(\imath,p)\in\mathsf I\times\Z\mid p-\xi_\cQ(\imath)\in2\Z\}$~\cite[\S2.1]{KKOP3}, \cite{HL15,FO21}. The Hernandez--Leclerc category $\cC^0_{\mathfrak g}$~\cite{HL10} is the smallest full subcategory of $\cC_{\mathfrak g}$ containing the modules $L(\imath,p)=V(\varpi_\imath)_{(-q)^p}$, $(\imath,p)\in\widehat\triangle_0$, and stable under tensor products, subquotients and extensions.

Every finite-dimensional module $X$ has a composition series $0=X_0\subset X_1\subset\cdots\subset X_\ell=X$, and the simple modules $X_j/X_{j-1}$, counted with multiplicity, are the \emph{composition factors} of $X$; in the Grothendieck ring, $[X]=\sum_j[X_j/X_{j-1}]$. The same terminology is used in any abelian category of finite length, for instance in the localized category of \S\ref{subsec-braider}.

\subsubsection{Universal $R$-matrices and the invariant $\Lambda$}\label{subsec:Runiv}

Universal $R$-matrices are used twice in this paper: to define the invariants $\Lambda(M,N)$ which give the quantum form of the KKOP seed, and in \S\ref{subsec-braider}, where suitably normalized $R$-matrices become the braiders of the localization.

Let $M,N\in\cC_{\mathfrak g}$. The universal $R$-matrix of $U_q'(\slh)$ converges $z$-adically on $M\otimes N_z$, and composed with the flip of tensor factors it gives an isomorphism of $\mathbf k((z))\otimes U_q'(\slh)$-modules
\begin{align}\label{eq-Runiv}
R^{\mathrm{univ}}_{M,N_z}\colon\mathbf k((z))\otimes_{\mathbf k[z^{\pm1}]}(M\otimes N_z)\xrightarrow{\ \sim\ }\mathbf k((z))\otimes_{\mathbf k[z^{\pm1}]}(N_z\otimes M),
\end{align}
the \emph{universal $R$-matrix of $M$ and $N$}~\cite{Dri87,FrRe92}, \cite{Kas02}, \cite{KKOP-comp}, \cite[\S1.3]{KKOP3}. It has the following properties.
\begin{enumerate}[label=(R\arabic*)]
\item\label{R:nat} \emph{Naturality.} For every morphism $f\colon N\to N'$, $(f_z\otimes\mathrm{id}_M)\circ R^{\mathrm{univ}}_{M,N_z}=R^{\mathrm{univ}}_{M,N'_z}\circ(\mathrm{id}_M\otimes f_z)$; similarly in $M$.
\item\label{R:hex} \emph{Hexagon.} $R^{\mathrm{univ}}_{M,(N\otimes N')_z}=(\mathrm{id}_{N_z}\otimes R^{\mathrm{univ}}_{M,N'_z})\circ(R^{\mathrm{univ}}_{M,N_z}\otimes\mathrm{id}_{N'_z})$.
\end{enumerate}
We regard $\mathbf k(z)\subseteq\mathbf k((z))$. Every $g\in\mathbf k(z)^\times$ can be written uniquely as $g=g_1/g_2$ with $g_1,g_2\in\mathbf k[z]$ coprime and $g_2$ monic; $g$ has a \emph{pole} (resp.\ \emph{zero}) at $z=a$ if $g_2(a)=0$ (resp.\ $g_1(a)=0$). We say that $R^{\mathrm{univ}}_{M,N_z}$ is \emph{rationally renormalizable} if there is $c_{M,N}(z)\in\mathbf k((z))^\times$ such that $R^{\mathrm{ren}}_{M,N_z}:=c_{M,N}(z)R^{\mathrm{univ}}_{M,N_z}$ maps $M\otimes N_z$ into $N_z\otimes M$, and its specialization at $z=x$ is non-zero for every $x\in\mathbf k^\times$. The \emph{renormalizing coefficient} $c_{M,N}(z)$ is then unique up to $\mathbf k[z^{\pm1}]^\times$, and
\begin{align}\label{eq-rmatrix}
\mathbf r_{M,N}:=R^{\mathrm{ren}}_{M,N_z}|_{z=1}\colon M\otimes N\to N\otimes M
\end{align}
is the \emph{$R$-matrix} of $M$ and $N$.
\begin{enumerate}[label=(R\arabic*),resume]
\item\label{R:ren} \emph{Rational renormalization}~\cite{Kas02,KKK18}, \cite[\S1.3]{KKOP3}. If $M$ and $N$ are simple, then $R^{\mathrm{univ}}_{M,N_z}$ is rationally renormalizable, and $\mathbf r_{M,N}\ne0$.
\end{enumerate}

Put $\widetilde p:=q^{2n}$ and $\varphi(z)=\prod_{s\ge0}(1-\widetilde p^{\,s}z)\in\mathbf k[[z]]$, and let $\cG\subseteq\mathbf k((z))^\times$ be the subgroup of elements $f(z)=c\,z^m\prod_{a\in\mathbf k^\times}\varphi(az)^{\eta_a}$ ($c\in\mathbf k^\times$, $m\in\Z$, almost all $\eta_a=0$). It contains $\mathbf k(z)^\times$, and every renormalizing coefficient lies in $\cG$~\cite[\S3]{KKOP-comp}. Define the homomorphism $\operatorname{Deg}\colon\cG\to\Z$ by~\cite{KKOP-comp}, \cite[\S1.3]{KKOP3}
\begin{align}\label{eq-Deg}
\operatorname{Deg}(f)=\sum_{a\in\widetilde p^{\,\Z_{\le0}}}\eta_a-\sum_{a\in\widetilde p^{\,\Z_{>0}}}\eta_a .
\end{align}
It vanishes on $\mathbf k^\times z^{\Z}$, so $\operatorname{Deg}(c_{M,N})$ does not depend on the choice of $c_{M,N}$.

\begin{definition}[{\cite{KKOP-comp}, \cite[Def.~1.1]{KKOP3}}]\label{def:Lambda}
For non-zero modules $M,N$ such that $R^{\mathrm{univ}}_{M,N_z}$ is rationally renormalizable, put
\begin{align}\label{eq-Lambda-def}
\Lambda(M,N)=\operatorname{Deg}\bigl(c_{M,N}(z)\bigr),\qquad \mathfrak d(M,N)=\tfrac12\bigl(\Lambda(M,N)+\Lambda(N,M)\bigr).
\end{align}
We only use $\mathfrak d(M,N)$ for simple $M,N$. KKOP define $\mathfrak d(M,N)=\frac12(\Lambda(M,N)+\Lambda(\mathscr D^{-1}M,N))$; for simple $M,N$ this agrees with the formula above by~\cite[Prop.~1.2(i)]{KKOP3}, and $\mathfrak d(M,N)\in\Z_{\ge0}$~\cite{KKOP-comp}, \cite[Prop.~1.2]{KKOP3}. A simple module $M$ is \emph{real} if $M\otimes M$ is simple, and two simple modules $M,N$ \emph{strongly commute} if $M\otimes N$ is simple.
\end{definition}

\begin{enumerate}[label=(R\arabic*),resume]
\item\label{R:comm} \emph{Commutation}~\cite{KKKO15,KKOP-comp}, \cite[Prop.~1.2]{KKOP3}. If $M$ is real simple and $N$ is simple, then $\mathfrak d(M,N)=0$ iff $\mathbf r_{M,N}$ is an isomorphism, iff $M$ and $N$ strongly commute.
\end{enumerate}
For a module $X$ of finite length, its \emph{head} $\operatorname{hd}(X)$ is its largest semisimple quotient. For simple $M,N$ with one of them real, $M\nabla N:=\operatorname{hd}(M\otimes N)$ is simple~\cite{KKKO15}, \cite{KKOP-comp}.

\subsubsection{Braid-word categories}\label{subsec_braid_cat}
For every $Q$-datum $\cQ$ there is a canonical complete duality datum $\mathbb D_{\cQ}=\{L^{\cQ}(\alpha_\imath)\}_{\imath\in\mathsf I}$, consisting of $n-1$ fundamental modules~ \cite{Qin24a}, \cite{Qin24b}, \cite{KKOP-PBW}, \cite[\S2, (2.5)]{KKOP3}. We fix such a $\mathbb D=\mathbb D_{\cQ}$ from now on.

Let $\mathrm{Br}_n$ be the braid group with generators $\sigma_1,\dots,\sigma_{n-1}$, let $\mathrm{Br}_n^+$ be the monoid of positive braids, and let $\mathbf i=(i_1,\dots,i_r)$ be a word for $\mathsf b\in\mathrm{Br}_n^+$, i.e.\ $\mathsf b=\sigma_{i_1}\cdots\sigma_{i_r}$ with $i_s\in\{1,\dots,n-1\}$. For integers $a\le b$, we write $[a,b]:=\{a,a+1,\dots,b\}$ for the set of integers $x$ with $a\le x\le b$. We call $s\in[1,r]$ a \emph{position} of $\mathbf i$, and $i_s$ its \emph{colour}. For a position $s$, let
\begin{align}\label{eq-spm}
s^-=\max\{s'<s\mid i_{s'}=i_s\},\qquad s^+=\min\{s'>s\mid i_{s'}=i_s\},
\end{align}
with $s^-=-\infty$ and $s^+=+\infty$ if the set is empty. A position $s$ is a \emph{first occurrence} if $s^-=-\infty$, and a \emph{last occurrence} if $s^+=+\infty$. KKOP define the \emph{affine cuspidal modules} $C_s=C^{\mathbb D,\mathbf i}_s$, $1\le s\le r$, via the braid group action on classes~\cite[(4.4)]{KKOP3}, \cite{KKOP-PBW,KKOP-braid}; these are real simple root modules. For positions $a\le b$ with $i_a=i_b$, let $a=s_1<s_2<\dots<s_l=b$ be the positions in $[a,b]$ of colour $i_a$. The \emph{affine determinantal module} is~\cite[Def.~4.8]{KKOP3}
\begin{align}\label{eq-Mab}
M[a,b]=\operatorname{hd}\bigl(C_{s_l}\otimes C_{s_{l-1}}\otimes\cdots\otimes C_{s_1}\bigr);
\end{align}
it is a real simple module~\cite[Thm.~5.16(iii)]{KKOP3}. For example, $M[a,a]=C_a$ and $M[a,a^+]=C_{a^+}\nabla C_a$. The \emph{braid-word category} $\cC(\mathsf b)=\cC^{\mathbb D}_{\mathfrak g}(\mathsf b)$ is the smallest full subcategory of $\cC^0_{\mathfrak g}$ containing $\mathbf 1$ and $C_1,\dots,C_r$, and stable under tensor products, subquotients and extensions. It depends only on $\mathsf b$~\cite[Cor.~5.34]{KKOP3}.

\subsubsection{The quiver of a word}

\begin{definition}[{\cite[\S6.1]{Qin24a},\cite[\S8.4]{KKOP3}, \cite[(2.14)]{Bi}}]\label{def:BFZ}
Let $\mathbf i=(i_1,\dots,i_r)$ be a word in which every colour $1,\dots,n-1$ occurs. Its \emph{quiver} $Q(\mathbf i)$ has vertex set $[1,r]$. The frozen vertices are the $n-1$ last occurrences, and all other positions are mutable. There are two types of arrows:
\begin{enumerate}
\item[$(Q1)$:] an arrow $s\to s^-$ for every $s$ which is not a first occurrence;
\item[$(Q2)$:] an arrow $s\to s'$ for every pair $s,s'$ with $|i_s-i_{s'}|=1$ and $s<s'<s^+<s'^+$. Here $+\infty$ is larger than every position, so there are type (Q2) arrows ending at last occurrences, but no type (Q2) arrow starts at a last occurrence.
\end{enumerate}

We denote by $Q_{\mathbf i}$ the signed adjacency matrix of $Q(\mathbf i)$; there are no arrows between frozen vertices.
\end{definition}

In the notation of~\cite[\S8.4]{KKOP3}, for $s'$ mutable, $Q_{\mathbf i}(s,s')=1$ if $s=s'^+$ or ($s<s'<s^+<s'^+$ and $|i_s-i_{s'}|=1$), and $Q_{\mathbf i}(s,s')=-1$ if $s'=s^+$ or ($s'<s<s'^+<s^+$ and $|i_s-i_{s'}|=1$).

\begin{example}\label{ex:BFZ}
Let $n=3$ and $\mathbf i=(1,2,1,1,2,1)$. Then $1^+=3$, $2^+=5$, $3^+=4$, $4^+=6$, and $5$, $6$ are the last occurrences. The type (Q1) arrows are $3\to1$, $4\to3$, $5\to2$ and $6\to4$, and the type (Q2) arrows are $1\to2$, $2\to4$ and $4\to5$. For instance, there is no arrow $1\to5$, since $5>1^+=3$. With rows indexed by $1,\dots,6$ and columns by the mutable vertices $1,2,3,4$,
\[
\bigl(Q_{\mathbf i}(s,s')\bigr)_{s\in[1,6],\ s'\in[1,4]}=\begin{pmatrix}
0&1&-1&0\\
-1&0&0&1\\
1&0&0&-1\\
0&-1&1&0\\
0&1&0&-1\\
0&0&0&1
\end{pmatrix}.
\]
The rows $3,5,4,6$, i.e.\ the rows $s'^+$ for $s'=1,2,3,4$, form the unimodular minor of Lemma~\ref{lem:unimod}.
\end{example}

Definition~\ref{def:BFZ} is the exchange matrix of~\cite[\S8.4, Thm.~9.7]{KKOP3} (Theorem~\ref{thm:KKOP}), as well as that of the right inductive weave of a braid variety in the form of~\cite[(2.14)]{Bi} (\S\ref{subsec:braidvar}). Furthermore, by \cite[Theorem 1.1]{CQW26}, this matrix coincides with the exchange matrix in \cite[\S6]{Qin24a}, following the work of Berenstein, Fomin, and Zelevinsky~\cite{BZ,BFZ05}. The type (Q2) arrows must be drawn even when they end at a last occurrence: in Lemma~\ref{lem:M} they correspond to the arrows of the skein quiver into the frozen vertices on the edge $(1,m)$.

\begin{lemma}\label{lem:Q2unique}
Let $x$ be a position and let $c'$ be a colour such that $|c'-i_x|=1$.
\begin{enumerate}[label={\rm (\alph*)}]
    \item If $x$ is a last occurrence, then there is no type (Q2) arrow starting at $x$.
    \item Otherwise, there exists a type (Q2) arrow from $x$ to a position of colour $c'$ if and only if the last position $y<x^+$ of colour $c'$ exists and satisfies $y>x$. In this case, $x\to y$ is the unique such arrow.
\end{enumerate}
Moreover, the quiver $Q(\mathbf{i})$ contains no multi-arrows and no cycles of length $2$.
\end{lemma}

\begin{proof}
Part (a) is part of Definition~\ref{def:BFZ}.

(b) Let $y$ be a position of colour $c'$. By Definition~\ref{def:BFZ}, there is a type (Q2) arrow $x\to y$ if and only if $x<y<x^+<y^+$. Since $y$ and $y^+$ are consecutive positions of colour $c'$, the condition $y<x^+<y^+$ holds if and only if $y$ is the last position of colour $c'$ before $x^+$. Hence there is at most one such $y$, and the arrow $x\to y$ exists if and only if, in addition, $y>x$.

For the last statement, note that a type (Q1) arrow joins two positions of the same colour, while a type (Q2) arrow joins two positions of adjacent colours; hence no two vertices are joined by arrows of different types. Two positions $x<y$ of the same colour are joined by a type (Q1) arrow only if $x=y^-$, and then only by the single arrow $y\to x$. Since every type (Q2) arrow $x\to y$ satisfies $x<y$ and is uniquely determined by the ordered pair $(x,y)$, there exists at most one type (Q2) arrow connecting any two adjacent positions of distinct colours. Therefore $Q(\mathbf i)$ has no multi-arrows and no cycles of length $2$.
\end{proof}

The following lemma will be used in Lemma~\ref{lem:Lst} to describe the cokernel of the exchange matrix of the skein seed, and in Proposition~\ref{prop:homog-simple}.

\begin{lemma}\label{lem:unimod}
Let $N$ be the number of mutable positions of $\mathbf i$.Consider the $N\times N$ submatrix of $Q_{\mathbf i}$ whose rows and columns are indexed by $\{s^+\mid s\text{ mutable}\}$ and $\{s\mid s\text{ mutable}\}$, respectively. When the rows and columns are ordered by $s^+$, this submatrix is triangular with diagonal entries equal to $1$. Since the rows $\{s^+\mid s\text{ mutable}\}$ correspond exactly to the positions that are not first occurrences, the column vectors ${\bf q}_s$ for $s$ mutable span a direct summand of $\Z^{[1,r]}$, with complement $\Z^{\{\text{first occurrences}\}}$.
\end{lemma}

\begin{proof}
For mutable $s$, $s^+$ is finite and is not a first occurrence, and conversely a position $x$ which is not a first occurrence is $s^+$ for the mutable position $s=x^-$. The entry in row $s^+$ and column $s$ is $1$, from the type (Q1) arrow $s^+\to s$. Let $s\ne s'$ be mutable with $Q_{\mathbf i}(s^+,s')\ne0$, and put $x=s^+$. If $x$ and $s'$ are joined by a type (Q1) arrow, then either $x\to x^-=s'$, which forces $s'=s$, or $s'\to s'^-=x$, and then $s'^+>s'>x$. If they are joined by a type (Q2) arrow, then $x<s'<x^+<s'^+$ or $s'<x<s'^+<x^+$. In all cases $s'^+>s^+$, which proves triangularity. Hence the restriction of the vectors ${\bf q}_s$ to the rows $\{s^+\}$ is a unimodular matrix, and $\Z^{[1,r]}=\operatorname{span}\{{\bf q}_s\}\oplus\Z^{\{\text{first occurrences}\}}$.
\end{proof}

\subsubsection{The word $u$ and its powers}\label{subsec-wordu}

Let $s_i=(i\ \,i{+}1)\in S_n$. The \emph{length} $\ell(w)$ of $w\in S_n$ is the smallest $\ell$ with $w=s_{i_1}\cdots s_{i_\ell}$, and such a word is a \emph{reduced word} of $w$. The longest element $w_0\in S_n$, $i\mapsto n+1-i$, has length $\binom n2$. The map $s_i\mapsto\sigma_i$ sends reduced words of $w_0$ to words of a single positive braid, the \emph{Garside element} $\Delta\in\mathrm{Br}_n$~\cite{Gar69}.

\begin{definition}\label{def:u}
Let
\begin{align}\label{eq-u}
u=(1\,2\,\cdots\,n{-}1)(1\,2\,\cdots\,n{-}2)\cdots(1\,2)(1),
\end{align}
which is a reduced word of $w_0$. The $r$-th \emph{block} of $u$ is $(1,2,\dots,n-r)$, $1\le r\le n-1$. For $l\ge1$, $u^l$ is the concatenation of $l$ copies of $u$; it is a word for $\Delta^l$.
\end{definition}

We label the positions of $u^l$ by triples $\langle\kappa,r,c\rangle$ with $1\le\kappa\le l$, $1\le r\le n-1$ and $1\le c\le n-r$, where $\langle\kappa,r,c\rangle$ is the occurrence of the colour $c$ in the $r$-th block of the $\kappa$-th copy. Positions are ordered lexicographically in $(\kappa,r)$, and by increasing $c$ within a block.

\begin{lemma}\label{lem:Quk}
For a colour $c$, the positions of colour $c$ in $u^l$ are given by $\langle\kappa,r,c\rangle$ with $1\le\kappa\le l$ and $1\le r\le n-c$. The last occurrence of $c$ is $\langle l,n-c,c\rangle$, and for every other position of colour $c$, we have
\begin{align}\label{eq-Quk-next}
\langle\kappa,r,c\rangle^+=\begin{cases}\langle\kappa,r+1,c\rangle & \text{if } r<n-c,\\ \langle\kappa+1,1,c\rangle & \text{if } r=n-c\text{ and }\kappa\le l-1.\end{cases}
\end{align}
Moreover, the quiver $Q(u^l)$ has exactly the following three types of arrows:
\begin{enumerate}[label={\rm (\alph*)}]
\item the type (Q1) arrows $x^+\to x$ for every position $x$ that is not a last occurrence, with $x^+$ given by~\eqref{eq-Quk-next};
\item the type (Q2) arrows $\langle\kappa,r,c\rangle\to\langle\kappa,r,c+1\rangle$ for $1\le\kappa\le l$, $1\le c\le n-2$, and $1\le r\le n-c-1$;
\item the type (Q2) arrows $\langle\kappa,r-1,c+1\rangle\to\langle\kappa,r,c\rangle$ for $1\le\kappa\le l$, $1\le c\le n-2$, and $2\le r\le n-c-1$, as well as $\langle\kappa-1,n-c-1,c+1\rangle\to\langle\kappa,1,c\rangle$ for $2\le\kappa\le l$ and $1\le c\le n-2$.
\end{enumerate}
\end{lemma}

\begin{proof}
Since the colour $c$ occurs exactly once in each of the blocks $r=1,\dots,n-c$ of every copy of $u$, this yields the description of the positions of colour $c$, the last occurrence, and~\eqref{eq-Quk-next}, thereby establishing (a). It remains to determine the type (Q2) arrows. Fix adjacent colours $a=c$ and $b=c+1$, with $1\le c\le n-2$. The condition $x<y<x^+<y^+$ involves only the subword of $u^l$ formed by the letters $a$ and $b$. Each copy of $u$ contributes the string $(ab)(ab)\cdots(ab)\,a$, consisting of $n-c-1$ pairs $ab$ (one in each block $r\le n-c-1$), followed by the $a$ of block $n-c$. By Lemma~\ref{lem:Q2unique}, a type (Q2) arrow from $x$ points to the last letter $y$ of the other colour preceding $x^+$, provided that $x$ is not a last occurrence and $y>x$.

Let $x=\langle\kappa,r,c\rangle$ be of colour $a$. If $r\le n-c-1$, the letter following $x$ in the subword is the $b$ of its block, which is immediately followed by $x^+$; this yields the arrows in (b). If $r=n-c$, then $x^+$ immediately follows $x$ in the subword, or $x$ is a last occurrence; in either case, there is no type (Q2) arrow from $x$ to a position of colour $b$.

Let $x=\langle\kappa,r,c+1\rangle$ be of colour $b$. If $r\le n-c-2$, then $x^+$ is the $b$ of block $r+1$, which is immediately preceded by the $a$ of block $r+1$; this yields the first family of arrows in (c). If $r=n-c-1$ and $\kappa\le l-1$, then $x^+$ is the first $b$ of the next copy, which is immediately preceded by its first $a$; this yields the second family of arrows in (c). If $r=n-c-1$ and $\kappa=l$, then $x$ is a last occurrence, and thus there is no type (Q2) arrow from $x$.
\end{proof}

\subsubsection{The KKOP seed}

Fix a word $\mathbf i=(i_1,\dots,i_r)$ for $\mathsf b$ in which every colour occurs, and let $s^{\min}$ be the first position of colour $i_s$. Put
\begin{align}\label{eq-Ms}
\mathsf M_s=M[s^{\min},s]\qquad(1\le s\le r).
\end{align}
The frozen modules of $\cC(\mathsf b)$ are the $\mathsf M_s$ with $s$ a last occurrence, i.e.\ the modules $F=M[a,b]$ with $a$ and $b$ the first and last occurrences of a colour. Let
\begin{align}\label{eq-Lambda-i}
\Lambda_{\mathbf i}=\bigl(\Lambda(\mathsf M_s,\mathsf M_{s'})\bigr)_{s,s'\in[1,r]}.
\end{align}
By~\cite[\S8.4]{KKOP3}, for $s\le s'$,
\begin{align}\label{eq-Lambda-formula}
\Lambda_{\mathbf i}(s,s')=-\bigl(\varpi_{i_s}-w_{\le s}\varpi_{i_s},\ \varpi_{i_{s'}}+w_{\le s'}\varpi_{i_{s'}}\bigr),\qquad w_{\le s}=s_{i_1}s_{i_2}\cdots s_{i_s}\in\mathsf W,
\end{align}
and $\Lambda_{\mathbf i}(s',s)=-\Lambda_{\mathbf i}(s,s')$.

The $(q,t)$-characters $[\mathsf M_s]_t$ of \S\ref{subsec:qGr} below lie in the quantum torus $\cY_t$ over $\Z[t^{\pm1/2}]$ of \S\ref{subsec:qGr}. Recall that $\Z[t^{\pm1/2}]\subseteq R$ via $t^{1/2}=\xi^{n/2}$ (see \eqref{eq-constants}), and let $\cD_t^R$ be the skew-field of fractions of the quantum torus $\cY_t\otimes_{\Z[t^{\pm1/2}]}R$, which is an Ore domain. Let $\mathcal F_t\subseteq\cD_t^R$ be the skew-field of fractions of the $R$-subalgebra of $\cD_t^R$ generated by the $[\mathsf M_s]_t^{\pm1}$, $1\le s\le r$ (the frozen classes are not invertible in $\cY_t$, so we work in $\cD_t^R$); by Theorem~\ref{thm:KKOP} below, this subalgebra is the based quantum torus $\bT(\mathsf s(\mathbf i))$.

\begin{definition}[{\cite[\S8.4]{KKOP3}}]\label{def:KKOPseed}
The \emph{KKOP seed} of $\mathbf i$ is the triple
\begin{align}\label{eKKOPseed}
\mathsf s(\mathbf i)=(Q_{\mathbf i},\ \Pi_{\mathbf i},\ M_{\mathbf i}),\qquad \Pi_{\mathbf i}=n\Lambda_{\mathbf i},\qquad M_{\mathbf i}({\bf k})=\Bigl[\prod_{s=1}^r[\mathsf M_s]_t^{k_s}\Bigr],
\end{align}
where $Q_{\mathbf i}$ is the matrix of Definition~\ref{def:BFZ}, and the Weyl-ordered product \eqref{Weyl-A} is taken with respect to $\Pi_{\mathbf i}$ and the variable $\xi$.
\end{definition}

Since $t^{1/2}=\xi^{n/2}$, we have $\xi^{\frac12{\bf k}\Pi_{\mathbf i}{\bf t}^T}=t^{\frac12{\bf k}\Lambda_{\mathbf i}{\bf t}^T}$. Hence $M_{\mathbf i}$ is the toric frame of~\cite[\S8.4]{KKOP3}, written in the variable $\xi$; the factor $n$ in $\Pi_{\mathbf i}$ is what allows us to use the same variable $\xi$ as for the skein seed (Proposition~\ref{prop:CHWseed}). KKOP use the same quiver convention as ours.

\subsection{Monoidal categorification of $\mathscr A^+(\mathsf s(\mathbf i))$}\label{subsec-moncat}

\subsubsection{Quantum Grothendieck rings}\label{subsec:qGr}

The Grothendieck ring $K(\cC^0_{\mathfrak g})$ is commutative. Frenkel and Reshetikhin~\cite{FR} constructed an injective ring homomorphism, the \emph{$q$-character map},
\begin{align}\label{eq-qchar}
\chi_q\colon K(\cC^0_{\mathfrak g})\hookrightarrow\cY:=\Z\bigl[Y_{\imath,p}^{\pm1}:(\imath,p)\in\widehat\triangle_0\bigr],
\end{align}
which records the dimensions of the $\ell$-weight spaces of a module~\cite[\S6.1]{KKOP3}. Dominant monomials $m$ parametrize the simple modules $L(m)$ of $\cC^0_{\mathfrak g}$, and $L(Y_{\imath,p})=L(\imath,p)$.

Let $\cY_t$ be the quantum torus of $(q,t)$-characters of~\cite{Nakajima,VV,Hernandez}. It is generated over $\Z[t^{\pm1/2}]$ by $\widetilde Y_{\imath,p}^{\pm1}$, $(\imath,p)\in\widehat\triangle_0$, subject to $\widetilde Y_{\imath,p}\,\widetilde Y_{\jmath,s}=t^{\mathcal N(\imath,p;\jmath,s)}\,\widetilde Y_{\jmath,s}\,\widetilde Y_{\imath,p}$ for an explicit skew-symmetric form $\mathcal N$~\cite{Hernandez,HL15,FO21}, \cite[\S6.1]{KKOP3}. The specialization $\mathrm{ev}_{t=1}\colon\cY_t\to\cY$, $t^{1/2}\mapsto1$, $\widetilde Y_{\imath,p}\mapsto Y_{\imath,p}$, is a ring homomorphism. The \emph{bar involution} of $\cY_t$ is the ring anti-involution which fixes every $\widetilde Y_{\imath,p}$ and sends $t^{1/2}$ to $t^{-1/2}$. For each dominant monomial $m$ there is a unique bar-invariant element $L_t(m)\in\cY_t$ which is unitriangular, with respect to the Nakajima order, in the basis of $(q,t)$-characters of standard modules; it is the \emph{$(q,t)$-character} of $L(m)$~\cite{Nakajima,Hernandez}, \cite[\S6.1]{KKOP3}, and we write $[M]_t:=L_t(m)$ for a simple module $M\cong L(m)$. In type $A^{(1)}_{n-1}$ every simple module is quantizable, i.e.\ $\mathrm{ev}_{t=1}([M]_t)=\chi_q(M)$ (\cite{Nakajima}; see~\cite[Rem.~6.5]{KKOP3}).

\begin{definition}[{\cite[Def.~9.11]{KKOP3}}]\label{def:Kt}
For a monoidal subcategory $\cC\subseteq\cC^0_{\mathfrak g}$, the \emph{quantum Grothendieck ring} $\cK_t(\cC)$ is defined as the $\Z[t^{\pm1/2}]$-subalgebra of $\cY_t$ generated by the classes $[L]_t$ for simple objects $L \in \cC$.
\end{definition}

The classes $[L]_t$, where $L$ ranges over the simple objects in $\cC$, constitute a $\Z[t^{\pm1/2}]$-basis of $\cK_t(\cC)$. The evaluation map $\mathrm{ev}_{t=1}$ induces a surjection from $\cK_t(\cC)$ onto $\chi_q(K(\cC))\cong K(\cC)$, sending each $[L]_t$ to the class of $L$. The structure constants for the products of simple classes lie in $\Z_{\ge0}[t^{\pm1/2}]$~\cite{VV,Nakajima}, \cite[Prop.~6.1]{KKOP3}:
\begin{align}\label{eq-positivity}
[S]_t[S']_t=\sum_{S''}c^{S''}_{S,S'}(t)[S'']_t,\qquad c^{S''}_{S,S'}(t)\in\Z_{\ge0}[t^{\pm1/2}],
\end{align}
where $c^{S''}_{S,S'}(1)$ coincides with the multiplicity of $S''$ as a composition factor of $S\otimes S'$.

\subsubsection{The monoidal categorification theorem}

The following theorem of Kashiwara--Kim--Oh--Park is the representation-theoretic counterpart of Theorem~\ref{thm:CHW}. It extends earlier quantum cluster structures on quantum Grothendieck rings of monoidal subcategories of $\cC^0_{\mathfrak g}$~\cite{Qin17,Bit21,Qin24a,Qin24b,KKOP2}, which quantize the cluster structures of Hernandez--Leclerc~\cite{HL10,HL16}.

\begin{theorem}[{\cite[Theorem 8.19]{Qin24a}, \cite[Theorem 1.6]{Qin24b},\cite[Prop.~8.14, Thm.~8.15, Thm.~9.7, Thm.~9.12]{KKOP3}}]\label{thm:KKOP}
Let $\mathbf i$ be a word for $\mathsf b\in\mathrm{Br}_n^+$ in which every color occurs.
\begin{enumerate}[label={\rm (\alph*)}]
\item $\sum_{s''}Q_{\mathbf i}(s'',s)\Lambda_{\mathbf i}(s'',s')=2\delta_{s,s'}$ for any mutable $s$ and $s'\in[1,r]$. Consequently, $\mathsf s(\mathbf i)$ is a quantum seed in $\mathcal F_t$, compatible with the scalar $2$.
\item The inclusion $\cK_t(\cC(\mathsf b))\otimes_{\Z[t^{\pm1/2}]}R\subseteq\cY_t\otimes_{\Z[t^{\pm1/2}]}R$ induces an identification
\begin{align}\label{eKKOP-iso}
\cK_t(\cC(\mathsf b))\otimes_{\Z[t^{\pm1/2}]}R=\mathscr A^+(\mathsf s(\mathbf i)) .
\end{align}
\item Every cluster monomial of $\mathscr A^+(\mathsf s(\mathbf i))$ (Definition~\ref{def-quan-cluster-algebra}(d)) is the $(q,t)$-character of a real simple module in $\cC(\mathsf b)$.
\end{enumerate}
\end{theorem}

\begin{proof}
By~\cite[Prop.~8.14, Thm.~8.15]{KKOP3}, the family $\{\mathsf M_s\}_{1\le s\le r}$, together with the matrix $Q_{\mathbf i}$ defined in~\cite[(8.5)]{KKOP3}, is $\Lambda$-admissible. Specifically, it coincides with the datum $B(\mathfrak C^+)$ from~\cite[Def.~9.5]{KKOP3} for the admissible chain of $i$-boxes $\mathfrak C^+=(\{1,s])_{s}$. In the notation of~\cite{KKOP3}, the $i$-box $\{1,s]$ denotes the set of positions of color $i_s$ up to index $s$, so that $M_{\{1,s\}}=\mathsf M_s$; this establishes (a).

Statements (b) and (c) follow from~\cite[Thm.~9.7]{KKOP3}, combined with the isomorphism in~\cite[Thm.~6.2]{KKOP3} and the definitions in~\cite[Def.~9.11, Thm.~9.12]{KKOP3}. These results hold for the quantum seed over $\Z[t^{\pm1/2}]$ with exchange matrix $Q_{\mathbf i}$, form $\Lambda_{\mathbf i}$, normalization $X^{\bf k}X^{\bf t}=t^{\frac12{\bf k}\Lambda_{\mathbf i}{\bf t}^T}X^{{\bf k}+{\bf t}}$, and initial variables $[\mathsf M_s]_t$, where the quantum cluster algebra $\cK_t(\cC(\mathsf b))$ is generated by cluster monomials. Since $t^{\frac12{\bf k}\Lambda_{\mathbf i}{\bf t}^T}=\xi^{\frac12{\bf k}\Pi_{\mathbf i}{\bf t}^T}$, this seed and $\mathsf s(\mathbf i)$ share the same toric frame, exchange matrix, and mutation rules, and hence the same cluster monomials. The ring $R$ is a free $\Z[t^{\pm1/2}]$-module with basis $\{\omega^{a/2} \mid 0\le a<n^2\}$; thus, tensoring with $R$ is exact. Consequently, $\mathscr A^+(\mathsf s(\mathbf i))$, being the $R$-span of the products of cluster monomials, is the image of $\cK_t(\cC(\mathsf b))\otimes R$. This proves (b) and (c) over $R$. Finally, Formula~\eqref{eq-Lambda-formula} holds for arbitrary words, with $w_{\le s}=s_{i_1}\cdots s_{i_s}\in\mathsf W$, by~\cite[Lemma~5.21, Cor.~5.23]{KKOP3}.
\end{proof}

By Lemma~\ref{lem:normal}, the frozen variables $[F]_t$ of $\mathsf s(\mathbf i)$ are normal in $\mathscr A^+(\mathsf s(\mathbf i))$, and $\mathscr A(\mathsf s(\mathbf i))$ is the Ore localization of $\cK_t(\cC(\mathsf b))\otimes R$ at the frozen classes. In \S\ref{subsec-braider} we realize this localization as the quantum Grothendieck ring of a monoidal localization of the subcategory $\cC^{\mathrm s}(\mathsf b)\subseteq\cC(\mathsf b)$ (Theorem~\ref{thm:Ktloc}).

\subsubsection{Block weights and homogeneity of simple classes}\label{subsec:homog}

The twist of Theorem~\ref{thm:main} is defined using a grading of the quantum torus of $\mathsf s(\mathbf i)$ by the cokernel of the exchange matrix. In this subsection we show that the $(q,t)$-character of every simple object of $\cC(\mathsf b)$ is homogeneous for this grading.

By~\cite{CM05,KKOP-root,KKOP-PBW}, \cite[Thm.~2.4]{KKOP3}, the category $\cC^0_{\mathfrak g}$ decomposes into blocks labelled by the root lattice $\mathsf Q$, $\cC^0_{\mathfrak g}=\bigoplus_{\beta\in\mathsf Q}\cC^0_{\mathfrak g,\beta}$, with $\cC^0_{\mathfrak g,\beta}\otimes\cC^0_{\mathfrak g,\gamma}\subseteq\cC^0_{\mathfrak g,\beta+\gamma}$; the labelling depends on $\cQ$. For a simple module $L$, let $\operatorname{wt}_{\cQ}(L)\in\mathsf Q$ be the label of its block. Every composition factor of $L\otimes L'$ has weight $\operatorname{wt}_{\cQ}L+\operatorname{wt}_{\cQ}L'$.

\begin{lemma}\label{lem:wtgrading}
Let $\cC=\cC(\mathsf b)$. There is a $\mathsf Q$-grading of the algebra $\cK_t(\cC)$ in which each $[S]_t$ is homogeneous of degree $\operatorname{wt}_{\cQ}(S)$.
\end{lemma}

\begin{proof}
Grade the free $\Z[t^{\pm1/2}]$-module $\cK_t(\cC)$ by declaring $[S]_t$ to have degree $\operatorname{wt}_{\cQ}(S)$. In \eqref{eq-positivity}, $c^{S''}_{S,S'}\neq0$ implies $c^{S''}_{S,S'}(1)>0$, by positivity, so $S''$ is a composition factor of $S\otimes S'$ and has weight $\operatorname{wt}_{\cQ}S+\operatorname{wt}_{\cQ}S'$.
\end{proof}

Let $\mathbf i$ be a word of length $r$ representing $\mathsf b$, in which every color occurs. By Lemma~\ref{lem:unimod}, the span $L_{\mathbf i}$ of the column vectors ${\bf q}_s$ of $Q_{\mathbf i}$ (for $s$ mutable) is a direct summand of $\Z^r$, and
\begin{align}\label{eq-Gamma-i}
\Gamma_{\mathbf i}:=\Z^r/L_{\mathbf i}\cong\Z^{\{\text{first occurrences}\}}
\end{align}
is free of rank $n-1$. Let $p_{\mathbf i}\colon\Z^r\to\Gamma_{\mathbf i}$ denote the canonical projection. An element $x=\sum_{\bf k}c_{\bf k}M_{\mathbf i}({\bf k})$ in the based quantum torus $\bT(\mathsf s(\mathbf i))$ is called \emph{homogeneous of degree $\gamma\in\Gamma_{\mathbf i}$} if $p_{\mathbf i}({\bf k})=\gamma$ for all ${\bf k}$ with $c_{\bf k}\ne0$. According to Theorem~\ref{thm-inclusion-quantum}, every element of $\cK_t(\cC(\mathsf b))$ resides in $\bT(\mathsf s(\mathbf i))$.

\begin{proposition}\label{prop:homog-simple}
With the above notation:
\begin{enumerate}[label={\rm (\alph*)}]
\item There exists a unique injective homomorphism $P_{\rm wt}\colon\Gamma_{\mathbf i}\to\mathsf Q$ such that $\operatorname{wt}_{\cQ}\mathsf M_s=P_{\rm wt}(p_{\mathbf i}({\bf e}_s))$ for all $1\le s\le r$;
\item For every simple object $S$ of $\cC(\mathsf b)$, the class $[S]_t$ is homogeneous of degree $P_{\rm wt}^{-1}(\operatorname{wt}_{\cQ}S)$.
\end{enumerate}
\end{proposition}

\begin{proof}
(a) Define $p_{\rm wt}\colon\Z^r\to\mathsf Q$ by $p_{\rm wt}({\bf e}_s)=\operatorname{wt}_{\cQ}\mathsf M_s$. For ${\bf c}\in\Z^r_{\ge0}$, the element $M_{\mathbf i}({\bf c})$ is a power of $\xi^{1/2}$ multiplied by $[\mathsf M_1]_t^{c_1}\cdots[\mathsf M_r]_t^{c_r}$; hence, it is homogeneous of degree $p_{\rm wt}({\bf c})$ with respect to the grading in Lemma~\ref{lem:wtgrading} (extended to $\cK_t(\cC)\otimes R$). Fix a mutable index $k$ and set ${\bf b}^{\pm}=[\pm{\bf q}_k]_+$. Multiplying the exchange relation~\eqref{eq-exchange} on the left by $M_{\mathbf i}({\bf e}_k)$ yields
\[
M_{\mathbf i}({\bf e}_k)A'_k=\xi^{a_+}M_{\mathbf i}({\bf b}^+)+\xi^{a_-}M_{\mathbf i}({\bf b}^-)
\]
for some $a_\pm\in\tfrac12\Z$. By Theorem~\ref{thm:KKOP}(c), $A'_k$ is the $(q,t)$-character of a real simple module, implying that the left-hand side is homogeneous. Since both terms on the right are non-zero and homogeneous of degrees $p_{\rm wt}({\bf b}^\pm)$, we must have $p_{\rm wt}({\bf q}_k)=p_{\rm wt}({\bf b}^+)-p_{\rm wt}({\bf b}^-)=0$. Consequently, $p_{\rm wt}$ factors through $p_{\mathbf i}$ as $P_{\rm wt}\circ p_{\mathbf i}$, and $P_{\rm wt}$ is unique because $p_{\mathbf i}$ is surjective.

Since both $\Gamma_{\mathbf i}$ and $\mathsf Q$ are free of rank $n-1$, the map $P_{\rm wt}$ is injective if and only if the weights $\operatorname{wt}_{\cQ}\mathsf M_s$ span $\mathsf Q\otimes\Q$. Let $s_1=s^{\min}<\dots<s_l=s$ denote the positions $\le s$ of color $i_s$. By \eqref{eq-Mab}, $\mathsf M_s$ is a quotient of $C_{s_l}\otimes\cdots\otimes C_{s_1}$, so $\operatorname{wt}_{\cQ}\mathsf M_s=\sum_{j}\operatorname{wt}_{\cQ}C_{s_j}$. This relation represents a unitriangular change of basis; thus, it suffices to show that the weights $\operatorname{wt}_{\cQ}C_s$ span $\mathsf Q\otimes\Q$. Let $\beta_s=s_{i_1}\cdots s_{i_{s-1}}(\alpha_{i_s})$. The cuspidal module $C_s$ is obtained by applying the braid group action $T_{i_1}\cdots T_{i_{s-1}}$ to the class of the module $L^{\cQ}(\alpha_{i_s})$ from the duality datum~\cite[(4.4)]{KKOP3}. The block weight of $L^{\cQ}(\alpha_\imath)$ is $\alpha_\imath$~\cite[Thm.~2.4(i)]{KKOP3}. The braid group action is defined on the bosonic extension, where $\operatorname{wt}(T_\imath x)=s_\imath\operatorname{wt}(x)$~\cite[\S3.4]{KKOP3}. Under the isomorphism of~\cite[Thm.~6.2]{KKOP3}, the generators $f_{\imath,0}$ (of weight $-\alpha_\imath$) map to the classes of $L^{\cQ}(\alpha_\imath)$. Therefore, for homogeneous elements, block weights correspond to the negatives of the weights in the bosonic extension, and $[C_s]$ is the image of $T_{i_1}\cdots T_{i_{s-1}}f_{i_s,0}$ by~\cite[(4.4)]{KKOP3}. It follows that $\operatorname{wt}_{\cQ}C_s=\beta_s$, reducing the problem to showing that the $\beta_s$ span $\mathsf Q\otimes\Q$. Since $s_\imath(x)\in x+\Q\alpha_\imath$, induction yields $\operatorname{span}(\beta_1,\dots,\beta_s)=\operatorname{span}(\alpha_{i_1},\dots,\alpha_{i_s})$. For $s=r$, this span is $\mathsf Q\otimes\Q$ because every color occurs in $\mathbf i$.

(b) Expand $[S]_t=\sum_{\bf a}c_{\bf a}M_{\mathbf i}({\bf a})$, and choose ${\bf d}\in\Z^r_{\ge0}$ such that ${\bf d}+{\bf a}\in\Z^r_{\ge0}$ whenever $c_{\bf a}\ne0$. Then
\[
M_{\mathbf i}({\bf d})[S]_t=\sum_{\bf a}c_{\bf a}\xi^{\frac12{\bf d}\Pi_{\mathbf i}{\bf a}^T}M_{\mathbf i}({\bf d}+{\bf a}).
\]
The left-hand side is homogeneous of degree $p_{\rm wt}({\bf d})+\operatorname{wt}_{\cQ}S$ by Lemma~\ref{lem:wtgrading}, while each $M_{\mathbf i}({\bf d}+{\bf a})$ is homogeneous of degree $p_{\rm wt}({\bf d}+{\bf a})$ by part (a). Since the elements $M_{\mathbf i}({\bf d}+{\bf a})$ are linearly independent, $c_{\bf a}\ne0$ implies $P_{\rm wt}(p_{\mathbf i}({\bf a}))=\operatorname{wt}_{\cQ}S$. The conclusion follows from the injectivity of $P_{\rm wt}$.
\end{proof}

\subsection{Braider and localization}\label{subsec-braider}

In this subsection we construct a monoidal localization $\widetilde\cC$, at the frozen modules, of the full subcategory $\cC^{\mathrm s}(\mathsf b)\subseteq\cC(\mathsf b)$ of subquotients of tensor products of cuspidal modules (Definition~\ref{def:Cs}), which has the same simple objects as $\cC(\mathsf b)$. We show that its quantum Grothendieck ring is the Ore localization of $\cK_t(\cC(\mathsf b))$ at the frozen classes (Theorem~\ref{thm:Ktloc}). Braiders in the sense of~\cite{KKOP-loc} are what is needed to invert objects of a non-braided monoidal category. The localizations of~\cite{KKOP-loc} concern quiver Hecke algebras, where the braiders come from the $R$-matrices of affinizations of determinantial modules. For braids given by reduced words, they could be transported to the quantum affine side through the monoidal equivalences of Fujita~\cite{Fuj} (see also~\cite{KKK15II}); to our knowledge, this has not been written down, and it does not cover non-reduced braids such as $\Delta^{m-2}$, $m\ge4$. We instead use renormalized universal $R$-matrices, and we make them coherent by restricting to $\cC^{\mathrm s}(\mathsf b)$. We do not know whether $\cC^{\mathrm s}(\mathsf b)=\cC(\mathsf b)$; the restriction is needed because, for non-split extensions, the renormalized $R$-matrices need not specialize well at $z=1$.

Throughout this subsection, $\mathbf i=(i_1,\dots,i_r)$ is a word for $\mathsf b$ in which every colour occurs, $\cC=\cC(\mathsf b)$, $C_1,\dots,C_r$ are its affine cuspidal modules, and $F_1,\dots,F_{n-1}$ are its frozen modules, where $F_i$ is the frozen module of colour $i$.

\subsubsection{Renormalized $R$-matrices on subquotients}

Let
\begin{align}\label{eq-H}
\cH=\bigl\{g=g_1/g_2\in\mathbf k(z)^\times\mid g_1(1)\ne0\text{ and }g_2(1)\ne0\bigr\}
\end{align}
be the group of non-zero rational functions having neither a zero nor a pole at $z=1$. The evaluation map $g\mapsto g(1)$ defines a homomorphism $\cH\to\mathbf k^\times$. A $\mathbf k((z))$-linear map between the modules $\mathbf k((z))\otimes_{\mathbf k[z^{\pm1}]}N_z$ and $\mathbf k((z))\otimes_{\mathbf k[z^{\pm1}]}N'_z$ is said to be \emph{defined over $\mathbf k[z^{\pm1}]$} if it maps $N_z$ into $N'_z$; its \emph{specialization} at $z=1$ is the induced map $N\to N'$. This specialization process commutes with composition and direct sums.

More generally, let $\mathcal O=\mathbf k[z]_{(z-1)}\subseteq\mathbf k(z)\subseteq\mathbf k((z))$ denote the local ring at $z=1$, so that $\cH=\mathcal O^\times$. A $\mathbf k((z))$-linear map as above is \emph{defined over $\mathcal O$} if it maps $\mathcal O\otimes_{\mathbf k[z^{\pm1}]}N_z$ into $\mathcal O\otimes_{\mathbf k[z^{\pm1}]}N'_z$. Its specialization at $z=1$ is then the induced map modulo the maximal ideal $(z-1)\mathcal O$. Any map defined over $\mathbf k[z^{\pm1}]$ is automatically defined over $\mathcal O$ with the same specialization. Furthermore, multiplying a map defined over $\mathcal O$ by $g\in\cH$ yields a map defined over $\mathcal O$ whose specialization is scaled by $g(1)$. Specialization from $\mathcal O$ also commutes with composition and direct sums. Consequently, the argument of Lemma~\ref{lem:Rsub}(b) applies verbatim over $\mathcal O$, as $\mathcal O\otimes N_z$ is a free $\mathcal O$-module.

\begin{lemma}\label{lem:Rsub}
Let $M$ be a real simple module, and let $T=X_1\otimes\cdots\otimes X_l$ with each $X_j$ simple and $\mathfrak d(M,X_j)=0$. Put $c_T=\prod_jc_{M,X_j}$. Then:
\begin{enumerate}[label={\rm (\alph*)}]
\item $c_T(z)R^{\mathrm{univ}}_{M,T_z}$ is defined over $\mathbf k[z^{\pm1}]$, and its specialization at $z=1$ is the isomorphism obtained by composing the $\mathbf r_{M,X_j}$ as in~\ref{R:hex};
\item For every subquotient $Y$ of $T$, $c_T(z)R^{\mathrm{univ}}_{M,Y_z}$ is defined over $\mathbf k[z^{\pm1}]$, and its specialization at $z=1$ is an isomorphism $M\otimes Y\xrightarrow{\sim}Y\otimes M$; these specializations are compatible with inclusions and projections;
\item If $S$ is a simple subquotient of $T$, then $c_T/c_{M,S}\in\cH$ and $\mathfrak d(M,S)=0$.
\end{enumerate}
\end{lemma}

\begin{proof}
(a) This follows directly from~\ref{R:hex} and~\ref{R:comm}.

(b) Let $Y'\subseteq T$ be a submodule. By~\ref{R:nat}, the map constructed in (a) restricts to a map on $c_T(z)R^{\mathrm{univ}}_{M,Y'_z}$. Choose a $\mathbf k$-linear complement $W$ of $Y'$ in $T$. Then $T_z\otimes M=(Y'_z\otimes M)\oplus(W_z\otimes M)$ as free $\mathbf k[z^{\pm1}]$-modules, which implies
\[
(T_z\otimes M)\cap\bigl(\mathbf k((z))\otimes Y'_z\otimes M\bigr)=Y'_z\otimes M.
\]
Hence, the restriction maps $M\otimes Y'_z$ into $Y'_z\otimes M$. Its specialization is the restriction of an isomorphism, hence injective, and therefore bijective by dimension counting. The case of quotients is treated similarly.

(c) By (b), we have $c_T(z)R^{\mathrm{univ}}_{M,S_z}=g(z)R^{\mathrm{ren}}_{M,S_z}$ with $g=c_T/c_{M,S}$, where both sides are defined over $\mathbf k[z^{\pm1}]$. Since some matrix entry $a(z)$ of $R^{\mathrm{ren}}_{M,S_z}$ satisfies $a(1)\ne0$, and $b(z)=g(z)a(z)$ is a Laurent polynomial, the ratio $g(z)=b(z)/a(z)$ has no pole at $z=1$. If $g$ had a zero at $z=1$, the specialization of $c_T(z)R^{\mathrm{univ}}_{M,S_z}$ would vanish, contradicting (b). Thus, $g\in\cH$, and $\mathbf r_{M,S}=g(1)^{-1}(c_T(z)R^{\mathrm{univ}}_{M,S_z})|_{z=1}$ is an isomorphism. Finally, $\mathfrak d(M,S)=0$ follows from~\ref{R:comm}.
\end{proof}

\subsubsection{Centrality of the frozen modules}\label{subsec:central}

\begin{proposition}\label{prop:central}
Let $F$ be a frozen module in $\cC$. Then:
\begin{enumerate}[label={\rm (\alph*)}]
\item $\mathfrak d(F,S)=0$ (i.e., $F\otimes S$ is simple) for every simple object $S$ in $\cC$.
\item The frozen modules pairwise strongly commute, and every module of the form $F^{\otimes\boldsymbol\lambda}:=F_1^{\otimes\lambda_1}\otimes\cdots\otimes F_{n-1}^{\otimes\lambda_{n-1}}$ with $\boldsymbol\lambda\in\Z^{n-1}_{\ge0}$ is real simple.
\item In $\cK_t(\cC)$, the relation $[F_i]_t[S]_t=t^{\Lambda(F_i,S)}[S]_t[F_i]_t$ holds for every simple object $S$.
\end{enumerate}
\end{proposition}

\begin{proof}
(a) According to~\cite[Thm.~5.16(iii),(iv)]{KKOP3}, the module $F=M[a,b]$ is real simple, and $\mathfrak d(C_s,M[a,b])=0$ whenever $a^-<s<b^+$. In our setting, $a^-=-\infty$ and $b^+=+\infty$, so $\mathfrak d(F,C_s)=0$ for all $s$. Since $\cC$ is the closure of $\{\mathbf 1,C_1,\dots,C_r\}$ under tensor products, subquotients, and extensions, every simple object in $\cC$ arises as a composition factor of a tensor product of cuspidal modules. The claim now follows from Lemma~\ref{lem:Rsub}(c).

(b) By part (a) and~\ref{R:comm}, the frozen modules pairwise strongly commute. We prove by induction on $|\boldsymbol\lambda|$ that $F^{\otimes\boldsymbol\lambda}$ is simple. Assume that $F^{\otimes\boldsymbol\mu}$ is simple and let $\boldsymbol\lambda=\boldsymbol\mu+\mathbf e_i$. Then $F_i\otimes F^{\otimes\boldsymbol\mu}$ is simple by (a), and it is isomorphic to $F^{\otimes\boldsymbol\lambda}$ via the isomorphisms $\mathbf r_{F_i,F_j}$ that reorder the factors. Applying this argument to $2\boldsymbol\lambda$ and reordering again, we obtain $F^{\otimes\boldsymbol\lambda}\otimes F^{\otimes\boldsymbol\lambda}\cong F^{\otimes2\boldsymbol\lambda}$. Since the latter is simple, $F^{\otimes\boldsymbol\lambda}$ is real.

(c) This is a direct consequence of~\cite[Lemma~6.8]{KKOP3}.
\end{proof}

\subsubsection{The category $\cC^{\mathrm s}$}

Our normalization of $R$-matrices is defined on tensor products of cuspidal modules and transported to their subquotients. We therefore work on the following subcategory.

\begin{definition}\label{def:Cs}
Let $\cC^{\mathrm s}=\cC^{\mathrm s}(\mathsf b)$ be the full subcategory of $\cC$ whose objects are the subquotients of finite direct sums of tensor products $C_{s_1}\otimes\cdots\otimes C_{s_l}$ ($l\ge0$) of affine cuspidal modules.
\end{definition}

\begin{lemma}\label{lem:Cs}
$\cC^{\mathrm s}$ is an abelian $\mathbf k$-linear monoidal category with exact tensor product, closed under subquotients, finite direct sums and tensor products. It has the same simple objects as $\cC$, so $K(\cC^{\mathrm s})=K(\cC)$ and $\cK_t(\cC^{\mathrm s})=\cK_t(\cC)$.
\end{lemma}

\begin{proof}
If $Y$ is a subquotient of $A$ and $Y'$ of $A'$, then $Y\otimes Y'$ is a subquotient of $A\otimes A'$, by exactness of $\otimes$, and $A\otimes A'$ is again a direct sum of tensor products of cuspidal modules. The statement on simple objects was shown in the proof of Proposition~\ref{prop:central}.
\end{proof}

\begin{remark}\label{rem:word}
The category $\cC^{\mathrm s}(\mathsf b)$, and hence the localization $\widetilde\cC(\mathsf b)$ below, is defined using the cuspidal modules of the fixed word $\mathbf i$; we suppress $\mathbf i$ from the notation. Only $\cC(\mathsf b)$ is known to be independent of the word~\cite[Proposition 8.21]{Qin24a},\cite[Cor.~5.34]{KKOP3}. For $\mathsf b=\Delta^{m-2}$ we always use the word $\mathbf j=u^{m-2}$.
\end{remark}

\subsubsection{Classes of normalizations}

A \emph{braider} in a monoidal category $\mathcal T$~\cite[Def.~2.1]{KKOP-loc} is a pair $(C,R_C)$ of an object $C$ and morphisms $R_C(X)\colon C\otimes X\to X\otimes C$, natural in $X$, with $R_C(\mathbf 1)=\mathrm{id}_C$ and $R_C(X\otimes Y)=(\mathrm{id}_X\otimes R_C(Y))\circ(R_C(X)\otimes\mathrm{id}_Y)$. It is \emph{central} if all $R_C(X)$ are isomorphisms. The natural candidate for $R_F(X)$ is the specialization at $z=1$ of $c(z)R^{\mathrm{univ}}_{F,X_z}$, but the factor $c(z)$ depends on $X$ and is determined only up to $\cH$. To choose it coherently, we group tensor products of cuspidal modules into classes.

Fix a frozen module $F$, and set $c_s=c_{F,C_s}(z)$ for $1\le s\le r$. For $\mathbf n\in\Z^r$, let $c_{\mathbf n}=\prod_s c_s^{n_s}$. For a tensor product $T$ of cuspidal modules, let $\mathbf n(T)\in\Z^r_{\ge0}$ denote the vector counting the occurrences of each $C_s$, and define $c_T=c_{\mathbf n(T)}$. Let $\mathcal L=\{\mathbf n\in\Z^r \mid c_{\mathbf n}\in\cH\}$, which is a subgroup of $\Z^r$, and let $\chi_0\colon\mathcal L\to\mathbf k^\times$ be the map $\mathbf n\mapsto c_{\mathbf n}(1)$. Since $\mathbf k^\times$ is divisible, $\chi_0$ extends to a homomorphism
\begin{align}\label{eq-chi}
\chi\colon\Z^r\to\mathbf k^\times,
\end{align}
which we fix henceforth. The \emph{class} of a tensor product $T$ is defined as the coset $\mathbf n(T)+\mathcal L \in \Z^r/\mathcal L$.

\begin{lemma}\label{lem:class}
\leavevmode
\begin{enumerate}[label={\rm (\alph*)}]
\item If a non-zero module is a subquotient of both $T$ and $T'$, then $T$ and $T'$ have the same class. Consequently, every simple object of $\cC$ has a well-defined class.
\item If no composition factor of $Y$ has the same class as a composition factor of $Y'$, then $\operatorname{Hom}(Y,Y')=0$.
\item Every object $X\in\cC^{\mathrm s}$ admits a unique decomposition $X=\bigoplus_{\gamma} X_{\gamma}$, where all composition factors of $X_{\gamma}$ have class $\gamma$, and $X_{\gamma}$ is a subquotient of a finite direct sum of tensor products of class $\gamma$.
\end{enumerate}
\end{lemma}

\begin{proof}
(a) Let $S$ be a simple subquotient of both $T$ and $T'$. By Lemma~\ref{lem:Rsub}(c), we have $c_T/c_{F,S} \in \cH$ and $c_{T'}/c_{F,S} \in \cH$. It follows that $c_T/c_{T'} \in \cH$, implying that $T$ and $T'$ belong to the same class.

(b) The image of any non-zero morphism from $Y$ to $Y'$ would be a common subquotient, which contradicts the assumption on their composition factors by part (a).

(c) Write $X=Y'/Y''$ with $Y''\subseteq Y'\subseteq A=\bigoplus_j T_j$. Grouping the summands $T_j$ by their classes yields a decomposition $A=\bigoplus_{\gamma} A_{\gamma}$, where, by (a) and the Jordan--H\"older theorem, all composition factors of $A_\gamma$ have class $\gamma$. We claim that every submodule $Y\subseteq A$ satisfies $Y=\bigoplus_{\gamma}(Y\cap A_{\gamma})$. By induction, it suffices to consider the case $A=A_{\gamma}\oplus A'$. Let $\pi_\gamma$ and $\pi'$ be the projections onto the respective summands, and set $Y_\gamma=Y\cap A_\gamma$ and $Y'_0=Y\cap A'$. Consider the map
\[ \pi_\gamma(Y)/Y_\gamma \to \pi'(Y)/Y'_0, \quad a+Y_\gamma \mapsto a'+Y'_0, \]
defined for $(a,a')\in Y$. This is a well-defined isomorphism. However, its source has only composition factors of class $\gamma$, while its target has none. By (b), this map must vanish, implying $Y=Y_\gamma\oplus Y'_0$. Applying this decomposition to $Y'$ and $Y''$ induces the decomposition of $X$, and uniqueness follows from (b).
\end{proof}

\subsubsection{Braiders}

By Lemma~\ref{lem:class}(a), every composition factor of an object $X\in\cC^{\mathrm s}$ has a well-defined class. We say that $X$ is \emph{of class $\gamma\in\Z^r/\mathcal L$} if all its composition factors have class $\gamma$, i.e.\ if $X=X_\gamma$ in the decomposition of Lemma~\ref{lem:class}(c). By Lemma~\ref{lem:class}(c), such an $X$ is a subquotient of a finite direct sum of tensor products of cuspidal modules of class $\gamma$.

\begin{definition}\label{def:braider}
Let $X\in\cC^{\mathrm s}$ be of class $\gamma$. Choose a finite direct sum $A$ of tensor products of class $\gamma$, a presentation $X\cong Y'/Y''$ with $Y''\subseteq Y'\subseteq A$, and a summand $T$ of $A$. Put
\begin{align}\label{eq-braider}
R_F(X)=\chi(\mathbf n(T))^{-1}\Bigl(c_T(z)\,R^{\mathrm{univ}}_{F,X_z}\Bigr)\Big|_{z=1}\colon F\otimes X\to X\otimes F .
\end{align}
For general $X$, put $R_F(X)=\bigoplus_{\gamma} R_F(X_{\gamma})$, with $X_\gamma$ as in Lemma~\ref{lem:class}(c).
\end{definition}

\begin{proposition}\label{prop:braider}
The morphisms $R_F(X)$ are well defined, independent of all choices, and are isomorphisms. They are natural in $X\in\cC^{\mathrm s}$, and satisfy $R_F(\mathbf 1)=\mathrm{id}_F$ and $R_F(X\otimes Y)=(\mathrm{id}_X\otimes R_F(Y))\circ(R_F(X)\otimes\mathrm{id}_Y)$. Hence $(F,R_F)$ is a central braider in $\cC^{\mathrm s}$ in the sense of~\cite[Def.~2.1]{KKOP-loc}.
\end{proposition}

\begin{proof}
\emph{Existence.} Every summand $T'$ of $A$ shares the same class as $T$, so $g=c_T/c_{T'}\in\cH$. By Lemma~\ref{lem:Rsub}(a), $c_{T'}R^{\mathrm{univ}}_{F,T'_z}$ is defined over $\mathbf k[z^{\pm1}]$, and thus over $\mathcal O$, specializing to an isomorphism. Consequently, $c_TR^{\mathrm{univ}}_{F,T'_z}=g\,c_{T'}R^{\mathrm{univ}}_{F,T'_z}$ is defined over $\mathcal O$ and also specializes to an isomorphism. (Note that it need not be defined over $\mathbf k[z^{\pm1}]$, as $g$ may have poles away from $z=1$.) Summing over all summands and applying the argument of Lemma~\ref{lem:Rsub}(b) over $\mathcal O$, we deduce that $c_TR^{\mathrm{univ}}_{F,X_z}$ is defined over $\mathcal O$ and specializes to an isomorphism. The specialization in~\eqref{eq-braider} is understood in this sense.

\emph{Independence.} For another choice $(A',T')$, we have $\mathbf n(T)-\mathbf n(T')\in\mathcal L$, and $(c_T/c_{T'})(1)=\chi_0(\mathbf n(T)-\mathbf n(T'))=\chi(\mathbf n(T))\chi(\mathbf n(T'))^{-1}$. It follows that
\[
\chi(\mathbf n(T))^{-1}(c_TR^{\mathrm{univ}}_{F,X_z})|_{z=1}=\chi(\mathbf n(T'))^{-1}(c_{T'}R^{\mathrm{univ}}_{F,X_z})|_{z=1}.
\]

\emph{Naturality.} By Lemma~\ref{lem:class}(b), every morphism $f\colon X\to X'$ decomposes into a sum of components $f_{\gamma}\colon X_{\gamma}\to X'_{\gamma}$. Since both $X_\gamma$ and $X'_\gamma$ are subquotients of the common module $A\oplus A'$, the independence established above allows us to compute $R_F(X_\gamma)$ and $R_F(X'_\gamma)$ using the same reference object $T$. Multiplying~\ref{R:nat} by $\chi(\mathbf n(T))^{-1}c_T(z)$ ensures that all maps are defined over $\mathcal O$. Specializing at $z=1$ then yields $(f_\gamma\otimes\mathrm{id}_F)\circ R_F(X_\gamma)=R_F(X'_\gamma)\circ(\mathrm{id}_F\otimes f_\gamma)$.

\emph{Hexagon and unit.} If the classes of $X$ and $Y$ are realized by $T$ and $T'$ respectively, then $X\otimes Y$ is realized by $T\otimes T'$, with $c_{T\otimes T'}=c_Tc_{T'}$. The identity then follows from~\ref{R:hex} and the multiplicativity of $\chi$. For the unit object, $\mathbf 1$ corresponds to the empty tensor product, so $\chi({\bf 0})=1$ and $R^{\mathrm{univ}}_{F,\mathbf 1_z}=\mathrm{id}$, since $(\mathrm{id}\otimes\varepsilon)(\mathcal R)=1$ for the universal $R$-matrix~\cite{Dri87, CP-book}.
\end{proof}

\begin{proposition}\label{prop:commfamily}
$\{(F_i,R_{F_i})\}_{1\le i\le n-1}$ is a real commuting family of braiders in $\cC^{\mathrm s}$~\cite[Def.~2.2]{KKOP-loc}, i.e.\ $R_{F_i}(F_i)\in\mathbf k^\times\mathrm{id}_{F_i\otimes F_i}$ and $R_{F_j}(F_i)\circ R_{F_i}(F_j)\in\mathbf k^\times\mathrm{id}_{F_i\otimes F_j}$.
\end{proposition}

\begin{proof}
Both maps are automorphisms of simple modules (Proposition~\ref{prop:central}(b)), whose endomorphism rings are $\mathbf k$.
\end{proof}

\subsubsection{The localized category}

For $\boldsymbol\lambda,\boldsymbol\mu\in\Z^{n-1}$ we write $\boldsymbol\lambda\ge\boldsymbol\mu$ if $\lambda_i\ge\mu_i$ for all $i$. An object $L$ of a monoidal category is \emph{invertible} if there is an object $L'$ with $L\otimes L'\cong\mathbf 1\cong L'\otimes L$.

\begin{proposition}\label{prop:loc}
Let $\widetilde\cC=\widetilde\cC(\mathsf b)=\cC^{\mathrm s}(\mathsf b)\bigl[F_1^{\otimes-1},\dots,F_{n-1}^{\otimes-1}\bigr]$ be the localization of~\cite[Thm.~2.7]{KKOP-loc} with respect to the family of Proposition~\ref{prop:commfamily}, and let $\Upsilon\colon\cC^{\mathrm s}\to\widetilde\cC$ be the canonical functor. Then:
\begin{enumerate}[label={\rm (\alph*)}]
\item $\widetilde\cC$ is a $\mathbf k$-linear abelian monoidal category with exact tensor product, $\Upsilon$ is exact and monoidal, and each $\Upsilon(F_i)$ is invertible. Moreover $(\widetilde\cC,\Upsilon)$ is universal among monoidal functors $\Psi$ from $\cC^{\mathrm s}$ such that each $\Psi(F_i)$ is invertible and each $\Psi(R_{F_i}(X))$ is an isomorphism~\cite[Thm.~2.7]{KKOP-loc}. Since the braiders $R_{F_i}$ are central, the second condition is automatic; hence $\widetilde\cC$ is, up to monoidal equivalence, independent of the choices of the coefficients $c_s$ and of the character $\chi$ made in~\eqref{eq-chi}.
\item There are invertible objects $\widetilde F^{\boldsymbol\lambda}$, $\boldsymbol\lambda\in\Z^{n-1}$, with $\widetilde F^{\boldsymbol\lambda}\cong\Upsilon(F^{\otimes\boldsymbol\lambda})$ for $\boldsymbol\lambda\ge{\bf 0}$ and $\widetilde F^{\boldsymbol\lambda}\otimes\widetilde F^{\boldsymbol\mu}\cong\widetilde F^{\boldsymbol\lambda+\boldsymbol\mu}$, and every object of $\widetilde\cC$ is isomorphic to $\Upsilon(X)\otimes\widetilde F^{\boldsymbol\lambda}$ for some $X\in\cC^{\mathrm s}$ and $\boldsymbol\lambda\in\Z^{n-1}$.
\item The simple objects of $\widetilde\cC$ are, up to isomorphism, the $\Upsilon(S)\otimes\widetilde F^{\boldsymbol\lambda}$ with $S$ simple in $\cC$ and $\boldsymbol\lambda\in\Z^{n-1}$, and $\Upsilon(S)\otimes\widetilde F^{\boldsymbol\lambda}\cong\Upsilon(S')\otimes\widetilde F^{\boldsymbol\lambda'}$ iff $F^{\otimes(\boldsymbol\delta+\boldsymbol\lambda)}\otimes S\cong F^{\otimes(\boldsymbol\delta+\boldsymbol\lambda')}\otimes S'$ for some $\boldsymbol\delta$ with $\boldsymbol\delta+\boldsymbol\lambda,\boldsymbol\delta+\boldsymbol\lambda'\ge{\bf 0}$. The objects $\Upsilon(S)\otimes\widetilde F^{\boldsymbol\lambda}$ with $S$ real are real.
\item Every object of $\widetilde\cC$ has finite length, and there is a ring isomorphism
\begin{align}\label{eq-Kloc}
K(\cC)\bigl[[F_1]^{-1},\dots,[F_{n-1}]^{-1}\bigr]\xrightarrow{\ \sim\ }K(\widetilde\cC),\qquad [S]\,[F]^{\boldsymbol\lambda}\mapsto[\Upsilon(S)\otimes\widetilde F^{\boldsymbol\lambda}].
\end{align}
\end{enumerate}
\end{proposition}

\begin{proof}
We recall the construction given in~\cite[\S2]{KKOP-loc}. The objects of $\widetilde\cC$ are pairs $(X,\boldsymbol\lambda)$ with $X\in\cC^{\mathrm s}$ and $\boldsymbol\lambda\in\Z^{n-1}$, and the morphism spaces are defined by
\[
\operatorname{Hom}_{\widetilde\cC}\bigl((X,\boldsymbol\lambda),(Y,\boldsymbol\mu)\bigr)=\varinjlim_{\boldsymbol\delta}\operatorname{Hom}_{\cC^{\mathrm s}}\bigl(F^{\otimes(\boldsymbol\delta+\boldsymbol\lambda)}\otimes X,\ Y\otimes F^{\otimes(\boldsymbol\delta+\boldsymbol\mu)}\bigr),
\]
where the colimit is taken over $\boldsymbol\delta$ such that $\boldsymbol\delta+\boldsymbol\lambda,\boldsymbol\delta+\boldsymbol\mu\ge{\bf 0}$, with transition maps induced by the braiding~\cite[(2.6), \S2.2--2.3]{KKOP-loc}. The tensor product is given by $(X,\boldsymbol\lambda)\otimes(Y,\boldsymbol\mu)=(X\otimes Y,\boldsymbol\lambda+\boldsymbol\mu)$, the embedding by $\Upsilon(X)=(X,{\bf 0})$, and we have $(\mathbf 1,\boldsymbol\lambda)\cong\Upsilon(F^{\otimes\boldsymbol\lambda})$ for $\boldsymbol\lambda\ge{\bf 0}$~\cite[Prop.~2.6]{KKOP-loc}. The universality is established in~\cite[Thm.~2.7]{KKOP-loc}, while the abelianness and exactness follow from~\cite[Prop.~2.9]{KKOP-loc}. This proves (a) and (b), identifying $\widetilde F^{\boldsymbol\lambda}$ with $(\mathbf 1,\boldsymbol\lambda)$.

(c) Let $S$ be a simple object. For any $\boldsymbol\nu\ge{\bf 0}$, the module $F^{\otimes\boldsymbol\nu}\otimes S\cong S\otimes F^{\otimes\boldsymbol\nu}$ is simple by Proposition~\ref{prop:central}. In particular, the transition maps in the colimit defining $\operatorname{End}_{\widetilde\cC}((S,\boldsymbol\lambda))$ are injective, since the functor $F^{\otimes\boldsymbol\delta}\otimes-$ is exact and faithful on finite-dimensional modules; consequently, $(S,\boldsymbol\lambda)\ne0$. A non-zero morphism $g\colon(Z,\boldsymbol\mu)\to(S,\boldsymbol\lambda)$ is represented by a non-zero morphism $f$ onto the simple module $S\otimes F^{\otimes(\boldsymbol\delta+\boldsymbol\lambda)}$, which is necessarily surjective. Moreover, $g$ corresponds to $\Upsilon(f)\otimes\mathrm{id}_{\widetilde F^{-\boldsymbol\delta}}$. Since $\Upsilon$ is exact and ${-}\otimes\widetilde F^{-\boldsymbol\delta}$ is an exact autoequivalence, $g$ is an epimorphism, proving that $(S,\boldsymbol\lambda)$ is simple.

Conversely, a composition series of $X$ in $\cC^{\mathrm s}$ induces a finite filtration of $\Upsilon(X)\otimes\widetilde F^{\boldsymbol\lambda}$ with subquotients of the form $\Upsilon(S_j)\otimes\widetilde F^{\boldsymbol\lambda}$. Thus, every simple object is of the stated form, and every object has finite length. A non-zero morphism between $(S,\boldsymbol\lambda)$ and $(S',\boldsymbol\lambda')$ is represented by a non-zero morphism between the simple modules $F^{\otimes(\boldsymbol\delta+\boldsymbol\lambda)}\otimes S$ and $F^{\otimes(\boldsymbol\delta+\boldsymbol\lambda')}\otimes S'$, which must be an isomorphism by Schur's lemma. Finally, we observe that $(\Upsilon(S)\otimes\widetilde F^{\boldsymbol\lambda})^{\otimes2}\cong\Upsilon(S\otimes S)\otimes\widetilde F^{2\boldsymbol\lambda}$.

(d) By~\eqref{eq-qchar}, $K(\cC)$ is a commutative domain that is free on the classes of simple objects. Moreover, multiplication by $[F_i]$ maps simple classes injectively to simple classes (Proposition~\ref{prop:central}). Hence, the localized ring $K(\cC)[[F]^{-1}]$ is free on the elements $[S][F]^{\boldsymbol\lambda}$, subject to the relations $[F^{\otimes\boldsymbol\mu}\otimes S][F]^{\boldsymbol\lambda}=[S][F]^{\boldsymbol\lambda+\boldsymbol\mu}$. These elements correspond precisely to the classes of simple objects in $\widetilde\cC$ by part (c). The resulting map is multiplicative because $\Upsilon$ is monoidal and exact.
\end{proof}

\subsubsection{The quantum Grothendieck ring of $\widetilde\cC$}\label{subsec:Ktloc}

The $(q,t)$-character of a frozen module is not invertible in $\cY_t$, so the $(q,t)$-characters of the objects of $\widetilde\cC$ live in the skew-field of fractions $\cD_t$ of $\cY_t$. The bar involution extends uniquely to an anti-automorphism of $\cD_t$. For $\boldsymbol\lambda\in\Z^{n-1}$ put $[F]_t^{\boldsymbol\lambda}=[F_1]_t^{\lambda_1}\cdots[F_{n-1}]_t^{\lambda_{n-1}}\in\cD_t$. By Proposition~\ref{prop:central}(c), the $[F_i]_t$ pairwise quasi-commute, and $[F]_t^{\boldsymbol\lambda}[S]_t\in t^{\frac12\Z}[S]_t[F]_t^{\boldsymbol\lambda}$ for every simple $S\in\cC$.

\begin{lemma}\label{lem:qtloc}
Let $X\cong\Upsilon(S)\otimes\widetilde F^{\boldsymbol\lambda}$ be a simple object of $\widetilde\cC$, with $S$ simple in $\cC$. There exists a unique $c\in\tfrac12\Z$ such that
\begin{equation}\label{eq-qtloc}
[X]_t:=t^{c}[S]_t[F]_t^{\boldsymbol\lambda}\in\cD_t
\end{equation}
is bar-invariant, and $[X]_t$ depends only on the isomorphism class of $X$. If $\boldsymbol\lambda\ge{\bf 0}$, then $[X]_t$ coincides with the $(q,t)$-character of the simple module $S\otimes F^{\otimes\boldsymbol\lambda}$; in particular, $[\Upsilon(S)]_t=[S]_t$.
\end{lemma}

\begin{proof}
Let $x=[S]_t[F]_t^{\boldsymbol\lambda}$. Since the bar involution is an anti-automorphism that fixes $[S]_t$ and $[F_i]_t$, we have $\overline{x}=t^e x$ for some $e\in\Z$. Consequently, $t^c x$ is bar-invariant if and only if $c=e/2$.

If $\boldsymbol\lambda\ge{\bf 0}$, the $(q,t)$-character of the simple module $S\otimes F^{\otimes\boldsymbol\lambda}$ is a bar-invariant element of the form $t^a[S]_t[F]_t^{\boldsymbol\lambda}$~\cite[\S6]{KKOP3}, and thus it must coincide with $[X]_t$.

To see that $[X]_t$ is well-defined, suppose that $(S,\boldsymbol\lambda)$ and $(S',\boldsymbol\lambda')$ present isomorphic objects in $\widetilde\cC$. By Proposition~\ref{prop:loc}(c), the elements $[S]_t[F]_t^{\boldsymbol\delta+\boldsymbol\lambda}$ and $[S']_t[F]_t^{\boldsymbol\delta+\boldsymbol\lambda'}$ agree up to a power of $t^{1/2}$. It follows that $[S]_t[F]_t^{\boldsymbol\lambda}$ and $[S']_t[F]_t^{\boldsymbol\lambda'}$ also differ only by a power of $t^{1/2}$, implying that their bar-invariant normalizations coincide.
\end{proof}

\begin{definition}\label{def:Ktloc}
We define the \emph{quantum Grothendieck ring} $\cK_t(\widetilde\cC)$ of $\widetilde\cC$ to be the $\Z[t^{\pm1/2}]$-subalgebra of $\cD_t$ generated by the $(q,t)$-characters $[X]_t$ associated with all simple objects $X$ in $\widetilde\cC$.
\end{definition}

The next theorem is the main result of this subsection. Part (c) identifies $\cK_t(\widetilde\cC)$ with the quantum cluster algebra $\mathscr A(\mathsf s(\mathbf i))$, parts (a) and (b) give the positive bases of Corollary~\ref{cor:positivity} and Theorem~\ref{thm:Aquant}, and part (d) gives the classical limits in Theorem~\ref{thm:main}.

\begin{theorem}\label{thm:Ktloc}
\leavevmode
\begin{enumerate}[label={\rm (\alph*)}]
\item The $(q,t)$-characters $[X]_t$ of the simple objects of $\widetilde\cC$, taken over isomorphism classes, form a $\Z[t^{\pm1/2}]$-basis of $\cK_t(\widetilde\cC)$.
\item For simple objects $X,Y$ of $\widetilde\cC$, the product decomposes as $[X]_t[Y]_t=\sum_Zc_{X,Y}^Z(t)\,[Z]_t$ with $c_{X,Y}^Z(t)\in\Z_{\ge0}[t^{\pm1/2}]$, where $c^Z_{X,Y}(1)$ equals the multiplicity of $Z$ as a composition factor of $X\otimes Y$.
\item We have $\cK_t(\widetilde\cC)=\cK_t(\cC)\bigl[[F_1]_t^{-1},\dots,[F_{n-1}]_t^{-1}\bigr]$, which is the Ore localization of $\cK_t(\cC)$ at the frozen classes. Consequently,
\begin{align}\label{eq-Ktloc-cluster}
\cK_t(\widetilde\cC)\otimes_{\Z[t^{\pm1/2}]}R=\mathscr A(\mathsf s(\mathbf i)) ,
\end{align}
viewed as subalgebras of $\cD_t^R$. The natural map $\cK_t(\widetilde\cC)\otimes_{\Z[t^{\pm1/2}]}R\to\cD_t^R$ is injective, since $R$ is a free $\Z[t^{\pm1/2}]$-module and every element of $\cK_t(\widetilde\cC)$ can be mapped into $\cK_t(\cC)\subseteq\cY_t$ by right multiplication by some $[F]_t^{\boldsymbol\delta}$.
\item The $\Z$-linear map $\cK_t(\widetilde\cC)\to K(\widetilde\cC)$ defined by $t^{a}[X]_t\mapsto[X]$ is a ring homomorphism. It induces an isomorphism $\cK_t(\widetilde\cC)\otimes_{\Z[t^{\pm1/2}]}\Z\cong K(\widetilde\cC)$ via $t^{1/2}\mapsto1$.
\end{enumerate}
\end{theorem}

\begin{proof}
(b) Let $X\cong\Upsilon(S)\otimes\widetilde F^{\boldsymbol\lambda}$ and $Y\cong\Upsilon(S')\otimes\widetilde F^{\boldsymbol\mu}$. By the quasi-commutation relations and~\eqref{eq-positivity}, we have
\[
[X]_t[Y]_t \in t^{\frac12\Z}\sum_{S''}c^{S''}_{S,S'}(t)\,[S'']_t[F]_t^{\boldsymbol\lambda+\boldsymbol\mu}.
\]
By Lemma~\ref{lem:qtloc}, each term $[S'']_t[F]_t^{\boldsymbol\lambda+\boldsymbol\mu}$ is a power of $t^{1/2}$ times $[Z'']_t$, where $Z''=\Upsilon(S'')\otimes\widetilde F^{\boldsymbol\lambda+\boldsymbol\mu}$. In $\widetilde\cC$, we have $X\otimes Y\cong\Upsilon(S\otimes S')\otimes\widetilde F^{\boldsymbol\lambda+\boldsymbol\mu}$, so its composition factors are precisely the $Z''$ with multiplicities $c^{S''}_{S,S'}(1)$. Moreover, $Z''\cong Z'''$ implies $S''\cong S'''$ by Proposition~\ref{prop:loc}(c) and the injectivity of $\chi_q$.

(a) Part (b) shows that the span of the $[X]_t$ is closed under multiplication, and thus coincides with $\cK_t(\widetilde\cC)$. For linear independence, consider a relation $\sum_jp_j(t)[X_j]_t=0$. Multiplying on the right by $[F]_t^{\boldsymbol\delta}$ for sufficiently large $\boldsymbol\delta$ yields a relation among the $(q,t)$-characters of pairwise non-isomorphic simple modules in $\cC$. These characters are known to be linearly independent~\cite[\S6.1]{KKOP3}, which proves the claim.

(c) By Proposition~\ref{prop:central}(c), the elements $[F_i]_t$ are normal in $\cK_t(\cC)$ and invertible in $\cD_t$. Hence, the inclusion $\cK_t(\cC)\hookrightarrow\cD_t$ extends to an injective homomorphism from the Ore localization. The image of this homomorphism is spanned by the elements $[S]_t[F]_t^{\boldsymbol\lambda}$, which coincide with the $[X]_t$ up to powers of $t^{1/2}$. By (a), this image is exactly $\cK_t(\widetilde\cC)$. Equation~\eqref{eq-Ktloc-cluster} then follows from \eqref{eKKOP-iso} and Lemma~\ref{lem:normal}, noting that tensoring with the free $\Z[t^{\pm1/2}]$-module $R$ commutes with Ore localization.

(d) By (a), $\cK_t(\widetilde\cC)\otimes\Z$ is free on the images of the $[X]_t$. The map sends these images bijectively to the classes of simple objects, which form a basis of $K(\widetilde\cC)$ by Proposition~\ref{prop:loc}(c). The multiplicativity of the map follows from (b).
\end{proof}

In particular, every frozen variable of $\mathsf s(\mathbf i)$ and its inverse coincides, up to a power of $t^{1/2}$, with the $(q,t)$-character of an invertible object $\widetilde F^{\pm\mathbf e_i}$ in $\widetilde\cC$.

\begin{remark}[Other types]\label{rem:ADE}
The constructions and proofs presented in this subsection apply verbatim to the braid-word categories $\cC(\mathsf b)$ associated with any untwisted quantum affine algebra of simply-laced type $ADE$, provided that every color occurs in the chosen word.
\end{remark}

\section{Monoidal categorification of the reduced skein algebra $\overline{\mathscr S}(\mathbb P_m)$}\label{sec-cat}

In this section, we present the proofs of Theorems~\ref{thmA} and~\ref{thmB}. In \S\ref{subsec-ident}, we establish that the sub-seed $\mathsf s^\qc_\core$ of the skein seed (Definition~\ref{def:sqc}) coincides with the KKOP seed $\mathsf s(u^{m-2})$ from Definition~\ref{def:KKOPseed}. Specifically, we show that the exchange matrices coincide via an explicit comparison of arrows (Lemma~\ref{lem:M}), while the quantum forms are shown to agree through a computation on $n-1$ vertices (Lemmas~\ref{lem:XiK} and~\ref{lem:XiS}). In \S\ref{subsec-loc-Delta}, we combine these results with the localization framework from Section~\ref{sec-KKOP} for $\mathsf b=\Delta^{m-2}$ to identify $\mathcal A_{n,m}$ with the quantum Grothendieck ring of $\widetilde\cC(\Delta^{m-2})$. Subsequently, in \S\ref{subsec-CI}, we adjoin the frozen variables associated with the edges $(1,2),\dots,(m-1,m)$ as invertible graded objects and prove the categorification theorem (Theorem~\ref{thm:main}). The results obtained in \S\ref{subsec-loc-Delta} will also be utilized in Section~\ref{sec-braid}.

Throughout this section, we fix integers $n\ge2$ and $m\ge3$, and we adopt the notation introduced in \S\ref{subsec-CHW} and \S\ref{subsec-qc}. This includes the vertex set $V=\Vc\sqcup \Viso$ (Definition~\ref{def:core}), the seeds $\mathsf s_{\rm skein}$ (Proposition~\ref{prop:CHWseed}) and $\mathsf s^\qc_{\rm skein}=(Q^\qc,\Lambda^\qc,M^\qc)$ (Definition~\ref{def:sqc}), and the subalgebra $\mathcal A_{n,m}$ (Definition~\ref{def:Anm}). We set
\begin{align}\label{eq-j-word}
\mathbf j=u^{m-2},\qquad r=(m-2)\tbinom n2,
\end{align}
where $u$ denotes the word defined in~\eqref{eq-u}. We label the positions of $\mathbf j$ by triples $\langle\kappa,r',c\rangle$ as described in \S\ref{subsec-wordu}, with $1\le\kappa\le m-2$.

\subsection{Identification of ${\mathsf s}^{\rm qc}_{\rm core}$ with the KKOP seed ${\mathsf s}(u^{m-2})$}\label{subsec-ident}

In this subsection, we prove that the sub-seed $\mathsf s^\qc_\core$ of the skein seed $\mathsf s^\qc_{\rm skein}$ coincides with the KKOP seed $\mathsf s(\mathbf j)$ of the word $\mathbf j=u^{m-2}$ of $\Delta^{m-2}$, and that the exchange matrix of $\mathsf s^\qc_{\rm skein}$ is that of $\mathsf s(\mathbf j)$ extended by $(m-1)(n-1)$ isolated frozen vertices.

\subsubsection{The exchange matrices}

The following lemma follows directly from Definition~\ref{def:CHWquiver} and Lemma~\ref{lem:sideedges}.

\begin{lemma}\label{lem:corequiver}
The arrows of the quiver $Q$ which join a vertex of $\Vc$ and a mutable vertex are exactly the arrows of the following three types, where $1\le i\le m-2$, and a vertex $(j,0)^i$ with $i\ge2$ stands for $(j,j)^{i-1}$:
\begin{align*}
&\text{type (A1):}\quad (j+1,s)^i\to(j,s)^i, && 1\le s\le j\le n-2;\\
&\text{type (A2):}\quad (j,s+1)^i\to(j,s)^i, && 0\le s\le j-1\le n-2,\ \text{ with } s\ge1 \text{ if } i=1;\\
&\text{type (A3):}\quad (j-1,s-1)^i\to(j,s)^i, && 1\le s\le j-1\le n-2,\ \text{ with } s\ge2 \text{ if } i=1.
\end{align*}
All of them have weight $1$, and no two of them join the same pair of vertices. In particular, $Q(x,y)=\#\{x\to y\}-\#\{y\to x\}$ for $x\in\Vc$ and $y\in V_{\rm mut}$, where only arrows of types (A1)--(A3) are counted.
\end{lemma}

\begin{lemma}\label{lem:M}
Define $\psi\colon\Vc\to[1,r]$ by
\begin{align}\label{eq-psi}
\psi\bigl((j,s)^i\bigr)=\langle i,\ s,\ n-j\rangle .
\end{align}
Then $\psi$ is a bijection onto the positions of $\mathbf j$, it maps the frozen vertices $(j,j)^{m-2}$ of $\Vc$ to the last occurrences, and
\begin{align}\label{eq-thmM}
Q(x,y)=Q^\qc(x,y)=Q_{\mathbf j}\bigl(\psi(x),\psi(y)\bigr)\qquad(x\in\Vc,\ y\in V_{\rm mut}),
\end{align}
where $Q_{\mathbf j}$ is the exchange matrix of $\mathsf s(\mathbf j)$ (Definitions~\ref{def:BFZ} and~\ref{def:KKOPseed}).
\end{lemma}

\begin{proof}
The first equality in~\eqref{eq-thmM} is Proposition~\ref{prop:sqc}(a). The triple $\langle i,s,n-j\rangle$ is a position, i.e.\ $1\le s\le n-1$ and $1\le n-j\le n-s$, iff $1\le s\le j\le n-1$. Hence $\psi$ is a bijection, with inverse $\langle\kappa,r',c\rangle\mapsto(n-c,\,r')^{\kappa}$; it identifies the triangle $\Delta_i$ with the $i$-th copy of $u$, the coordinate $s$ with the block, and $j$ with the colour $n-j$. By Lemma~\ref{lem:Quk}, the last occurrence of the colour $n-j$ is $\langle m-2,j,n-j\rangle=\psi((j,j)^{m-2})$.

Under $\psi$, the arrows in Lemma~\ref{lem:corequiver} correspond to those in Lemma~\ref{lem:Quk}. Specifically, by~\eqref{eq-Quk-next} and $(j,j)^{i}=(j,0)^{i+1}$, the type (A2) arrows correspond to those in Lemma~\ref{lem:Quk}(a). Furthermore, a direct comparison of the index ranges shows that the type (A1) and (A3) arrows correspond to those in Lemma~\ref{lem:Quk}(b) and~(c), respectively. All these arrows have weight $1$, and neither quiver contains a pair of vertices joined twice or in both directions (Lemmas~\ref{lem:corequiver} and~\ref{lem:Q2unique}). This proves~\eqref{eq-thmM}.
\end{proof}

\subsubsection{The cokernel of the exchange matrix}

Let
\begin{align}\label{eq-Lst}
\Lst=\{x_j\mid 1\le j\le n-1\}\subseteq\Vc,\qquad x_j:=(j,1)^1 .
\end{align}
By~\eqref{eq-psi}, $\psi(x_j)=\langle 1,1,n-j\rangle$ is the first occurrence of the colour $n-j$ in $\mathbf j$, so $\psi$ maps $\Lst$ bijectively onto the first occurrences; this is the reason for the subscript ``first''. Let $L_\core\subseteq\Z^{\Vc}$ be the span of the column vectors ${\bf q}^\qc_y$, $y\in V_{\rm mut}$, of $Q^\qc$; they are supported on $\Vc$ by Proposition~\ref{prop:sqc}(a).

\begin{lemma}\label{lem:Lst}
We have $\Z^{\Vc}=L_\core\oplus\Z^{\Lst}$. Consequently, the quotient $\Z^{\Vc}/L_\core$ is identified with $\Z^{\Lst}$, and we denote by
\begin{align}\label{eq-p-Lst}
p\colon\Z^{\Vc}\to\Z^{\Lst}
\end{align}
the projection along $L_\core$. Under $\psi$, $p$ is the projection $p_{\mathbf j}$ of~\eqref{eq-Gamma-i}.
\end{lemma}

\begin{proof}
By Lemma~\ref{lem:M}, for every mutable $y$, $\psi$ maps the column vector ${\bf q}^\qc_y$ to the column vector of $Q_{\mathbf j}$ at $\psi(y)$, and it maps $\Lst$ onto the set of first occurrences of $\mathbf j$. Hence the lemma follows from Lemma~\ref{lem:unimod}, applied to the word $\mathbf j$.
\end{proof}

\subsubsection{The quantum forms on $\Lst$}

The next two lemmas compute the quantum forms of $\mathsf s(\mathbf j)$ and of $\mathsf s^\qc_{\rm skein}$ on the vertices of $\Lst$. For $n=2$ the set $\Lst$ has one element, and there is nothing to compute.

\begin{lemma}\label{lem:XiK}
For $1\le j<j'\le n-1$, we have
\[
\Pi_{\mathbf j}\bigl(\psi(x_j),\psi(x_{j'})\bigr)=n\Lambda\bigl(\mathsf M_{\psi(x_j)},\mathsf M_{\psi(x_{j'})}\bigr)=-n,
\]
where $\Pi_{\mathbf j}$ is defined in~\eqref{eKKOPseed} and $\mathsf M_s$ in~\eqref{eq-Ms}.
\end{lemma}

\begin{proof}
The first block $(1,2,\dots,n-1)$ of the word $u$ occupies the first $n-1$ positions of $\mathbf j$. Consequently, $\psi(x_j)$ corresponds to the position $c=n-j$ of colour $c$, with the associated prefix $w_{\le c}=s_1s_2\cdots s_c$. Observing that $s_i\varpi_c=\varpi_c$ for $i\ne c$ and $s_1\cdots s_{c-1}\alpha_c=\alpha_1+\dots+\alpha_c=\varepsilon_1-\varepsilon_{c+1}$, we obtain
\[
\varpi_c-w_{\le c}\varpi_c=\varepsilon_1-\varepsilon_{c+1} \quad \text{and} \quad \varpi_c+w_{\le c}\varpi_c=2\varpi_c-\varepsilon_1+\varepsilon_{c+1}.
\]
For indices $s<s'\le n-1$, applying formula~\eqref{eq-Lambda-formula} yields
\[
\Lambda_{\mathbf j}(s,s')=-\bigl(\varepsilon_1-\varepsilon_{s+1},\ 2\varpi_{s'}-\varepsilon_1+\varepsilon_{s'+1}\bigr)=-\bigl(\varepsilon_1-\varepsilon_{s+1},\ -\varepsilon_1+\varepsilon_{s'+1}\bigr)=1,
\]
where we have used the facts that $(\alpha_i,\varpi_{s'})=0$ for $i\le s<s'$ and $(\varepsilon_a,\varepsilon_b)=\delta_{ab}-\frac1n$. Since $j<j'$ implies $n-j>n-j'$, the skew-symmetry of $\Lambda_{\mathbf j}$ gives $\Lambda_{\mathbf j}(n-j,n-j')=-\Lambda_{\mathbf j}(n-j',n-j)=-1$.
\end{proof}

On the skein side, we employ~\eqref{eq-PDelta}. For any two small vertices $p,p'$ of a triangle, we write $P_\Delta(p,p')=n\,\varsigma(p,p')$. Let $P_x=(x,n-x,0)$ and $Q_y=(0,n-y,y)$ denote the barycentric coordinates of the small vertices on the sides $c=0$ and $a=0$, respectively. Note that the corner cases $P_0, P_n, Q_0, Q_n$ are treated as zero. A direct computation yields
\begin{align}\label{eq:green}
\varsigma(P_x,P_{x'})=\varsigma(Q_y,Q_{y'})=0,\qquad \varsigma(P_x,Q_y)=-\varsigma(Q_y,P_x)=g(x,y):=\min(x,y)\bigl(n-\max(x,y)\bigr).
\end{align}
Indeed, the terms $\varsigma(P_x,P_{x'})$ and $\varsigma(Q_y,Q_{y'})$ vanish as computed from the coordinates $(n-x,0),(n-x',0)$ and $(0,n-y),(0,n-y')$. The term $\varsigma(P_x,Q_y)$ is determined by the first two coordinates as $x(n-y)$ when $x\le y$, and by the last two coordinates as $y(n-x)$ when $y<x$. The kernel $g$ is clearly symmetric.

\begin{lemma}\label{lem:XiS}
For $1\le j<j'\le n-1$, $\Lambda^\qc({\bf e}_{x_j},{\bf e}_{x_{j'}})=-n$. More precisely, writing $X_j=\mathsf G_{\rm iso}({\bf e}_{x_j})$, in $n\Lambda^\qc({\bf e}_{x_j},{\bf e}_{x_{j'}})=\sum_{l=1}^{m-2}P_{\Delta_l}\bigl(\mathrm{sk}_{\Delta_l}(X_j),\mathrm{sk}_{\Delta_l}(X_{j'})\bigr)$ the term of $\Delta_1$ equals $-n^2$, and all other terms vanish.
\end{lemma}

\begin{proof}
By~\eqref{eq-LambdaCHW} and Definition~\ref{def:sqc}, $n\Lambda^\qc({\bf e}_{x_j},{\bf e}_{x_{j'}})=P(X_j,X_{j'})$, which is the displayed sum by~\eqref{eq-P}. By~\eqref{eq-Pnatural}, $X_j={\bf e}_{x_j}-\deg_{\rm iso}({\bf e}_{x_j})$ with $\deg_{\rm iso}({\bf e}_{x_j})\in\Z^{\Viso}$, and by~\eqref{eq-Pnatural-degree}, $\operatorname{res}\varrho(X_j)=\theta(\deg_{\rm iso}(X_j))=0$, i.e.\ the components of $\varrho(X_j)$ on $E_0,\dots,E_{m-2}$ vanish. We read off the skeletons from \S\ref{subsec-CHW}. The vertex $x_j$, with coordinates $(n-j,j-1,1)$, has the skeleton point $(n-j,j-1,1)$ in $\Delta_1$ and, if $m\ge4$, the point $P_{n-1}$ in $\Delta_2$. A vertex $(j'',0)^1\in \Viso$ has the single point $P_{n-j''}$ in $\Delta_1$. A vertex $(n,s)^i\in \Viso$ has the point $Q_s$ in $\Delta_i$ and, if $i<m-2$, the point $P_{n-s}$ in $\Delta_{i+1}$. By Lemma~\ref{lem:legwt}, a skeleton point $P_x$ in $\Delta_l$, $l\ge2$, comes from a leg of weight $n-x$ ending on $E_l$, and a point $Q_y$ in $\Delta_l$ from a leg of weight $n-y$ ending on $E_l$.

\emph{The triangles $\Delta_l$, $l\ge2$.} Write $\mathrm{sk}_{\Delta_l}(X_j)=\sum_x\mu_j(x)P_x+\sum_y\nu_j(y)Q_y$. Every leg ending on $E_l$ has a skeleton point in $\Delta_l$, so the $E_l$-component of $\varrho(X_j)$ is $\sum_x\bigl(\mu_j(x)+\nu_j(x)\bigr)f(n-x)$. It vanishes, and the $f(w)$ form a basis of $\mathsf P$, so $\nu_j=-\mu_j$. By~\eqref{eq:green},
\[
\varsigma\bigl(\mathrm{sk}_{\Delta_l}X_j,\mathrm{sk}_{\Delta_l}X_{j'}\bigr)=\sum_{x,y}\bigl(-\mu_j(x)\mu_{j'}(y)+\mu_j(y)\mu_{j'}(x)\bigr)g(x,y)=0,
\]
by the symmetry of $g$.

\emph{The triangle $\Delta_1$.} Here $\mathrm{sk}_{\Delta_1}(X_j)$ consists of $(n-j,j-1,1)$ with coefficient $1$, of points $P_x$ from the vertices $(n-x,0)^1$, and of points $Q_y$ from the vertices $(n,y)^1$. The legs on $E_0$ are the leg of weight $n-j$ of $x_j$ and the legs of weight $n-j''$ of the $(j'',0)^1$, so balance on $E_0$ forces the coefficient of $(j,0)^1$ to be $-1$ and the others to vanish. The legs on $E_1$ are the leg of weight $j-1$ of $x_j$, the leg of weight $j$ of $(j,0)^1$, and the legs of weight $n-y$ of the $(n,y)^1$; so balance on $E_1$ forces the coefficients of $Q_{n-j+1}$ and $Q_{n-j}$ to be $-1$ and $+1$, and the others to vanish. Hence
\begin{align}\label{eq-skDelta1}
\mathrm{sk}_{\Delta_1}(X_j)=A_j-B_j-D_j+E_j,
\end{align}
where $A_j=(n-j,j-1,1)$, $B_j=(n-j,j,0)$, $D_j=(0,j-1,n-j+1)$ and $E_j=(0,j,n-j)$,
with $D_1=Q_n$ read as $0$. For $j<j'$ the values of $\varsigma$ are as follows; in each case the indicated pair of coordinates is comparable:
\[
\begin{array}{c|cccc}
\varsigma & A_{j'} & B_{j'} & D_{j'} & E_{j'}\\\hline
A_j & j-j' & n-j' & (j-1)(n-j'+1)-(j'-1) & (j-1)(n-j')-j'\\
B_j & j & 0 & j(n-j'+1) & j(n-j')\\
D_j & -(j-1)(n-j') & -(j-1)(n-j') & 0 & 0\\
E_j & -j(n-j') & -j(n-j') & 0 & 0
\end{array}
\]
The entries are computed from the coordinates $(b,c)$ for the pairs $(A_j,A_{j'})$, $(A_j,D_{j'})$, $(A_j,E_{j'})$ and the row $B_j$; from $(c,a)$ for $(A_j,B_{j'})$; and from $(a,b)$ for the rows $D_j$ and $E_j$. With the signs $+,-,-,+$ of $A,B,D,E$, the four rows contribute
\[
\bigl(j-j'\bigr)-\bigl(n-j'\bigr)-\bigl((j-1)(n-j'+1)-(j'-1)\bigr)+\bigl((j-1)(n-j')-j'\bigr)=-n,
\]
$-\bigl(j-0-j(n-j'+1)+j(n-j')\bigr)=0$, $-\bigl(-(j-1)(n-j')+(j-1)(n-j')\bigr)=0$ and $-j(n-j')+j(n-j')=0$. So $\varsigma=-n$, and the term of $\Delta_1$ is $n\varsigma=-n^2$.
\end{proof}

\subsubsection{The identification}

We now combine the comparison of the exchange matrices (Lemma~\ref{lem:M}) with the comparison of the quantum forms on $\Lst$ (Lemmas~\ref{lem:XiK} and~\ref{lem:XiS}). Theorem~\ref{thm:seedmatch}, which is Theorem~\ref{thmA} of the introduction, asserts that the sub-seed $\mathsf s^\qc_\core$ is the KKOP seed $\mathsf s(\mathbf j)$. Corollary~\ref{cor:seedmatch} then establishes an isomorphism of quantum cluster algebras that matches cluster monomials. This serves as the bridge between the skein side and the representation side: specifically, in \S\ref{subsec-loc-Delta}, it identifies $\mathcal A_{n,m}$ with the quantum Grothendieck ring of $\widetilde\cC$; in \S\ref{subsec-CI}, it provides the categorification of $\overline{\cS}(\bP_m)$; and in Section~\ref{sec-braid}, it yields the quantum cluster structure of $\mathcal A_{n,m}$.

\begin{theorem}\label{thm:seedmatch}
For $x,x'\in\Vc$ and $y\in V_{\rm mut}$,
\begin{align}\label{eq-seedmatch}
Q^\qc(x,y)=Q_{\mathbf j}\bigl(\psi(x),\psi(y)\bigr),\qquad \Lambda^\qc(x,x')=\Pi_{\mathbf j}\bigl(\psi(x),\psi(x')\bigr),
\end{align}
and $Q^\qc(g,y)=0$ for $g\in \Viso$. In other words, under the relabelling $\psi$, the sub-seed $\mathsf s^\qc_\core$ of $\mathsf s^\qc_{\rm skein}$ and the KKOP seed $\mathsf s(\mathbf j)$ have the same exchange matrix and the same quantum form, and the exchange matrix of $\mathsf s^\qc_{\rm skein}$ is obtained from that of $\mathsf s^\qc_\core$ by adjoining the $(m-1)(n-1)$ isolated frozen vertices of $\Viso$.
\end{theorem}

\begin{proof}
The statements on exchange matrices are Lemma~\ref{lem:M} and Proposition~\ref{prop:sqc}(a). Let $D=\Lambda^\qc|_{\Vc\times\Vc}-\psi^*\Pi_{\mathbf j}$, an anti-symmetric form on $\Z^{\Vc}$, where $(\psi^*\Pi_{\mathbf j})(x,x')=\Pi_{\mathbf j}(\psi(x),\psi(x'))$. By Proposition~\ref{prop:sqc}(c) and Theorem~\ref{thm:KKOP}(a), both $\Lambda^\qc|_{\Vc\times\Vc}$ and $\psi^*\Pi_{\mathbf j}$ are compatible with the same exchange matrix and the same scalar $2n$. Hence $D({\bf q}^\qc_y,{\bf e}_x)=2n\delta_{y,x}-2n\delta_{y,x}=0$ for all $y\in V_{\rm mut}$ and $x\in\Vc$, so $D$ vanishes as soon as one of its arguments lies in $L_\core$. By Lemma~\ref{lem:Lst}, it suffices to show $D({\bf e}_{x_j},{\bf e}_{x_{j'}})=0$ for $j<j'$, which follows from Lemmas~\ref{lem:XiK} and~\ref{lem:XiS}.
\end{proof}

\begin{corollary}\label{cor:seedmatch}
There is a unique isomorphism of $R$-algebras
\begin{align}\label{eq-PhiA}
\Phi_{\mathcal A}\colon\mathscr A(\mathsf s(\mathbf j))\xrightarrow{\ \sim\ }\mathscr A(\mathsf s^\qc_\core)=\mathcal A_{n,m},\qquad M_{\mathbf j}({\bf k}\circ\psi^{-1})\mapsto M^\qc({\bf k})\quad({\bf k}\in\Z^{\Vc}),
\end{align}
and it maps the cluster variables, cluster monomials and frozen variables of every mutated seed of $\mathsf s(\mathbf j)$ to those of the corresponding mutated seed of $\mathsf s^\qc_\core$. It restricts to an isomorphism $\mathscr A^+(\mathsf s(\mathbf j))\cong\mathscr A^+(\mathsf s^\qc_\core)$.
\end{corollary}

\begin{proof}
By Theorem~\ref{thm:seedmatch}, after the relabelling $\psi$, the seeds $\mathsf s(\mathbf j)$ and $\mathsf s^\qc_\core$ have the same vertex set, the same mutable vertices, the same skew-symmetric matrix and agreeing exchange matrices. Hence Lemma~\ref{lem:seediso} gives a unique isomorphism $\mathscr A(\mathsf s(\mathbf j))\to\mathscr A(\mathsf s^\qc_\core)$ with $M_{\mathbf j}({\bf k}\circ\psi^{-1})\mapsto M^\qc({\bf k})$, which matches the cluster variables, cluster monomials and frozen variables of corresponding mutated seeds and restricts to $\mathscr A^+(\mathsf s(\mathbf j))\cong\mathscr A^+(\mathsf s^\qc_\core)$. Finally, by Proposition~\ref{prop:skewLaurent} we have $\mathscr A(\mathsf s^\qc_\core)=\mathcal A_{n,m}$.
\end{proof}

\subsection{Localization of ${\mathcal C}^{\mathrm s}(\Delta^{m-2})$ at the frozen simple objects and categorification of $\mathcal A_{n,m}$}\label{subsec-loc-Delta}

We now apply the framework of Section~\ref{sec-KKOP} to the case $\mathsf b=\Delta^{m-2}$ with the word $\mathbf j=u^{m-2}$. Let
\begin{align}\label{eq-C-Delta}
\cC=\cC^{\mathbb D_{\cQ}}(\Delta^{m-2}),\qquad \widetilde\cC=\widetilde\cC(\Delta^{m-2})=\cC^{\mathrm s}\bigl[F_1^{\otimes-1},\dots,F_{n-1}^{\otimes-1}\bigr]
\end{align}
denote the braid-word category defined in \S\ref{subsec_braid_cat} and its monoidal localization from Proposition~\ref{prop:loc}, respectively. Here, $F_1,\dots,F_{n-1}$ are the frozen modules $\mathsf M_s$, where $s$ runs over the last occurrences in $\mathbf j$. Let
\begin{align}\label{eq-Kprime}
\cK'=\cK_t(\widetilde\cC)\otimes_{\Z[t^{\pm1/2}]}R .
\end{align}

The next theorem follows from the isomorphism~\eqref{eq-Ktloc-cluster} of Section~\ref{sec-KKOP}, for $\mathbf i=\mathbf j$, and Corollary~\ref{cor:seedmatch}. It identifies the algebra $\mathcal A_{n,m}$ of balanced webs with the quantum Grothendieck ring of $\widetilde\cC$, with no twist, because $\mathsf s^\qc_\core$ and $\mathsf s(\mathbf j)$ have the same quantum form. It also shows that $\cK'$ is graded by the cokernel $\Z^{\Lst}$ of the exchange matrix (Lemma~\ref{lem:Lst}), and that the $(q,t)$-characters of all simple objects are homogeneous for this grading. This grading is used to define the twist in \S\ref{subsec-CI}, and the theorem gives Theorem~\ref{thmD}(c).

\begin{theorem}\label{thm:locDelta}
There is an isomorphism of $R$-algebras
\begin{align}\label{eq-PhiA-loc}
\Phi_{\mathcal A}\colon\cK'=\cK_t(\cC)\bigl[[F_1]_t^{-1},\dots,[F_{n-1}]_t^{-1}\bigr]\otimes_{\Z[t^{\pm1/2}]}R\xrightarrow{\ \sim\ }\mathcal A_{n,m},\qquad [\mathsf M_{\psi(x)}]_t\mapsto A^\qc_x\quad(x\in\Vc),
\end{align}
with no twist, which maps the cluster monomials of the mutated seeds of $\mathsf s(\mathbf j)$ to the corresponding cluster monomials of $\mathsf s^\qc_\core$. These cluster monomials are, up to powers of $t^{1/2}$, the $(q,t)$-characters of real simple objects of $\widetilde\cC$ (of the form $\Upsilon(M)\otimes\widetilde F^{\boldsymbol\lambda}$, with $M$ a real simple module of $\cC$; Theorem~\ref{thm:KKOP}(c) and Lemma~\ref{lem:qtloc}). Moreover:
\begin{enumerate}[label={\rm (\alph*)}]
\item $\cK'=\bigoplus_{\gamma\in\Z^{\Lst}}\cK'_\gamma$ is a $\Z^{\Lst}$-graded algebra, where $\cK'_\gamma$ is the set of elements whose expansion in the based quantum torus of $\mathsf s(\mathbf j)$ only involves monomials $M_{\mathbf j}({\bf k}\circ\psi^{-1})$ with $p({\bf k})=\gamma$, for the projection $p$ of~\eqref{eq-p-Lst};
\item the $(q,t)$-character of every simple object of $\widetilde\cC$ is homogeneous for this grading;
\item $\Phi_{\mathcal A}$ maps $\cK'_\gamma$ into the span of the $M^\qc({\bf k})$, ${\bf k}\in\Z^{\Vc}$, with $p({\bf k})=\gamma$.
\end{enumerate}
\end{theorem}

\begin{proof}
The isomorphism is the composite of~\eqref{eq-Ktloc-cluster} for $\mathbf i=\mathbf j$ with Corollary~\ref{cor:seedmatch}, and the statement on cluster monomials follows from Theorem~\ref{thm:KKOP}(c). By Lemma~\ref{lem:Lst}, the grading of~\eqref{eq-Gamma-i} for $\mathbf j$ is the grading by $p$. Let $X\cong\Upsilon(S)\otimes\widetilde F^{\boldsymbol\lambda}$ be simple in $\widetilde\cC$. By Lemma~\ref{lem:qtloc}, $[X]_t=t^c[S]_t[F]_t^{\boldsymbol\lambda}$. Here $[S]_t$ is homogeneous by Proposition~\ref{prop:homog-simple}(b), the $[F_i]_t^{\pm1}$ are monomials of the based quantum torus, and a product of homogeneous elements is homogeneous, with degrees adding, by~\eqref{eq-M-relation}. This proves (b). By Theorem~\ref{thm:Ktloc}(a), the $[X]_t$ form a basis of $\cK'$, so (b) implies that $\cK'=\bigoplus_\gamma\cK'_\gamma$, and $\cK'_\gamma\cK'_{\gamma'}\subseteq\cK'_{\gamma+\gamma'}$; this is (a). Part (c) holds because $\Phi_{\mathcal A}$ is the restriction of the isomorphism of based quantum tori of Lemma~\ref{lem:seediso}.
\end{proof}

\subsection{Graded objects of $\widetilde{\mathcal C}(\Delta^{m-2})$ and monoidal categorification of $\overline{\mathscr S}(\mathbb P_m)$}\label{subsec-CI}

\subsubsection{The category $\widetilde\cC_{\Viso}$}

The category $\widetilde\cC_{\Viso}$ of $\Z^{\Viso}$-graded objects of $\widetilde\cC$ is defined as follows. Its objects are families $X=(X_{\boldsymbol\mu})_{\boldsymbol\mu\in\Z^{\Viso}}$ of objects of $\widetilde\cC$ with $X_{\boldsymbol\mu}=0$ for all but finitely many $\boldsymbol\mu$, morphisms are families of morphisms, composition is componentwise, and
\begin{align}\label{eq-CI-tensor}
(X\otimes Y)_{\boldsymbol\mu}=\bigoplus_{\boldsymbol\alpha+\boldsymbol\beta=\boldsymbol\mu}X_{\boldsymbol\alpha}\otimes Y_{\boldsymbol\beta},
\end{align}
with the associativity constraint induced by that of $\widetilde\cC$ and the unit $\mathbf 1$ in degree $0$. For $X\in\widetilde\cC$ and $\boldsymbol\mu\in\Z^{\Viso}$, write $X\langle\boldsymbol\mu\rangle$ for the object with $X$ in degree $\boldsymbol\mu$ and $0$ elsewhere; then $X\langle\boldsymbol\mu\rangle\otimes Y\langle\boldsymbol\nu\rangle=(X\otimes Y)\langle\boldsymbol\mu+\boldsymbol\nu\rangle$, and every object is the finite direct sum $\bigoplus_{\boldsymbol\mu}X_{\boldsymbol\mu}\langle\boldsymbol\mu\rangle$. Then the category $\widetilde{\mathcal{C}}_{\Viso}$ is an abelian monoidal category with an exact tensor product, since kernels, cokernels, and the tensor product $\otimes$ are all computed componentwise. Its simple objects are of the form $S\langle\boldsymbol{\mu}\rangle$, where $S$ is a simple object in $\widetilde{\mathcal{C}}$. In particular, since the unit object $\mathbf{1}$ is simple in $\widetilde{\mathcal{C}}$, each $\mathbf{1}\langle\boldsymbol{\mu}\rangle$ is simple and invertible, with inverse given by $\mathbf{1}\langle-\boldsymbol{\mu}\rangle$. Denote
\begin{align}\label{eq-Zg}
Z_g=[\mathbf{1}\langle\mathbf{e}_g\rangle]\quad (g\in \Viso), \qquad Z^{\boldsymbol{\mu}}=\prod_g Z_g^{\mu_g}=[\mathbf{1}\langle\boldsymbol{\mu}\rangle].
\end{align}
Then, via the identification $[S\langle\boldsymbol{\mu}\rangle]=[S]Z^{\boldsymbol{\mu}}$, the Grothendieck ring is given by $K(\widetilde{\mathcal{C}}_{\Viso})=K(\widetilde{\mathcal{C}})[Z_g^{\pm1}: g\in \Viso]$.

For a simple object $X\langle\boldsymbol\mu\rangle$ of $\widetilde\cC_{\Viso}$, its \emph{$(q,t)$-character} is $[X\langle\boldsymbol\mu\rangle]_t=[X]_tZ^{\boldsymbol\mu}$ in the Laurent polynomial ring $\cD_t[Z_g^{\pm1}:g\in \Viso]$ with central variables $Z_g$, and the \emph{quantum Grothendieck ring} $\cK_t(\widetilde\cC_{\Viso})$ is the $\Z[t^{\pm1/2}]$-algebra generated by these elements. By Theorem~\ref{thm:Ktloc}, we have
\begin{align}\label{eq-KtCI}
\cK_t(\widetilde\cC_{\Viso})=\cK_t(\widetilde\cC)[Z_g^{\pm1}:g\in \Viso].
\end{align}
Moreover, the $(q,t)$-characters of the simple objects of $\widetilde\cC_{\Viso}$ form a $\Z[t^{\pm1/2}]$-basis with structure constants in $\Z_{\ge0}[t^{\pm1/2}]$, whose values at $t^{1/2}=1$ are composition multiplicities, and its specialization at $t^{1/2}=1$ is $K(\widetilde\cC_{\Viso})$.
 Define the $R$-algebra
\begin{align}\label{eq:Kext}
\cK_{\rm ext}=\cK_t(\widetilde\cC_{\Viso})\otimes_{\Z[t^{\pm1/2}]}R=\cK'[Z_g^{\pm1}:g\in \Viso]=\bigoplus_{\gamma\in\Z^{\Lst},\ \mathbf c\in\Z^{\Viso}}\cK'_\gamma Z^{\mathbf c},
\end{align}
which is graded by $\Gamma:=\Z^{\Lst}\oplus\Z^{\Viso}$, with $\cK'_\gamma Z^{\mathbf c}$ in degree $(\gamma,\mathbf c)$ (Theorem~\ref{thm:locDelta}(a)). By Theorem~\ref{thm:locDelta}(b), the $(q,t)$-character of every simple object of $\widetilde\cC_{\Viso}$ is homogeneous.

\subsubsection{The twist}

We regard $\Gamma=\Z^{\Lst}\oplus\Z^{\Viso}$ as a subgroup of $\Z^V$.

\begin{definition}\label{def:Xi}
Let $\Xi$ be the anti-symmetric integer form on $\Gamma=\Z^{\Lst}\oplus\Z^{\Viso}$ defined by
\begin{align}\label{eq-Xi}
\Xi\bigl((\gamma,\mathbf c),(\gamma',\mathbf c')\bigr)&=\Lambda^\qc(\gamma+\mathbf c,\gamma'+\mathbf c')-\Lambda^\qc(\gamma,\gamma')\nonumber\\&=\Lambda^\qc(\mathbf c,\gamma')-\Lambda^\qc(\mathbf c',\gamma)+\Lambda^\qc(\mathbf c,\mathbf c'),
\end{align}
where $\Lambda^\qc$ is the skew-symmetric matrix of $\mathsf s^\qc_{\rm skein}$ (Definition~\ref{def:sqc}) and we use the notation~\eqref{eq-bilinear}. Thus $\Xi$ vanishes on $\Z^{\Lst}\times\Z^{\Lst}$, and $\Xi({\bf e}_g,{\bf e}_z)=\Lambda^\qc(g,z)$ for $g\in \Viso$ and $z\in\Lst\sqcup \Viso$. For a $\Gamma$-graded $R$-algebra $A$, the \emph{twist} $A^\Xi$ is the graded $R$-module $A$ with the multiplication given by
\begin{align}\label{eq-twist-product}
x\ast y=\xi^{\frac12\Xi(\alpha,\beta)}\,xy\qquad(x\in A_\alpha,\ y\in A_\beta),
\end{align}
\end{definition}

Since $\Xi$ is bilinear, $A^\Xi$ is an associative $R$-algebra.

For $g\in \Viso$ and ${\bf k}\in\Z^{\Vc}$, the value $\Lambda^\qc({\bf e}_g,{\bf k})$ depends only on $p({\bf k})$. Indeed, by~\eqref{eq-compatible}, we have $\Lambda^\qc({\bf e}_g,{\bf q}^\qc_y)=-\Lambda^\qc({\bf q}^\qc_y,{\bf e}_g)=0$ for every mutable vertex $y$. Consequently, $\Lambda^\qc({\bf e}_g,\cdot)$ vanishes on $L_\core$. Hence, for any $\mathbf c\in\Z^{\Viso}$ and ${\bf k}\in\Z^{\Vc}$,

\begin{align}\label{eq-Lambda-depends-p}
\Lambda^\qc(\mathbf c,{\bf k})=\Lambda^\qc(\mathbf c,p({\bf k})).
\end{align}

\subsubsection{The main theorem}

The following is the main theorem of this section; it is Theorem~\ref{thmB} of the introduction. By Proposition~\ref{prop:skewLaurent}, $\overline{\cS}(\bP_m)$ is obtained from $\mathcal A_{n,m}$ by adjoining the frozen variables $A_g$, $g\in\Viso$, and their inverses, and by Theorem~\ref{thm:locDelta}, $\mathcal A_{n,m}$ is the quantum Grothendieck ring of $\widetilde\cC$. On the categorical side, we realize $A_g$ by the invertible object $\mathbf 1\langle\mathbf e_g\rangle$ of $\widetilde\cC_{\Viso}$, whose class $Z_g$ is central. The $A_g$ are not central in $\overline{\cS}(\bP_m)$, but each of them quasi-commutes with every homogeneous element by a power of $\xi$ which only depends on the degree of this element in $\Z^{\Lst}\oplus\Z^{\Viso}$, by~\eqref{eq-Lambda-depends-p}. The twist by $\Xi$ records exactly these powers.

\begin{theorem}\label{thm:main}\label{thm:categorification}
Let $n\ge2$ and $m\ge3$. There is an isomorphism of $R$-algebras
\begin{align}\label{eq-Phi}
\Phi\colon\cK_{\rm ext}^{\ \Xi}=\bigl(\cK_t(\widetilde\cC_{\Viso})\otimes_{\Z[t^{\pm1/2}]}R\bigr)^{\Xi}\xrightarrow{\ \sim\ }\overline{\cS}(\bP_m),\qquad \Phi(aZ^{\mathbf c})=\xi^{\frac12\Lambda^\qc(\mathbf c,\gamma)}\,\Phi_{\mathcal A}(a)\,A^{\mathbf c}
\end{align}
for $a\in\cK'_\gamma$ and $\mathbf c\in\Z^{\Viso}$, where $\Phi_{\mathcal A}$ is the isomorphism of Theorem~\ref{thm:locDelta} and $A^{\mathbf c}$ is~\eqref{eq-Ac}. In particular, $\Phi([\mathsf M_{\psi(x)}]_t)=A^\qc_x$ for $x\in\Vc$, and $\Phi(Z_g)=A_g$ for $g\in \Viso$. Moreover:
\begin{enumerate}[label={\rm (\alph*)}]
\item For every mutation sequence $\mu$ at mutable vertices, $\Phi$ maps each cluster monomial of $\mu(\mathsf s(\mathbf j))$, which is up to a power of $t^{1/2}$ the $(q,t)$-character of a real simple object of $\widetilde\cC$, to the corresponding cluster monomial of $\mu(\mathsf s^\qc_{\rm skein})$, which is the corresponding cluster monomial of $\mu(\mathsf s_{\rm skein})$ times a Laurent monomial in the frozen variables, up to a power of $\xi^{1/2}$. More generally, $\Phi$ maps the $(q,t)$-characters of the real simple objects $(\Upsilon(M)\otimes\widetilde F^{\boldsymbol\lambda})\langle\boldsymbol\mu\rangle$ of $\widetilde\cC_{\Viso}$, with $M$ as above, to these cluster monomials times Laurent monomials in the frozen variables, up to powers of $\xi^{1/2}$.
\item \emph{(Classical limit.)} Setting $\omega^{1/2}=1$, $\Phi$ induces a ring isomorphism
\begin{align}\label{eq-Phi-classical}
K(\widetilde\cC_{\Viso})\xrightarrow{\ \sim\ }\overline\cS_1(\bP_m)=\overline{\cS}(\bP_m)\otimes_R\Z,
\end{align}
where $K(\widetilde\cC_{\Viso})$ is the Grothendieck ring of the abelian monoidal category $\widetilde\cC_{\Viso}$, i.e.\ the free abelian group on the isomorphism classes of its simple objects, with product $[X][Y]=[X\otimes Y]$ (\S\ref{subsec-CI}), under which every cluster monomial of $\overline\cS_1(\bP_m)$ is the class of a real simple object of $\widetilde\cC_{\Viso}$, and every frozen variable is the class of an invertible object.
\end{enumerate}
\end{theorem}

\begin{proof}
By~\eqref{eq:Kext}, $\Phi$ is a well-defined $R$-linear map. According to Theorem~\ref{thm:locDelta}(c), for any $a\in\cK'_\gamma$, $\Phi_{\mathcal A}(a)$ is a linear combination of elements $M^\qc({\bf k})$ with ${\bf k}\in\Z^{\Vc}$ and $p({\bf k})=\gamma$. In view of~\eqref{eq-M-relation} and~\eqref{eq-Lambda-depends-p}, we have $A^{\mathbf c}M^\qc({\bf k})=\xi^{\Lambda^\qc(\mathbf c,\gamma)}M^\qc({\bf k})A^{\mathbf c}$ for such ${\bf k}$, which implies
\begin{align}\label{eq-commute-Ac}
A^{\mathbf c}\,\Phi_{\mathcal A}(a)=\xi^{\Lambda^\qc(\mathbf c,\gamma)}\,\Phi_{\mathcal A}(a)\,A^{\mathbf c}\qquad(a\in\cK'_\gamma).
\end{align}
Let $a\in\cK'_\gamma$ and $b\in\cK'_{\gamma'}$, and let $\mathbf c,\mathbf c'\in\Z^{\Viso}$. We compute
\begin{align*}
\Phi(aZ^{\mathbf c})\,\Phi(bZ^{\mathbf c'})&=\xi^{\frac12\Lambda^\qc(\mathbf c,\gamma)+\frac12\Lambda^\qc(\mathbf c',\gamma')}\,\Phi_{\mathcal A}(a)\,A^{\mathbf c}\,\Phi_{\mathcal A}(b)\,A^{\mathbf c'}\\
&=\xi^{\frac12\Lambda^\qc(\mathbf c,\gamma)+\frac12\Lambda^\qc(\mathbf c',\gamma')+\Lambda^\qc(\mathbf c,\gamma')+\frac12\Lambda^\qc(\mathbf c,\mathbf c')}\,\Phi_{\mathcal A}(ab)\,A^{\mathbf c+\mathbf c'}\\
&=\xi^{\frac12\Lambda^\qc(\mathbf c,\gamma')-\frac12\Lambda^\qc(\mathbf c',\gamma)+\frac12\Lambda^\qc(\mathbf c,\mathbf c')}\,\Phi\bigl(ab\,Z^{\mathbf c+\mathbf c'}\bigr)\\
&=\Phi\bigl(aZ^{\mathbf c}\ast bZ^{\mathbf c'}\bigr).
\end{align*}
Here the first equality is the definition~\eqref{eq-Phi}. The second equality follows from~\eqref{eq-commute-Ac} for $b\in\cK'_{\gamma'}$, which gives $A^{\mathbf c}\Phi_{\mathcal A}(b)=\xi^{\Lambda^\qc(\mathbf c,\gamma')}\Phi_{\mathcal A}(b)A^{\mathbf c}$, from the relation $A^{\mathbf c}A^{\mathbf c'}=\xi^{\frac12\Lambda^\qc(\mathbf c,\mathbf c')}A^{\mathbf c+\mathbf c'}$, which is~\eqref{eq-M-relation}, and from $\Phi_{\mathcal A}(a)\Phi_{\mathcal A}(b)=\Phi_{\mathcal A}(ab)$. The third equality follows from~\eqref{eq-Phi} applied to $ab\in\cK'_{\gamma+\gamma'}$ (Theorem~\ref{thm:locDelta}(a)), i.e.\ $\Phi(abZ^{\mathbf c+\mathbf c'})=\xi^{\frac12\Lambda^\qc(\mathbf c+\mathbf c',\gamma+\gamma')}\Phi_{\mathcal A}(ab)A^{\mathbf c+\mathbf c'}$, and from the bilinearity of $\Lambda^\qc$. The fourth equality follows from~\eqref{eq-twist-product} and~\eqref{eq-Xi}: the elements $aZ^{\mathbf c}$ and $bZ^{\mathbf c'}$ have degrees $(\gamma,\mathbf c)$ and $(\gamma',\mathbf c')$, the product $aZ^{\mathbf c}\cdot bZ^{\mathbf c'}$ in $\cK_{\rm ext}$ equals $ab\,Z^{\mathbf c+\mathbf c'}$ because the $Z_g$ are central, and $\Xi((\gamma,\mathbf c),(\gamma',\mathbf c'))=\Lambda^\qc(\mathbf c,\gamma')-\Lambda^\qc(\mathbf c',\gamma)+\Lambda^\qc(\mathbf c,\mathbf c')$. Thus, $\Phi$ is an algebra homomorphism.

To show that $\Phi$ is bijective, observe that by~\eqref{eq:Kext}, $\cK_{\rm ext}=\bigoplus_{\mathbf c}\cK'Z^{\mathbf c}$. On the summand $\cK'Z^{\mathbf c}$, the map $\Phi$ acts as $aZ^{\mathbf c}\mapsto\sum_\gamma\xi^{\frac12\Lambda^\qc(\mathbf c,\gamma)}\Phi_{\mathcal A}(a_\gamma)A^{\mathbf c}$, where $a=\sum_\gamma a_\gamma$ with $a_\gamma\in\cK'_\gamma$. Since $\Phi_{\mathcal A}$ is bijective and respects the $\Z^{\Lst}$-gradings (Theorem~\ref{thm:locDelta}(c)), this map induces a bijection onto $\mathcal A_{n,m}A^{\mathbf c}$, noting that right multiplication by $A^{\mathbf c}$ is injective. Furthermore, Proposition~\ref{prop:skewLaurent} gives $\overline{\cS}(\bP_m)=\bigoplus_{\mathbf c}\mathcal A_{n,m}A^{\mathbf c}$. The values on generators follow directly from~\eqref{eq-Phi}, since $[\mathsf M_{\psi(x)}]_t$ and $1$ have degrees $p({\bf e}_x)$ and $0$, respectively, and $\Lambda^\qc({\bf e}_g,0)=0$; specifically, for $\Phi([\mathsf M_{\psi(x)}]_t)$, we take $\mathbf c=0$.

(a) By Theorem~\ref{thm:locDelta}, $\Phi_{\mathcal A}$ maps the $(q,t)$-characters of these real simple modules to the cluster monomials of the mutated seeds of $\mathsf s^\qc_\core$. By Proposition~\ref{prop:sqc}(c), these coincide with the cluster monomials of the mutated seeds of $\mathsf s^\qc_{\rm skein}$, and Corollary~\ref{cor:Pchange} relates them to those of $\mathsf s_{\rm skein}$. Regarding the final statement, $[(\Upsilon(M)\otimes\widetilde F^{\boldsymbol\lambda})\langle\boldsymbol\mu\rangle]_t$ is equal to a power of $t^{1/2}$ times $[M]_t[F]_t^{\boldsymbol\lambda}Z^{\boldsymbol\mu}$ (Lemma~\ref{lem:qtloc}). Since $\Phi$ maps $[F_i]_t$ and $Z_g$ to the frozen variables of $\mathsf s^\qc_{\rm skein}$, and the twist modifies products of homogeneous elements only by powers of $\xi^{1/2}$, the assertion follows. These objects are real simple by Proposition~\ref{prop:loc}(c).

(b) We apply the functor $-\otimes_R\Z$, with $\omega^{1/2}\mapsto1$, to $\Phi$. On the right-hand side, this yields $\overline\cS_1(\bP_m)$. On the left-hand side, the twist becomes trivial because $\xi^{\frac12\Xi}\mapsto1$, and the restriction of the map $\omega^{1/2}\mapsto1$ to $\Z[t^{\pm1/2}]$ sends $t^{1/2}$ to $1$. Consequently, $\cK_{\rm ext}\otimes_R\Z=\cK_t(\widetilde\cC_{\Viso})\otimes_{\Z[t^{\pm1/2}]}\Z$, which is identified with $K(\widetilde\cC_{\Viso})$ by Theorem~\ref{thm:Ktloc}(d) and~\eqref{eq-KtCI}. The statement concerning cluster monomials follows from part (a) and Lemma~\ref{lem:specialize}, applied to $\mathsf s_{\rm skein}$ and $\mathsf s^\qc_{\rm skein}$ (via Theorem~\ref{thm:CHW} and Proposition~\ref{prop:sqc}(b)), which identifies the cluster monomials of $\overline\cS_1(\bP_m)$ with the specializations of the quantum ones. By (a), such a specialization corresponds to the class of an object $(\Upsilon(M)\otimes\widetilde F^{\boldsymbol\lambda})\langle\boldsymbol\mu\rangle$ with $M$ real simple. By Proposition~\ref{prop:loc}(c), this object is itself real simple; hence, no invertible factor needs to be split off. Finally, the frozen variables are the images of the classes of the invertible objects $\widetilde F^{\pm\mathbf e_i}\langle0\rangle$ and $\mathbf 1\langle\pm\mathbf e_g\rangle$.
\end{proof}

The following consequence is Corollary~\ref{thmC} of the introduction. Its restriction to the $\Viso$-degree $0$ gives the positive basis of $\mathcal A_{n,m}$ (Theorem~\ref{thm:Aquant}(c)).

\begin{corollary}[Positivity]\label{cor:positivity}
The images under $\Phi$ of the $(q,t)$-characters $[X]_tZ^{\boldsymbol\mu}$ of the simple objects $X\langle\boldsymbol\mu\rangle$ of $\widetilde\cC_{\Viso}$ form an $R$-basis of $\overline{\cS}(\bP_m)$ whose structure constants lie in $\Z_{\ge0}[\xi^{\pm1/2}]$.
\end{corollary}

\begin{proof}
By Theorem~\ref{thm:Ktloc}(a),(b) and~\eqref{eq-KtCI}, these elements form a basis of $\cK_{\rm ext}$ with structure constants in $\Z_{\ge0}[t^{\pm1/2}]\subseteq\Z_{\ge0}[\xi^{\pm1/2}]$. Every basis element is homogeneous (Theorem~\ref{thm:locDelta}(b)), so by~\eqref{eq-twist-product} the twist multiplies each structure constant by a power of $\xi^{1/2}$. Then the result follows by applying $\Phi$.
\end{proof}

\begin{remark}[The quantum parameters]\label{rem:parameters}
In Theorem~\ref{thm:main}, we regard $\Z[t^{\pm1/2}]$ as a subring of $R$ via $t^{1/2}=\xi^{n/2}$; that is, we identify the parameter $t$ of the $(q,t)$-characters with the skein parameter $t=\omega^{n^2}$ defined in~\eqref{eq-constants}. We choose this identification because it ensures the compatibility of the two quantum forms. Specifically, the KKOP seed is compatible with the scalar $2$ for $\Lambda_{\mathbf j}$ in the variable $t^{1/2}$, while the skein seed is compatible with the scalar $2n$ for $\Lambda_{\rm skein}$ in the variable $\xi^{1/2}=t^{1/(2n)}$. Consequently, setting $\Pi_{\mathbf j}=n\Lambda_{\mathbf j}$ allows both seeds to be written in the same variable $\xi$ and to be compatible with the same scalar $2n$, which is precisely what enables the comparison in Theorem~\ref{thm:seedmatch}.
\end{remark}

\section{A quantized braid variety subalgebra of $\overline{\mathscr S}(\mathbb P_m)$}\label{sec-braid}

Casals, Gorsky, Gorsky, Le, Shen, and Simental demonstrated that the coordinate ring of a braid variety admits a cluster algebra structure~\cite{CGGLSS}. Independently, Galashin, Lam, Sherman-Bennett, and Speyer constructed a cluster structure on braid varieties via $3$D plabic graphs~\cite{GLSBS}. In this section, we prove that the subalgebra $\mathcal A_{n,m}$ of balanced stated $n$-webs (Definition~\ref{def:Anm}) is a quantization of the braid variety $X(\Delta^{m-1})$ (Theorem~\ref{thm:Aquant} and Corollary~\ref{cor:Aquantization}), and that the classical reduced skein algebra is isomorphic to the coordinate ring of $X(\Delta^{m-1})$ tensored with a torus (Corollary~\ref{cor:braidvar}). We begin by reviewing braid varieties and their cluster structure in \S\ref{subsec:braidvar}, prove the main theorem in \S\ref{subsec-Aquant}, and derive the consequences in \S\ref{subsec-conseq}. The proofs rely on Proposition~\ref{prop:skewLaurent}, Theorem~\ref{thm:seedmatch}, and Theorem~\ref{thm:locDelta}.

Throughout this section, we keep the notation of Section~\ref{sec-cat}, and abbreviate $\overline{\cS}(\bP_m)$ to $\overline{\cS}$.

\subsection{Braid varieties and their cluster structure}\label{subsec:braidvar}

Let $G=\mathrm{SL}_n(\C)$, and let $B\subseteq G$ be the Borel subgroup of upper triangular matrices. Flags in $\C^n$ are identified with points of the quotient $G/B$. By the Bruhat decomposition $G=\bigsqcup_{w\in S_n}BwB$~\cite[\S28.3]{Hum75}, two flags $gB$ and $g'B$ are \emph{in relative position $w$} if $g^{-1}g'\in BwB$; for a simple transposition $s_i$, this signifies that the two flags differ precisely in their $i$-dimensional subspace. The \emph{Demazure product} on $S_n$ is the associative product defined by $w\ast s_i=ws_i$ if $\ell(ws_i)>\ell(w)$ and $w\ast s_i=w$ otherwise. For a word $\beta=(i_1,\dots,i_l)$, we set $\delta(\beta)=s_{i_1}\ast\cdots\ast s_{i_l}$~\cite[\S3]{KM04}. The \emph{braid variety} associated with $\beta$ is defined as
\begin{align}\label{eq-braidvar}
X(\beta)=\bigl\{(B_1,\dots,B_{l+1})\mid B_1=B,\ B_{l+1}=\delta(\beta)B,\ B_j\text{ and }B_{j+1}\text{ are in relative position }s_{i_j}\bigr\},
\end{align}
see~\cite[\S3]{CGGLSS}. Up to isomorphism, $X(\beta)$ depends only on the positive braid represented by $\beta$. When $\delta(\beta)=w_0$, the variety $X(\beta)$ is smooth and affine, with dimension $l-\binom n2$.

Casals, Gorsky, Gorsky, Le, Shen, and Simental proved that $\C[X(\beta)]$ is a cluster algebra with invertible frozen variables, which coincides with its upper cluster algebra; moreover, its seeds are constructed via Demazure weaves~\cite[Thm.~1.1]{CGGLSS}. Galashin, Lam, Sherman-Bennett, and Speyer subsequently provided a cluster structure using $3$D plabic graphs~\cite{GLSBS}, and it was shown that these two structures agree~\cite{CGGSSS}.

Let $\eta$ be a positive braid word in which every colour appears, and let $\Delta\eta$ denote the word obtained by concatenating a reduced word for $w_0$ with $\eta$. Then $\delta(\Delta\eta)=w_0$. Among the seeds of $\C[X(\Delta\eta)]$ arising from Demazure weaves, there exists a distinguished seed given by the right inductive weave of $\Delta\eta$~\cite[\S4]{CGGLSS}, which we refer to as the \emph{CGGLSS seed attached to $\Delta\eta$}. By~\cite[Prop.~4.44, Cor.~4.45]{CGGLSS}, the quiver $Q^{\mathrm w}(\eta)$ of this seed has the positions of $\eta$ as its vertices, with the last occurrence of each colour being frozen. Furthermore, there is an arrow $s\to s^+$ for every $s$ such that $s^+$ is finite, and for adjacent colours, there is an arrow $s'\to s$ whenever $s<s'<s^+<s'^+$. Comparing this with Definition~\ref{def:BFZ}, we see that $Q^{\mathrm w}(\eta)$ is the quiver $Q(\eta)$ with all arrows reversed:
\begin{align}\label{eq-weave-quiver}
Q^{\mathrm w}(\eta)(s,s')=-Q_\eta(s,s')\qquad\text{for all positions }s,s'.
\end{align}
The matrix $Q_\eta$ is given in~\cite{Bi}; it determines $Q^{\mathrm w}(\eta)$ because the convention in~\cite{Bi} interprets a positive entry $b_{ss'}$ as an arrow $s'\to s$.

\begin{lemma}\label{lem:braidseed}
For $x\in\Vc$ and $y\in V_{\rm mut}$, $Q^\qc(x,y)=-Q^{\mathrm w}(u^{m-2})(\psi(x),\psi(y))$. That is, up to the relabelling $\psi$ of Lemma~\ref{lem:M} and a global sign, the exchange matrix of $\mathsf s^\qc_\core$ is the exchange matrix of the CGGLSS seed of $X(\Delta^{m-1})$ attached to the word $\Delta\cdot u^{m-2}$.
\end{lemma}

\begin{proof}
By Lemma~\ref{lem:M}, $Q^\qc(x,y)=Q_{\mathbf j}(\psi(x),\psi(y))$, and by~\eqref{eq-weave-quiver} for $\eta=u^{m-2}=\mathbf j$, $Q_{\mathbf j}(\psi(x),\psi(y))=-Q^{\mathrm w}(u^{m-2})(\psi(x),\psi(y))$; this proves the first statement. Since $u$ is a reduced word of $w_0$, $\Delta\cdot u^{m-2}$ is a word for $\Delta^{m-1}$, so $X(\Delta\cdot u^{m-2})\cong X(\Delta^{m-1})$.
\end{proof}

A global sign does not change the classical cluster algebra, since the exchange relations are symmetric under $Q\mapsto-Q$. For an exchange matrix $Q$ on a vertex set $\mathcal V$, we write $\mathscr A_1(Q)$ for the classical cluster algebra over $\Z$, with invertible frozen variables, of a classical seed with exchange matrix $Q$ and an algebraically independent initial cluster; up to isomorphism, it depends only on $Q$. Its extension of scalars $\mathscr A_1(Q)\otimes_\Z\C$ is the classical cluster algebra over $\C$ of the same seed, since $\C$ is flat over $\Z$ (Lemma~\ref{lem:specialize}).

\subsection{The balanced subalgebra}\label{subsec-Aquant}

\begin{remark}[Torus invariants]\label{rem:TI}
Let $\mathbf k'$ be an infinite field and $R\to\mathbf k'$ a ring homomorphism. The torus $H^{m-1}=\operatorname{Hom}(\mathsf P^{m-1},\mathbf k'^\times)\cong(\mathbf k'^\times)^{(m-1)(n-1)}$ acts on $\overline{\cS}\otimes_R\mathbf k'$ by algebra automorphisms: by~\eqref{eq-coarse-grading}, $\lambda\in H^{m-1}$ acts on $\overline{\cS}[\nu]\otimes\mathbf k'$ by the scalar $\lambda(\nu)$, i.e.\ it rescales the states of stated webs on the edges $(1,2),\dots,(m-1,m)$. Since $\mathsf P^{m-1}$ is free and $\mathbf k'$ is infinite, $(\overline{\cS}\otimes\mathbf k')^{H^{m-1}}=\overline{\cS}[0]\otimes\mathbf k'=\mathcal A_{n,m}\otimes\mathbf k'$. Classically, $H^{m-1}$ is the torus of decorations at these edges, and this matches the torus $H^{m-1}$ in the description $X(\Delta^{m-1})\times H^{m-1}$ of decorated local systems by Ishibashi--Oya--Shen~\cite[\S6.2.2]{IOS}.
\end{remark}

By Theorem~\ref{thm:locDelta}(a),(c), $\mathcal A_{n,m}=\Phi_{\mathcal A}(\cK')$ is $\Z^{\Lst}$-graded, with $\mathcal A_{n,m,\gamma}:=\Phi_{\mathcal A}(\cK'_\gamma)$ spanned, inside the based quantum torus of $\mathsf s^\qc_\core$, by elements whose expansion only involves the $M^\qc({\bf k})$ with $p({\bf k})=\gamma$. The following theorem is the main result of this section, and it is Theorem~\ref{thmD} of the introduction. Part (a) is the skew Laurent decomposition of Proposition~\ref{prop:skewLaurent}; part (b) shows that $\mathcal A_{n,m}$ is a quantum cluster algebra, equal to its own upper cluster algebra, whose seed is the KKOP seed; part (c) is the categorification of Theorem~\ref{thm:locDelta}; and part (d) identifies the classical limit with the coordinate ring of the braid variety $X(\Delta^{m-1})$, using Lemmas~\ref{lem:specialize} and~\ref{lem:braidseed}.

In summary, we have the following result.

\begin{theorem}\label{thm:Aquant}
Let $n\ge2$ and $m\ge3$.
\begin{enumerate}[label={\rm (\alph*)}]
\item We have
\begin{align}\label{eq-Aquant-a}
\overline{\cS}=\bigoplus_{\mathbf c\in\Z^{\Viso}}\mathcal A_{n,m}\,A^{\mathbf c},
\end{align}
which is a skew Laurent extension of $\mathcal A_{n,m}$ in the frozen variables $A_g$, $g\in \Viso$.
\item $\mathcal A_{n,m}=\mathscr A(\mathsf s^\qc_\core)=\mathscr U(\mathsf s^\qc_\core)$. In particular, $\mathcal A_{n,m}$ is a quantum cluster algebra whose seed coincides with the KKOP seed $\mathsf s(u^{m-2})$.
\item There is an isomorphism of $R$-algebras $\Phi_{\mathcal A}\colon\cK_t(\widetilde\cC)\otimes_{\Z[t^{\pm1/2}]}R\xrightarrow{\sim}\mathcal A_{n,m}$, with no twist, which matches cluster monomials. In particular, $\mathcal A_{n,m}$ is a free $R$-module, and the simple objects of $\widetilde\cC$ give a basis of $\mathcal A_{n,m}$ with structure constants in $\Z_{\ge0}[t^{\pm1/2}]$.
\item At $\omega^{1/2}=1$, there is an isomorphism
\begin{align}\label{eq-Aquant-d}
\mathcal A_{n,m}\otimes_R\C\ \cong\ \C[X(\Delta^{m-1})],
\end{align}
compatible with the cluster structure of~\cite{CGGLSS}.
\end{enumerate}
\end{theorem}

\begin{proof}
(a) This is Proposition~\ref{prop:skewLaurent}.

(b) The first equality follows directly from Proposition~\ref{prop:skewLaurent}. For the upper cluster algebra, by Theorem~\ref{thm-inclusion-quantum} we have $\mathscr A(\mathsf s^\qc_\core)\subseteq\mathscr U(\mathsf s^\qc_\core)$. Conversely, by Proposition~\ref{prop:sqc}(a),(c), each based quantum torus of a mutated seed of $\mathsf s^\qc_{\rm skein}$ is the corresponding based quantum torus of $\mathsf s^\qc_\core$ extended by the Laurent monomials $A^{\mathbf c}$, so $\mathscr U(\mathsf s^\qc_\core)\subseteq\mathscr U(\mathsf s^\qc_{\rm skein})=\overline{\cS}$, using Proposition~\ref{prop:sqc}(b). Let $x\in\mathscr U(\mathsf s^\qc_\core)$, and write $x=\sum_{\bf k}c_{\bf k}M^\qc({\bf k})$ with ${\bf k}\in\Z^{\Vc}$. Choose ${\bf d}\in\Z^{\Vc}_{\ge0}$ with ${\bf d}+{\bf k}\ge0$ whenever $c_{\bf k}\ne0$. Then $M^\qc({\bf d})$ and $M^\qc({\bf d})x$ are linear combinations of products of the $A^\qc_y$, $y\in\Vc$, so they lie in $\overline{\cS}[0]$ by Proposition~\ref{prop:sqc}(d). Writing $x=\sum_\nu x_\nu$ with $x_\nu\in\overline{\cS}[\nu]$, we have $M^\qc({\bf d})x_\nu\in\overline{\cS}[\nu]$, and it is non-zero if $x_\nu\ne0$, since $\overline{\cS}$ is a domain. Hence $x\in\overline{\cS}[0]=\mathcal A_{n,m}$. The statement on seeds is Theorem~\ref{thm:seedmatch}.

(c) The isomorphism follows directly from Theorem~\ref{thm:locDelta}, and it is the restriction of $\Phi$ by~\eqref{eq-Phi} with $\mathbf c=0$.
Since $R$ is free over $\Z[t^{\pm1/2}]$, $\mathcal A_{n,m}$ is a free $R$-module and the images of the $(q,t)$-characters of the simple objects of $\widetilde\cC$ form an $R$-basis of $\mathcal A_{n,m}$ with structure constants in $\Z_{\ge0}[t^{\pm1/2}]$ (Theorem~\ref{thm:Ktloc}(a),(b)).

(d) By (b), Lemma~\ref{lem:specialize} applies to $\mathsf s^\qc_\core$. Taking $\mathbf k'=\C$, it shows that $\mathcal A_{n,m}\otimes_R\C\cong\mathscr A_1(\mathsf s^\qc_\core)\otimes_\Z\C$ is the classical cluster algebra over $\C$, with invertible frozen variables, of the classical limit of $\mathsf s^\qc_\core$, and that each quantum cluster variable specializes to the corresponding classical cluster variable. The exchange matrix of this classical seed is $Q^\qc$ restricted to $\Vc$. By Lemma~\ref{lem:braidseed}, under the relabelling $\psi$ it is the negative of the exchange matrix of the CGGLSS seed of $X(\Delta^{m-1})$ attached to $\Delta\cdot u^{m-2}$; in particular, both seeds have the last occurrences as frozen vertices, and they have the same arrows into them, including the type (Q2) arrows. Since the classical cluster algebras of $Q$ and $-Q$ coincide, and $\C[X(\Delta^{m-1})]$ is the cluster algebra, with invertible frozen variables, of the CGGLSS seed~\cite[Thm.~1.1, Prop.~4.44, Cor.~4.45]{CGGLSS}, we obtain~\eqref{eq-Aquant-d}, compatibly with the cluster structures.
\end{proof}

We record further properties of $\mathcal A_{n,m}$, which were used in the introduction.

\begin{proposition}\label{prop:Aquant-more}
\leavevmode
\begin{enumerate}[label={\rm (\alph*)}]
\item $\mathcal A_{n,m}$ is stable under the reflection $\ast$, and $A_g\,a=\xi^{\Lambda^\qc(g,\gamma)}a\,A_g$ for $a\in\mathcal A_{n,m,\gamma}$ and $g\in \Viso$.
\item The natural map $\mathcal A_{n,m}\otimes_R\Z\to\overline\cS_1(\bP_m)$, $\omega^{1/2}\mapsto1$, is injective, its image is the component $\overline{\cS}[0]\otimes_R\Z$ of degree zero for the grading~\eqref{eq-coarse-grading}, and $\overline\cS_1(\bP_m)=(\mathcal A_{n,m}\otimes_R\Z)[A_g^{\pm1}:g\in \Viso]$ is a Laurent polynomial ring over $\mathcal A_{n,m}\otimes_R\Z$. Moreover, $\mathcal A_{n,m}\otimes_R\Z=\mathscr A_1(\mathsf s^\qc_\core)$.
\item $\mathcal A_{n,m}\otimes_R\C$ is the subalgebra of invariants of the torus $H^{m-1}=\operatorname{Hom}(\mathsf P^{m-1},\C^\times)$ of Remark~\ref{rem:TI} acting on $\overline\cS_1(\bP_m)\otimes_\Z\C$.
\end{enumerate}
\end{proposition}

\begin{proof}
(a) The reflection $\ast$ reverses height orders and switches crossings, but it does not change the states or orientations of endpoints. Hence it preserves edge weights, and $\mathcal A_{n,m}$ is stable under $\ast$. The commutation rule is~\eqref{eq-commute-Ac}.

(b) By Theorem~\ref{thm:Aquant}(a) and~\eqref{eq-coarse-grading}, $\mathcal A_{n,m}\otimes_R\Z=\overline{\cS}[0]\otimes_R\Z$ is a direct summand of $\overline{\cS}\otimes_R\Z=\overline\cS_1(\bP_m)$, so the natural map is injective, and $\overline\cS_1(\bP_m)=\bigoplus_{\mathbf c}(\mathcal A_{n,m}\otimes_R\Z)A^{\mathbf c}$. At $\omega^{1/2}=1$, the factor $\xi^{\Lambda^\qc(g,\gamma)}$ in (a) equals $1$, so the $A_g$ commute with each other and with $\mathcal A_{n,m}\otimes_R\Z$. Hence $\overline\cS_1(\bP_m)$ is the Laurent polynomial ring $(\mathcal A_{n,m}\otimes_R\Z)[A_g^{\pm1}:g\in \Viso]$. The last statement follows from Theorem~\ref{thm:Aquant}(b) and Lemma~\ref{lem:specialize}.

(c) Tensoring with $\C$ over $\Z$ and applying Remark~\ref{rem:TI} with $\mathbf k'=\C$, the image of $\mathcal A_{n,m}\otimes_R\C$ in $\overline\cS_1(\bP_m)\otimes_\Z\C=\overline{\cS}\otimes_R\C$ is the subalgebra of $H^{m-1}$-invariants.
\end{proof}

\subsection{Consequences}\label{subsec-conseq}

The following corollary is the classical identification of the reduced skein algebra with a braid variety times a torus; it is Corollary~\ref{thmE} of the introduction.

\begin{corollary}\label{cor:braidvar}\label{thm:braidvar}
There is an isomorphism of $\C$-algebras
\begin{align}\label{eq-braidvar-cor}
\overline\cS_1(\bP_m)\otimes_\Z\C\ \cong\ \C[X(\Delta^{m-1})]\otimes\C\bigl[Z_g^{\pm1}:g\in \Viso\bigr],\qquad |\Viso|=(m-1)(n-1),
\end{align}
which identifies the classical cluster structure of~\cite{CHW}, after the change of frozen coordinates $\mathsf G_{\rm iso}$ of Lemma~\ref{lem:balP}, with the CGGLSS cluster structure of $X(\Delta^{m-1})$ times the torus, the variables $Z_g$ being the frozen variables $A_g$, $g\in \Viso$.
\end{corollary}

\begin{proof}
By Proposition~\ref{prop:Aquant-more}(b), tensored with $\C$, we have $\overline\cS_1(\bP_m)\otimes_\Z\C=(\mathcal A_{n,m}\otimes_R\C)[A_g^{\pm1}:g\in\Viso]$, and by Theorem~\ref{thm:Aquant}(d), $\mathcal A_{n,m}\otimes_R\C\cong\C[X(\Delta^{m-1})]$, compatibly with the cluster structures. By Theorem~\ref{thm:CHW}, Proposition~\ref{prop:sqc}(b) and Lemma~\ref{lem:specialize}, the classical cluster structure of~\cite{CHW}, after the change of frozen coordinates $\mathsf G_{\rm iso}$, is that of the classical limit of $\mathsf s^\qc_{\rm skein}$, whose cluster variables at the vertices of $\Vc$ are those of $\mathsf s^\qc_\core$ (Proposition~\ref{prop:sqc}(c)) and whose remaining frozen variables are the $A_g$, $g\in\Viso$.
\end{proof}

\begin{remark}\label{rem:IOS}
The dimensions of both sides of \eqref{eq-braidvar-cor} match: $\dim X(\Delta^{m-1})+|\Viso|=(m-2)\binom n2+(m-1)(n-1)=|V|$, which coincides with the Gelfand--Kirillov dimension of $\overline{\cS}(\bP_m)$ computed by L\^e--Yu~\cite{LY}. Geometrically, $\operatorname{Spec}(\overline\cS_1(\bP_m)\otimes_\Z\C)$ is expected to be the moduli space of decorated $\mathrm{SL}_n$-local systems on $\bP_m$ subject to reduction conditions at the corners; Ishibashi--Oya--Shen identify the generic part of this moduli space with $X(\Delta^{m-1})\times H^{m-1}$~\cite[\S6.2.2]{IOS}. Corollary~\ref{cor:braidvar} serves as the skein-theoretic counterpart, and its proof is algebraic: the torus factor is generated by the frozen variables $A_g$ for $g\in \Viso$, while the braid-variety factor corresponds to the subalgebra of $H^{m-1}$-invariants (Proposition~\ref{prop:Aquant-more}(c)). Consequently, $H^{m-1}$ acts freely on $\operatorname{Spec}(\overline\cS_1(\bP_m)\otimes_\Z\C)\cong X(\Delta^{m-1})\times H^{m-1}$ by translation on the second factor, with quotient $X(\Delta^{m-1})$.
\end{remark}

We now make precise in what sense $\mathcal A_{n,m}$ is a quantization, following~\cite{Bi}.

\begin{definition}\label{def:quantization}
Let $\mathcal B$ be a commutative $\C$-algebra with a cluster structure, i.e.\ $\mathcal B=\mathscr A_1(Q_1)\otimes_\Z\C$ for a classical seed with exchange matrix $Q_1$. A {\rm cluster quantization} of $\mathcal B$ over $R$ is a quantum cluster algebra $\mathscr A(\mathsf s)$ over $R$, $\mathsf s=(Q,\Pi,M)$, such that $\mathscr A(\mathsf s)$ is a free $R$-module, $Q$ and $Q_1$ agree up to a global sign (Remark~\ref{rem:frozen-entries}), and the specialization $\omega^{1/2}\mapsto1$ induces an isomorphism $\mathscr A(\mathsf s)\otimes_R\C\cong\mathcal B$ which maps quantum cluster variables to the corresponding cluster variables.
\end{definition}

\begin{corollary}\label{cor:Aquantization}
$\mathcal A_{n,m}$ is a cluster quantization of $\C[X(\Delta^{m-1})]$ over $R$. It has a basis with structure constants in $\Z_{\ge0}[t^{\pm1/2}]$, it is stable under the reflection $\ast$, and $\mathcal A_{n,m}\otimes_R\C$ is the subalgebra of $H^{m-1}$-invariants of $\overline\cS_1(\bP_m)\otimes_\Z\C$.
\end{corollary}

\begin{proof}
By Theorem~\ref{thm:Aquant}(b) and (c), $\mathcal A_{n,m}=\mathscr A(\mathsf s^\qc_\core)$ is a quantum cluster algebra which is a free $R$-module. By Lemma~\ref{lem:braidseed}, its exchange matrix agrees, up to a global sign, with that of the CGGLSS seed of $X(\Delta^{m-1})$, and by Theorem~\ref{thm:Aquant}(d), the specialization $\omega^{1/2}\mapsto1$ gives an isomorphism $\mathcal A_{n,m}\otimes_R\C\cong\C[X(\Delta^{m-1})]$ which maps quantum cluster variables to the corresponding cluster variables. Hence $\mathcal A_{n,m}$ is a cluster quantization of $\C[X(\Delta^{m-1})]$ in the sense of Definition~\ref{def:quantization}. The remaining statements are Theorem~\ref{thm:Aquant}(c) and Proposition~\ref{prop:Aquant-more}(a),(c).
\end{proof}

\bibliography{ref}

@misc{Bi,
  author = {Bi, Y.},
  title  = {Monoidal categorification and quantization of braid varieties},
  note   = {arXiv:2609.02962},
  year   = {2026}
}

@article{BZ,
  author  = {Berenstein, A. and Zelevinsky, A.},
  title   = {Quantum cluster algebras},
  journal = {Adv. Math.},
  volume  = {195},
  year    = {2005},
  pages   = {405--455}
}

@article{CGGLSS,
  author  = {Casals, R. and Gorsky, E. and Gorsky, M. and Le, I. and Shen, L. and Simental, J.},
  title   = {Cluster structures on braid varieties},
  journal = {J. Amer. Math. Soc.},
  volume  = {38},
  number  = {2},
  year    = {2025},
  pages   = {369--479}
}

@article{CGGS,
  author  = {Casals, R. and Gorsky, E. and Gorsky, M. and Simental, J.},
  title   = {Algebraic weaves and braid varieties},
  journal = {Amer. J. Math.},
  volume  = {146},
  year    = {2024},
  pages   = {1469--1576},
  note    = {arXiv:2012.06931}
}

@misc{CGGSSS,
  author = {Casals, R. and Galashin, P. and Gorsky, M. and Shen, L. and Sherman-Bennett, M. and Simental, J.},
  title  = {Comparing cluster algebras on braid varieties},
  note   = {arXiv:2508.03816},
  year   = {2025}
}

@misc{CHW,
  author = {Cao, P. and Huang, M. and Wang, Z.},
  title  = {Quantum cluster algebra realization for stated {${\rm SL}_n$}-skein algebras and rotation-invariant bases for polygons},
  note   = {arXiv:2605.12114},
  year   = {2026}
}

@article{CN22,
  author  = {Casals, R. and Ng, L.},
  title   = {Braid loops with infinite monodromy on the {L}egendrian contact {DGA}},
  journal = {J. Topol.},
  volume  = {15},
  year    = {2022},
  pages   = {1927--2016}
}

@article{CP91,
  author  = {Chari, V. and Pressley, A.},
  title   = {Quantum affine algebras},
  journal = {Comm. Math. Phys.},
  volume  = {142},
  year    = {1991},
  pages   = {261--283}
}

@book{CP-book,
  author    = {Chari, V. and Pressley, A.},
  title     = {A guide to quantum groups},
  publisher = {Cambridge University Press},
  year      = {1994}
}

@incollection{CP95,
  author    = {Chari, V. and Pressley, A.},
  title     = {Quantum affine algebras and their representations},
  booktitle = {Representations of groups ({B}anff, 1994)},
  series    = {CMS Conf. Proc.},
  volume    = {16},
  publisher = {Amer. Math. Soc.},
  year      = {1995},
  pages     = {59--78}
}

@incollection{FR,
  author    = {Frenkel, E. and Reshetikhin, N.},
  title     = {The {$q$}-characters of representations of quantum affine algebras and deformations of {$\mathcal W$}-algebras},
  booktitle = {Recent developments in quantum affine algebras and related topics},
  series    = {Contemp. Math.},
  volume    = {248},
  publisher = {Amer. Math. Soc.},
  year      = {1999},
  pages     = {163--205}
}

@article{FO21,
  author  = {Fujita, R. and Oh, S.-j.},
  title   = {Q-data and representation theory of untwisted quantum affine algebras},
  journal = {Comm. Math. Phys.},
  volume  = {384},
  year    = {2021},
  pages   = {1351--1407}
}

@article{Fraser,
  author  = {Fraser, C.},
  title   = {Quasi-homomorphisms of cluster algebras},
  journal = {Adv. in Appl. Math.},
  volume  = {81},
  year    = {2016},
  pages   = {40--77}
}

@misc{GLSBS,
  author = {Galashin, P. and Lam, T. and Sherman-Bennett, M. and Speyer, D.},
  title  = {Braid variety cluster structures, {I}: 3{D} plabic graphs},
  note   = {arXiv:2210.04778},
  year   = {2022}
}

@article{Hernandez,
  author  = {Hernandez, D.},
  title   = {Algebraic approach to {$q,t$}-characters},
  journal = {Adv. Math.},
  volume  = {187},
  year    = {2004},
  pages   = {1--52}
}

@article{HL10,
  author  = {Hernandez, D. and Leclerc, B.},
  title   = {Cluster algebras and quantum affine algebras},
  journal = {Duke Math. J.},
  volume  = {154},
  year    = {2010},
  pages   = {265--341}
}

@article{HL15,
  author  = {Hernandez, D. and Leclerc, B.},
  title   = {Quantum {G}rothendieck rings and derived {H}all algebras},
  journal = {J. Reine Angew. Math.},
  volume  = {701},
  year    = {2015},
  pages   = {77--126}
}

@misc{HW,
  author = {Huang, M. and Wang, Z.},
  title  = {Quantum cluster realization for projected stated {${\rm SL}_n$}-skein algebras},
  note   = {arXiv:2509.25938},
  year   = {2025}
}

@article{IOS,
  author  = {Ishibashi, T. and Oya, H. and Shen, L.},
  title   = {{$\mathscr A=\mathscr U$} for cluster algebras from moduli spaces of {$G$}-local systems},
  journal = {Adv. Math.},
  volume  = {431},
  year    = {2023},
  pages   = {Paper No. 109256}
}

@article{Kal06,
  author  = {K{\'a}lm{\'a}n, T.},
  title   = {Braid-positive {L}egendrian links},
  journal = {Int. Math. Res. Not.},
  volume  = {2006},
  year    = {2006},
  pages   = {Art. ID 14874, 29 pp.}
}

@article{Kas02,
  author  = {Kashiwara, M.},
  title   = {On level-zero representations of quantized affine algebras},
  journal = {Duke Math. J.},
  volume  = {112},
  year    = {2002},
  pages   = {117--175}
}

@article{KKK18,
  author  = {Kang, S.-J. and Kashiwara, M. and Kim, M.},
  title   = {Symmetric quiver {H}ecke algebras and {$R$}-matrices of quantum affine algebras},
  journal = {Invent. Math.},
  volume  = {211},
  year    = {2018},
  pages   = {591--685},
  note    = {arXiv:1304.0323}
}

@article{KKOP-comp,
  author  = {Kashiwara, M. and Kim, M. and Oh, S.-j. and Park, E.},
  title   = {Monoidal categorification and quantum affine algebras},
  journal = {Compos. Math.},
  volume  = {156},
  year    = {2020},
  pages   = {1039--1077},
  note    = {arXiv:1910.08307}
}

@article{KKOP-loc,
  author  = {Kashiwara, M. and Kim, M. and Oh, S.-j. and Park, E.},
  title   = {Localizations for quiver {H}ecke algebras},
  journal = {Pure Appl. Math. Q.},
  volume  = {17},
  year    = {2021},
  pages   = {1465--1548},
  note    = {arXiv:1901.09319}
}

@article{KKOP2,
  author  = {Kashiwara, M. and Kim, M. and Oh, S.-j. and Park, E.},
  title   = {Monoidal categorification and quantum affine algebras {II}},
  journal = {Invent. Math.},
  volume  = {236},
  number  = {2},
  year    = {2024},
  pages   = {837--924}
}

@misc{KKOP3,
  author = {Kashiwara, M. and Kim, M. and Oh, S.-j. and Park, E.},
  title  = {Monoidal categorification and quantum affine algebras {III}},
  note   = {arXiv:2509.14552},
  year   = {2025}
}

@article{LS,
  author  = {L{\^e}, T. T. Q. and Sikora, A. S.},
  title   = {Stated {${\rm SL}(n)$}-skein modules and algebras},
  journal = {J. Topol.},
  volume  = {17},
  number  = {3},
  year    = {2024},
  pages   = {Paper No. e12350}
}

@misc{LY,
  author = {L{\^e}, T. T. Q. and Yu, T.},
  title  = {Quantum traces for {${\rm SL}_n$}-skein algebras},
  note   = {arXiv:2303.08082},
  year   = {2023}
}

@article{Nakajima,
  author  = {Nakajima, H.},
  title   = {Quiver varieties and {$t$}-analogs of {$q$}-characters of quantum affine algebras},
  journal = {Ann. of Math.},
  volume  = {160},
  year    = {2004},
  pages   = {1057--1097}
}

@misc{Qin24a,
  author = {Qin, F.},
  title  = {Analogs of the dual canonical bases for cluster algebras from {L}ie theory},
  note   = {arXiv:2407.02480},
  year   = {2024}
}

@article {CQW26,
    AUTHOR = {Contu, Alessandro and Qin, Fan and Wei, Qiaoling},
     TITLE = {From {$i$}-boxes to signed words},
   JOURNAL = {Pacific J. Math.},
  FJOURNAL = {Pacific Journal of Mathematics},
    VOLUME = {342},
      YEAR = {2026},
    NUMBER = {1},
     PAGES = {63--78},
}

@article {BFZ05,
    AUTHOR = {Berenstein, Arkady and Fomin, Sergey and Zelevinsky, Andrei},
     TITLE = {Cluster algebras. {III}. {U}pper bounds and double {B}ruhat
              cells},
   JOURNAL = {Duke Math. J.},
  FJOURNAL = {Duke Mathematical Journal},
    VOLUME = {126},
      YEAR = {2005},
    NUMBER = {1},
     PAGES = {1--52},
}

@incollection{VV,
  author    = {Varagnolo, M. and Vasserot, E.},
  title     = {Perverse sheaves and quantum {G}rothendieck rings},
  booktitle = {Studies in memory of {I}ssai {S}chur},
  series    = {Progr. Math.},
  volume    = {210},
  publisher = {Birkh{\"a}user},
  year      = {2003},
  pages     = {345--365}
}

@article{Dri85,
  author  = {Drinfeld, V. G.},
  title   = {Hopf algebras and the quantum {Y}ang--{B}axter equation},
  journal = {Soviet Math. Dokl.},
  volume  = {32},
  year    = {1985},
  pages   = {254--258}
}

@article{Jim85,
  author  = {Jimbo, M.},
  title   = {A {$q$}-difference analogue of {$U(\mathfrak g)$} and the {Y}ang--{B}axter equation},
  journal = {Lett. Math. Phys.},
  volume  = {10},
  year    = {1985},
  pages   = {63--69}
}

@inproceedings{Dri87,
  author    = {Drinfeld, V. G.},
  title     = {Quantum groups},
  booktitle = {Proceedings of the {I}nternational {C}ongress of {M}athematicians ({B}erkeley, {C}alif., 1986)},
  volume    = {1},
  publisher = {Amer. Math. Soc., Providence, RI},
  year      = {1987},
  pages     = {798--820}
}

@article{FrRe92,
  author  = {Frenkel, I. B. and Reshetikhin, N. Yu.},
  title   = {Quantum affine algebras and holonomic difference equations},
  journal = {Comm. Math. Phys.},
  volume  = {146},
  year    = {1992},
  pages   = {1--60}
}

@article{KKKO15,
  author  = {Kang, S.-J. and Kashiwara, M. and Kim, M. and Oh, S.-j.},
  title   = {Simplicity of heads and socles of tensor products},
  journal = {Compos. Math.},
  volume  = {151},
  year    = {2015},
  pages   = {377--396}
}

@article{KKKO18,
  author  = {Kang, S.-J. and Kashiwara, M. and Kim, M. and Oh, S.-j.},
  title   = {Monoidal categorification of cluster algebras},
  journal = {J. Amer. Math. Soc.},
  volume  = {31},
  year    = {2018},
  pages   = {349--426}
}

@article{KKOP-PBW,
  author  = {Kashiwara, M. and Kim, M. and Oh, S.-j. and Park, E.},
  title   = {{PBW} theory for quantum affine algebras},
  journal = {J. Eur. Math. Soc. (JEMS)},
  volume  = {26},
  year    = {2024},
  pages   = {2679--2743}
}

@article{KKOP-braid,
  author  = {Kashiwara, M. and Kim, M. and Oh, S.-j. and Park, E.},
  title   = {Braid group action on the module category of quantum affine algebras},
  journal = {Proc. Japan Acad. Ser. A Math. Sci.},
  volume  = {97},
  number  = {3},
  year    = {2021},
  pages   = {13--18}
}

@article{KKOP-root,
  author  = {Kashiwara, M. and Kim, M. and Oh, S.-j. and Park, E.},
  title   = {Simply-laced root systems arising from quantum affine algebras},
  journal = {Compos. Math.},
  volume  = {158},
  number  = {1},
  year    = {2022},
  pages   = {168--210}
}

@article{CM05,
  author  = {Chari, V. and Moura, A.},
  title   = {Characters and blocks for finite-dimensional representations of quantum affine algebras},
  journal = {Int. Math. Res. Not.},
  volume  = {2005},
  number  = {5},
  year    = {2005},
  pages   = {257--298}
}

@article{Qin17,
  author  = {Qin, F.},
  title   = {Triangular bases in quantum cluster algebras and monoidal categorification conjectures},
  journal = {Duke Math. J.},
  volume  = {166},
  year    = {2017},
  pages   = {2337--2442}
}

@article{Bit21,
  author  = {Bittmann, L.},
  title   = {A quantum cluster algebra approach to representations of simply laced quantum affine algebras},
  journal = {Math. Z.},
  volume  = {298},
  year    = {2021},
  pages   = {1449--1485}
}

@article{HL16,
  author  = {Hernandez, D. and Leclerc, B.},
  title   = {A cluster algebra approach to {$q$}-characters of {K}irillov--{R}eshetikhin modules},
  journal = {J. Eur. Math. Soc. (JEMS)},
  volume  = {18},
  year    = {2016},
  pages   = {1113--1159}
}

@article{Le18,
  author  = {L{\^e}, T. T. Q.},
  title   = {Triangular decomposition of skein algebras},
  journal = {Quantum Topol.},
  volume  = {9},
  year    = {2018},
  pages   = {591--632}
}

@article{CL22,
  author  = {Costantino, F. and L{\^e}, T. T. Q.},
  title   = {Stated skein algebras of surfaces},
  journal = {J. Eur. Math. Soc. (JEMS)},
  volume  = {24},
  year    = {2022},
  pages   = {4063--4142}
}

@article{Kup96,
  author  = {Kuperberg, G.},
  title   = {Spiders for rank 2 {L}ie algebras},
  journal = {Comm. Math. Phys.},
  volume  = {180},
  year    = {1996},
  pages   = {109--151}
}

@article{Sik05,
  author  = {Sikora, A. S.},
  title   = {Skein theory for {$SU(n)$}-quantum invariants},
  journal = {Algebr. Geom. Topol.},
  volume  = {5},
  year    = {2005},
  pages   = {865--897}
}

@article{CKM14,
  author  = {Cautis, S. and Kamnitzer, J. and Morrison, S.},
  title   = {Webs and quantum skew {H}owe duality},
  journal = {Math. Ann.},
  volume  = {360},
  year    = {2014},
  pages   = {351--390}
}

@article{FG06,
  author  = {Fock, V. and Goncharov, A.},
  title   = {Moduli spaces of local systems and higher {T}eichm{\"u}ller theory},
  journal = {Publ. Math. Inst. Hautes {\'E}tudes Sci.},
  volume  = {103},
  year    = {2006},
  pages   = {1--211}
}

@article{FZ02,
  author  = {Fomin, S. and Zelevinsky, A.},
  title   = {Cluster algebras {I}: {F}oundations},
  journal = {J. Amer. Math. Soc.},
  volume  = {15},
  year    = {2002},
  pages   = {497--529}
}

@article{Gar69,
  author  = {Garside, F. A.},
  title   = {The braid group and other groups},
  journal = {Quart. J. Math. Oxford Ser. (2)},
  volume  = {20},
  year    = {1969},
  pages   = {235--254}
}

@book{Hum75,
  author    = {Humphreys, James E.},
  title     = {Linear Algebraic Groups},
  series    = {Graduate Texts in Mathematics},
  volume    = {21},
  publisher = {Springer-Verlag, New York-Heidelberg},
  year      = {1975},
}

@article{KM04,
  author  = {Knutson, Allen and Miller, Ezra},
  title   = {Subword complexes in {C}oxeter groups},
  journal = {Adv. Math.},
  volume  = {184},
  number  = {1},
  pages   = {161--176},
  year    = {2004},
}

@article{Thu14,
  author  = {Thurston, D. P.},
  title   = {Positive basis for surface skein algebras},
  journal = {Proc. Natl. Acad. Sci. USA},
  volume  = {111},
  number  = {27},
  pages   = {9725--9732},
  year    = {2014}
}

@article{QW21,
  author  = {Queffelec, H. and Wedrich, P.},
  title   = {Khovanov homology and categorification of skein modules},
  journal = {Quantum Topol.},
  volume  = {12},
  number  = {1},
  pages   = {129--209},
  year    = {2021}
}

@article{Prz91,
  author  = {Przytycki, J. H.},
  title   = {Skein modules of 3-manifolds},
  journal = {Bull. Polish Acad. Sci. Math.},
  volume  = {39},
  number  = {1-2},
  pages   = {91--100},
  year    = {1991}
}

@article{Tur91,
  author  = {Turaev, V. G.},
  title   = {Skein quantization of {P}oisson algebras of loops on surfaces},
  journal = {Ann. Sci. \'Ecole Norm. Sup. (4)},
  volume  = {24},
  number  = {6},
  pages   = {635--704},
  year    = {1991}
}

@article{BW11,
  author  = {Bonahon, F. and Wong, H.},
  title   = {Quantum traces for representations of surface groups in {${\rm SL}_2(\mathbb C)$}},
  journal = {Geom. Topol.},
  volume  = {15},
  number  = {3},
  pages   = {1569--1615},
  year    = {2011}
}

@article{Mul16,
  author  = {Muller, G.},
  title   = {Skein and cluster algebras of marked surfaces},
  journal = {Quantum Topol.},
  volume  = {7},
  number  = {3},
  pages   = {435--503},
  year    = {2016}
}

@article{IY23,
  author  = {Ishibashi, T. and Yuasa, W.},
  title   = {Skein and cluster algebras of unpunctured surfaces for {$\mathfrak{sl}_3$}},
  journal = {Math. Z.},
  volume  = {303},
  number  = {3},
  pages   = {Paper No. 72},
  year    = {2023}
}

@book{Cohn,
  author    = {P. M. Cohn},
  title     = {Skew fields: Theory of general division rings},
  series    = {Encyclopedia of Mathematics and its Applications},
  volume    = {57},
  publisher = {Cambridge University Press},
  address   = {Cambridge},
  year      = {1995}
}

@article{goodearl2017quantum,
  author  = {K. R. Goodearl and M. T. Yakimov},
  title   = {Quantum cluster algebra structures on quantum nilpotent algebras},
  journal = {Mem. Amer. Math. Soc.},
  volume  = {247},
  number  = {1169},
  year    = {2017}
}

@article{CHL,
  author  = {W. Chang and M. Huang and J.-R. Li},
  title   = {Quasi-homomorphisms of quantum cluster algebras},
  journal = {J. Algebra},
  volume  = {638},
  pages   = {506--534},
  year    = {2024}
}

@article{KWQ,
  author  = {Y. Kimura and F. Qin and Q. Wei},
  title   = {Twist automorphisms and {P}oisson structures},
  journal = {SIGMA Symmetry Integrability Geom. Methods Appl.},
  volume  = {19},
  pages   = {Paper No. 105, 39 pp.},
  year    = {2023}
}

@article{KKK15II,
  author  = {S.-J. Kang and M. Kashiwara and M. Kim},
  title   = {Symmetric quiver {H}ecke algebras and {R}-matrices of quantum affine algebras, {II}},
  journal = {Duke Math. J.},
  volume  = {164},
  number  = {8},
  pages   = {1549--1602},
  year    = {2015}
}

@article{Fuj,
  author  = {R. Fujita},
  title   = {Affine highest weight categories and quantum affine {S}chur--{W}eyl duality of {D}ynkin quiver types},
  journal = {Represent. Theory},
  volume  = {26},
  pages   = {211--263},
  year    = {2022}
}

@article{GLS-spec,
  author  = {C. Gei{\ss} and B. Leclerc and J. Schr{\"o}er},
  title   = {Quantum cluster algebras and their specializations},
  journal = {J. Algebra},
  volume  = {558},
  pages   = {411--422},
  year    = {2020}
}

@article{APS04,
  author  = {M. M. Asaeda and J. H. Przytycki and A. S. Sikora},
  title   = {Categorification of the {K}auffman bracket skein module of {$I$}-bundles over surfaces},
  journal = {Algebr. Geom. Topol.},
  volume  = {4},
  pages   = {1177--1210},
  year    = {2004}
}

@misc{AKS25,
  author        = {D. G. L. Allegretti and H. K. Kim and P. Shan},
  title         = {Monoidal categorification of genus zero skein algebras},
  year          = {2025},
  eprint        = {2505.13332},
  archivePrefix = {arXiv},
  note          = {arXiv:2505.13332}
}

@misc{Qin24b,
  author        = {F. Qin},
  title         = {Based cluster algebras of infinite rank},
  year          = {2024},
  eprint        = {2409.02881v4},
  archivePrefix = {arXiv},
  note          = {arXiv:2409.02881v4}
}

\end{document}